\documentclass[11pt]{amsart} 
\usepackage{amscd}           
\usepackage{amssymb}         

\newtheorem{thm}{Theorem}[section]
\newtheorem{lem}[thm]{Lemma}

\newtheorem{cor}[thm]{Corollary}

\newtheorem{prop}[thm]{Proposition}

\theoremstyle{definition}

\renewcommand{\thecase}{}

\newtheorem{conj}[thm]{Conjecture}
\newtheorem{conv}[thm]{Convention}

\newtheorem{defn}[thm]{Definition}

\newtheorem{rmk}[thm]{Remark} 
\renewcommand{\thestep}{}

\theoremstyle{remark}

\makeatletter
\def\alphenumi{
  \def\theenumi{\alph{enumi}}
  \def\p@enumi{\theenumi}
  \def\labelenumi{(\@alph\c@enumi)}}
\makeatother

\makeatletter
\def\thecase{\@arabic\c@case}
\makeatother

\numberwithin{equation}{section}

\makeatletter
\def\thestep{\@arabic\c@step}
\makeatother

\renewcommand\emptyset{\varnothing}

\newcommand\embed{\hookrightarrow}

\newcommand\barM{{\bar{M}}}

\newcommand\barZ{{\bar{Z}}}

\newcommand\ubarRR{{\underline{\mathbb{R}}}}

\newcommand\AAA{\mathbb{A}}
\newcommand\BB{\mathbb{B}}
\newcommand\CC{\mathbb{C}}

\newcommand\FF{\mathbb{F}}

\newcommand\LL{\mathbb{L}}

\newcommand\NN{\mathbb{N}}

\newcommand\PP{\mathbb{P}}

\newcommand\RR{\mathbb{R}}

\newcommand\ZZ{\mathbb{Z}}

\newcommand\bdelta{{\boldsymbol{\delta}}}

\newcommand\bga{{\boldsymbol{\gamma}}}
\newcommand\bgamma{{\boldsymbol{\gamma}}}

\newcommand\bl{{\boldsymbol{\ell}}}

\newcommand\bvarphi{{\boldsymbol{\varphi}}}

\newcommand\bsD{{\boldsymbol{\mathcal{D}}}}

\newcommand\bD{{\mathbf{D}}}

\newcommand\bK{{\mathbf{K}}}

\newcommand\bL{{\mathbf{L}}}

\newcommand\bS{{\mathbf{S}}}

\newcommand\bW{{\mathbf{W}}}
\newcommand\bx{{\mathbf{x}}}

\newcommand\by{{\mathbf{y}}}

\newcommand{\rd}{\partial}

\newcommand\HG{{}_2 F_1}

\newcommand\half{{{\frac{1}{2}}}}

\newcommand\quarter{{{\frac{1}{4}}}}

\newcommand\fg{{\mathfrak{g}}}

\newcommand\fs{{\mathfrak{s}}}

\newcommand\ft{{\mathfrak{t}}}

\newcommand\eps{\varepsilon}
\newcommand\ga{\gamma}

\newcommand\La{\Lambda}
\newcommand\ka{\kappa}
\newcommand\om{\omega}
\newcommand\Om{\Omega}
\newcommand\si{\sigma}

\newcommand\su{{\mathfrak{s}\mathfrak{u}}}

\newcommand\BS{\operatorname{BS}}

\newcommand\ES{\operatorname{ES}}

\newcommand\PU{\operatorname{PU}}

\newcommand\SO{\operatorname{SO}}

\newcommand\SU{\operatorname{SU}}
\newcommand\U{\operatorname{U}}

\newcommand\less{\setminus}

\newcommand{\8}{\infty}

\newcommand\asd{{\operatorname{asd}}}

\newcommand\Aut{\operatorname{Aut}}

\newcommand\Coker{\operatorname{Coker}}

\newcommand\ind{\operatorname{Index}}

\newcommand\Ker{\operatorname{Ker}}

\newcommand\Map{\operatorname{Map}}

\newcommand\PD{\operatorname{PD}}

\newcommand\Ran{\operatorname{Ran}}
\newcommand\rank{\operatorname{rank}}

\newcommand\red{\operatorname{red}}

\newcommand\SW{SW}
\newcommand\Sym{\operatorname{Sym}}
\newcommand\Tor{\operatorname{Tor}}

\newcommand\Ad{{\mathrm{Ad}\,}}

\newcommand\even{{\mathrm{even}}}

\newcommand\id{{\mathrm{id}}}

\newcommand\odd{{\mathrm{odd}}}

\newcommand\spinc{\text{$\text{spin}^c$ }}
\newcommand\spinu{\text{$\text{spin}^u$ }}
\newcommand\Spinc{\text{$\text{Spin}^c$}}

\newcommand\sA{{\mathcal{A}}}
\newcommand\sB{{\mathcal{B}}}
\newcommand\sC{{\mathcal{C}}}
\newcommand\sD{{\mathcal{D}}}

\newcommand\sG{{\mathcal{G}}}

\newcommand\sM{{\mathcal{M}}}

\newcommand\sO{{\mathcal{O}}}

\newcommand\sU{{\mathcal{U}}}
\newcommand\sV{{\mathcal{V}}}
\newcommand\sW{{\mathcal{W}}}
\newcommand\sZ{{\mathcal{Z}}}

\newcommand\tsC{{\tilde\sC}}

\newcommand\tM{{\tilde M}}
\newcommand\tN{{\tilde N}}

\begin{document}
\title[PU(2) Monopoles. II]
{PU(2) Monopoles. II: Top-level Seiberg-Witten Moduli Spaces
and Witten's Conjecture in Low Degrees}
\author[Paul M. N. Feehan]{Paul M. N. Feehan}
\address{Department of Mathematics\\
Ohio State University\\
Columbus, OH 43210}
\email{feehan@math.ohio-state.edu}
\urladdr{http://www.math.ohio-state.edu/$\sim$feehan/} 
\author[Thomas G. Leness]{Thomas G. Leness}
\address{Department of Mathematics\\
Florida International University\\
Miami, FL 33199}
\email{lenesst@fiu.edu}
\urladdr{http://www.fiu.edu/$\sim$lenesst/} 
\date{July 31, 2000, \texttt{dg-ga/9712005}. First version: December 8, 1997,
\texttt{dg-ga/9712005 (v1)}, \S 4--7.}
\thanks{The first author was supported in part by an NSF Mathematical 
Sciences Postdoctoral Fellowship under grant DMS 9306061 and by NSF grant
DMS 9704174} 

\begin{abstract}
  In this article, a continuation of \cite{FL2a}, we complete the
proof---for a broad class of four-manifolds---of Witten's conjecture that
the Donaldson and Seiberg-Witten series coincide, at least through terms of
degree less than or equal to $c-2$, where $c=-\frac{1}{4}(7\chi+11\sigma)$
and $\chi$ and $\sigma$ are the Euler characteristic and signature of the
four-manifold. We use our computations of Chern classes for the virtual
normal bundles for the Seiberg-Witten strata from the companion article
\cite{FL2a}, a comparison of all the orientations, and the PU(2) monopole
cobordism to compute pairings with the links of level-zero Seiberg-Witten
moduli subspaces of the moduli space of PU(2) monopoles.  These
calculations then allow us to compute low-degree Donaldson invariants in
terms of Seiberg-Witten invariants and provide a partial verification of
Witten's conjecture.
\end{abstract}

\maketitle


\section{Introduction}
\label{sec:Intro2}
\subsection{Main results}
\label{subsec:Statement}
The purpose of the present article, a continuation of \cite{FL2a}, is to
complete the proof that Witten's conjecture \cite{Witten} relating the
Donaldson and Seiberg-Witten invariants holds in ``low degrees'' for a
broad class of four-manifolds, using the $\PU(2)$-monopole cobordism of
Pidstrigatch and Tyurin \cite{PTLocal}. We assume throughout that $X$
is a closed, connected, smooth four-manifold with an orientation for which
$b_2^+(X)>0$. The Seiberg-Witten (SW) invariants (see
\S \ref{subsec:SWinvariants}) comprise a function, $SW_X:\Spinc(X)\to\ZZ$,
where $\Spinc(X)$ is the set of isomorphism classes of \spinc structures on
$X$. For $w\in H^2(X;\ZZ)$, define
\begin{equation}
\label{eq:SWSeries}
\bS\bW_X^{w}(h)
=
\sum_{\fs \in \Spinc(X)}(-1)^{\half(w^{2}+c_{1}(\fs)\cdot w)}
SW_X(\fs)e^{\langle c_{1}(\fs),h\rangle},
\quad
h\in H_2(X;\RR),
\end{equation}
by analogy with the structure of the Donaldson series $\bD_X^w(h)$
\cite[Theorem 1.7]{KMStructure}. There is a map $c_1:\Spinc(X)\to
H^2(X;\ZZ)$ and the image of the support of $SW_X$ is the set $B$ of
SW-basic classes \cite{Witten}. A four-manifold $X$ has SW-simple type when
$b_1(X)=0$ if $c_1(\fs)^2=2\chi+3\sigma$ for all $c_1(\fs)\in B$,
where $\chi$ and $\sigma$ are the Euler characteristic and signature of
$X$.  Let $B^\perp\subset H^2(X;\ZZ)$ denote the orthogonal complement
of $B$ with respect to the intersection form $Q_X$ on $H^2(X;\ZZ)$. Let
$c(X)=-\frac{1}{4}(7\chi+11\sigma)$. As stated in \cite{FL2a}, we have:

\begin{thm}
\label{thm:DSWSeriesRelation}
Let $X$ be four-manifold with $b_1(X)=0$ and odd $b_2^+(X)\geq 3$. Assume
$X$ is abundant, SW-simple type, and effective. For any $\La\in B^\perp$
and $w\in H^2(X;\ZZ)$ for which $\La^2=2-(\chi+\si)$ and $w-\Lambda\equiv
w_2(X)\pmod{2}$, and any $h\in H_2(X;\RR)$, one has
\begin{equation}
\label{eq:LowDegreeEquality}
\begin{aligned}
\bD^{w}_X(h) 
&\equiv 0\equiv \bS\bW^{w}_X(h) \pmod{h^{c(X)-2}},
\\
\bD^{w}_X(h)
&\equiv 
2^{2-c(X)}e^{\frac{1}{2}Q_X(h,h)}\bS\bW^{w}_X(h) 
\pmod{h^{c(X)}}.
\end{aligned}
\end{equation}
\end{thm}

The order-of-vanishing assertion for the series $\bD^{w}_X(h)$ and
$\bS\bW^{w}_X(h)$ in equation \eqref{eq:LowDegreeEquality} was proved in
joint work with Kronheimer and Mrowka \cite{FKLM}, based on the results in
an earlier version \cite{FL2} of this article and its companion \cite{FL2a}.

The background material underlying the statement of Theorem
\ref{thm:DSWSeriesRelation}---including the definition and significance of
``abundant'' and ``effective'' four-manifolds---was explained in \cite[\S
1]{FL2a}, so we refer the reader to \cite{FL2a} for details and just
briefly mention here some aspects of the statement which may be less
familiar.

As customary, $b_2^+(X)$ denotes the dimension of a maximal
positive-definite linear subspace $H^{2,+}(X;\RR)$ for the intersection
pairing $Q_X$ on $H^2(X;\RR)$. It is implicit in the statement of Theorem
\ref{thm:DSWSeriesRelation} that we have selected an orientation for
$H^1(X;\RR)\oplus H^{2,+}(X;\RR)$, and the Donaldson and Seiberg-Witten
invariants are computed with respect to this choice.

A four-manifold is {\em abundant\/} if the restriction of $Q_X$ to
$B^\perp$ contains a hyperbolic sublattice \cite[Definition
1.2]{FL2a}. This condition ensures that there exist classes $\Lambda\in
B^\perp$ with prescribed even square, such as $\La^2=2-(\chi+\si)$. All
compact, complex algebraic, simply connected surfaces with $b_2^+\geq 3$
are abundant. There exist simply connected four-manifolds with $b_2^+ \geq
3$ which are not abundant, but which nonetheless admit classes $\Lambda\in
B^\perp$ with prescribed even squares \cite[p. 175]{FKLM}.

As described in \cite[Definition 1.3]{FL2a}, a four-manifold is
{\em effective\/} if it satisfies Conjecture 3.1 in \cite{FKLM},
restated as Conjecture \ref{conj:Multiplicity} in this article.
This conjecture asserts that the pairings of Donaldson-type cohomology
classes with the link of a Seiberg-Witten moduli subspace of the
(compactified) moduli space of $\PU(2)$ monopoles are multiples of its
Seiberg-Witten invariants, so these pairings are zero when the Seiberg-Witten
invariants for that Seiberg-Witten moduli space are trivial. 

For any $w\in H^{2}(X;\ZZ)$, one can define a Donaldson
invariant (see \S \ref{subsubsec:DefnDInvariants} for a detailed
description) as a real-linear function \cite[p. 595]{KMStructure}
$$
D^{w}_{X}:\AAA(X) \to \RR,
$$
where \cite{KMStructure}
\begin{equation}
  \label{eq:AAAXdefn}
\AAA(X) = \Sym\left( H_{\even}(X;\RR)\right)\otimes
\La^\bullet(H_{\odd}(X;\RR))
\end{equation}
is the graded algebra. If $z\in\AAA(X)$ is a
monomial then $D_X^w(z)=0$ unless its degree satisfies
\begin{equation}
\label{eq:Mod8}
\deg(z) \equiv -2w^2 -\textstyle{\frac{3}{2}}(\chi+\sigma) \pmod{8}.
\end{equation}
Recall from \cite[Equation (1.5)]{KMStructure} that the
Donaldson series is defined by
\begin{equation}
  \label{eq:DSeries}
  \bD^{w}_{X}(h) = D^{w}_{X}((1+{\textstyle{\frac{1}{2}}} x)e^{h})
=
\sum_{d\geq 0} {\textstyle{\frac{1}{d!}}}D_X^w(h^d) 
+ {\textstyle{\frac{1}{2d!}}}D_X^w(xh^d),
\quad h\in H_2(X;\RR).
\end{equation}
A four-manifold with $b_1(X)=0$ and odd $b_2^+(X)\geq 3$
has {\em Kronheimer-Mrowka (KM) simple type\/}
\cite{KMStructure} if for some $w$ and all $z\in \AAA(X)$,
$$
D^{w}_{X}(x^{2}z)=4D^{w}_{X}(z).
$$
If as in Theorem \ref{thm:DSWSeriesRelation}, we do not assume that $X$ has
KM-simple type, then equation
\eqref{eq:DSeries} only defines $\bD^{w}_{X}(h)$ as a formal power series
and one cannot necessarily recover all invariants of the form
$D_X^w(x^mh^{d-2m})$ from the series \eqref{eq:DSeries}.  According to
\cite[Theorem 1.7]{KMStructure}, when $X$ has KM-simple type
the series $\bD^{w}_{X}(h)$ is an analytic function of $h$ and there are
finitely many characteristic cohomology classes $K_{1},\ldots,K_{m}$ (the
KM-basic classes) and constants $a_{1},\ldots,a_{m}$ (independent of $w$)
so that
$$
\bD^{w}_X(h)
=
 e^{\half Q_X(h,h)}\sum_{r=1}^{s}(-1)^{\half(w^{2}+K_{r}\cdot w)}
a_{r}e^{\langle K_{r},h\rangle},
\quad h\in H_2(X;\RR).
$$
More generally \cite{KMGeneralType}, a four-manifold $X$ has {\em finite
type\/} or {\em type $\tau$\/} if
$$
D_X^w((x^2-4)^\tau z) = 0,
$$ 
for some $\tau\in\NN$ and all $z\in\AAA(X)$. Kronheimer and Mrowka
conjectured that all four-manifolds $X$ with $b_2^+(X)>1$ have finite type
and state an analogous formula for the series $\bD^{w}_X(h)$; proofs of
different parts of their conjecture have been reported Fr\o yshov
\cite[Corollary 1]{FroyshovFiniteType}, Mu\~noz
\cite[Corollary 7.2 \& Proposition 7.6]{MunozFFHomology}, and Wieczorek
\cite[Theorem 1.3]{WWFiniteType}.

For a four-manifold $X$ with $b_1(X)=0$ and odd $b_2^+(X)\geq 3$, Witten's
conjecture \cite{Witten} asserts that $X$ has KM-simple type if and only if
it has SW-simple type; if $X$ has simple type, then
\begin{equation}
\label{eq:WittenConjSeries}
    \bD^{w}_{X}(h)=2^{2-c(X)}e^{\half Q_X(h,h)}\bS\bW^{w}_{X}(h),
\quad h\in H_2(X;\RR).
\end{equation}
Equation \eqref{eq:LowDegreeEquality} therefore tells us that Witten's
formula holds, modulo terms of degree greater than or equal to $c(X)$, at
least for four-manifolds satisfying the hypotheses of Theorem
\ref{thm:DSWSeriesRelation}. Equation
\eqref{eq:LowDegreeEquality} is proved by considering Seiberg-Witten moduli
spaces in the top level, $\ell=0$, of the compactified $\PU(2)$ monopole
moduli space; in order to prove that equation \eqref{eq:WittenConjSeries}
holds modulo $h^d$ for arbitrary $d\geq c(X)$ (and the same $w$,
$\Lambda$), one needs to compute the contributions of Seiberg-Witten moduli
spaces in arbitrary levels $\ell\geq 0$. In \cite{FLLevelOne} we use the
case $\ell=1$ to show that equation
\eqref{eq:WittenConjSeries} holds mod $h^{c(X)+2}$.

Equation \eqref{eq:LowDegreeEquality} is a special case of a more general
formula for Donaldson invariants which we now describe; the hypotheses
still include an important restriction which guarantees that the only
Seiberg-Witten moduli spaces with non-trivial invariants lie in the top
level of the $\PU(2)$-monopole moduli space.  When $b_1(X)\geq 0$,
the Seiberg-Witten invariants for $(X,\fs)$, with $\fs\in\Spinc(X)$, are
defined collectively as a real-linear function (see \S
\ref{subsec:SWinvariants} for a detailed description),
\begin{equation}
\label{eq:IntroSWFunction}
SW_{X,\fs}:\BB(X) \to \RR,
\end{equation}
where the graded algebra is given by
\begin{equation}
  \label{eq:SWGradedAlgebra}
  \BB(X) = \RR[x]\otimes\Lambda^\bullet(H_1(X;\RR)).
\end{equation}
Here, $\Lambda^\bullet(H_1(X;\RR))$ is the exterior algebra on $H_1(X;\RR)$,
with $\gamma\in H_1(X;\RR)$ having degree one, and $\RR[x]$ is the
polynomial algebra with generator $x$ of degree two. If $z\in\BB(X)$ then
$SW_{X,\fs}(z)=0$ unless its degree satisfies
$$
\deg(z) = d_s(\fs),
$$
where $d_s(\fs)$ is the dimension of the Seiberg-Witten moduli space $M_{\fs}$,
\begin{equation}
\label{eq:SWDimension}
d_s(\fs) = \textstyle{\frac{1}{4}}(c_1(\fs)^2 - 2\chi - 3\sigma).
\end{equation}
If $z=x^m\in\BB(X)$ and $2m=d_s(\fs)$, then as customary \cite{KMThom},
\cite{MorganSWNotes}, \cite{Witten} one has
\begin{equation}
  \label{eq:RestrictedSWInvtDefn}
SW_{X}(\fs)
=
SW_{X,\fs}(z).
\end{equation}
As in \cite[\S 1]{OzsvathSzaboAdjunct}, when $b_1(X)\geq 0$ we call
$c_1(\fs)\in H^2(X;\ZZ)$ an {\em SW-basic class\/} if the Seiberg-Witten
function \eqref{eq:IntroSWFunction} is non-trivial. If $b_2^+(X)=1$ or
$b_1(X)>0$, there are examples of four-manifolds whose basic classes have
positive-dimensional Seiberg-Witten moduli spaces: $\CC\PP^2$ and its
blow-ups give examples with $b_2^+(X)=1$
\cite{SalamonSWBook}, and connected sums of $S^1\times S^3$ and a
four-manifold with non-trivial Seiberg-Witten invariants provide examples
with $b_1(X)>0$ \cite[\S 2]{OzsvathSzaboAdjunct}.

For $\Lambda\in H^{2}(X;\ZZ)$, define
\begin{equation}
\label{eq:PositiveDiracIndexFunction}
i(\La)=\Lambda^2 + c(X) + \chi + \sigma.
\end{equation}
If $S(X)\subset\Spinc(X)$ is the subset yielding non-trivial Seiberg-Witten
functions \eqref{eq:IntroSWFunction}, let
\begin{equation}
\label{eq:SWInTopLevelFunction}
r(\Lambda,c_1(\fs))
=
-(c_1(\fs)-\Lambda)^2 - \textstyle{\frac{3}{4}}(\chi+\sigma)
\quad\text{and}\quad
r(\Lambda)
=
\mathop{\min}\limits_{\fs\in S(X)} r(\Lambda,c_1(\fs)).
\end{equation}
See Remark \ref{rmk:WhereDoTheSWSpacesLie} for a discussion of the
significance of $r(\Lambda,c_1(\fs))$ and $r(\Lambda)$.  We then have:

\begin{thm}
\label{thm:Main}
Let $X$ be a four-manifold with $b_2^+(X) \geq 1$. Assume
$\alpha\smile\alpha'=0$ for every $\alpha,\alpha'\in H^1(X;\ZZ)$ and that
$X$ is effective.  Suppose $\Lambda, w\in H^{2}(X;\ZZ)$ are classes such
that $w-\Lambda\equiv w_{2}(X) \pmod 2$ and, if $b_2^+(X)=1$, the class
$w\pmod 2$ admits no torsion integral lifts.  Let $z=x^{\delta_0}\vartheta
h^{\delta_2}$, where $h\in H_2(X;\RR)$, $\vartheta \in
\Lambda^{\delta_1}(H_1(X;\RR))$, and $x\in H_0(X;\ZZ)$ is the positive
generator, and write $\deg(z)=2\delta$, for $\delta\in\frac{1}{2}\ZZ$.
\medskip

\noindent{\em{\bf (a)}} If $\delta < i(\Lambda)$ and
$\delta < r(\Lambda)$, then
\begin{equation}
\label{eq:Dzero}
D^w_X(z)  = 0.
\end{equation}

\noindent{\em{\bf (b)}} If $\delta < i(\Lambda)$ and
$\delta=r(\Lambda)$, then 
\begin{equation}
\label{eq:Main}
\begin{aligned}
D^w_X(z) &= 2^{1-\quarter(i(\La)-\delta)-\delta_2-2\delta_0}
(-1)^{\delta_0+\delta_1+\half(\sigma-w^2)+1}
\\
&\quad\times
\sum_{\{\fs\in S(X):\  r(\Lambda,c_1(\fs))=r(\Lambda)\}}
(-1)^{\frac{1}{2}(w^{2}+c_{1}(\fs)\cdot(w-\Lambda))}
\\
&\quad\times
H_{\chi,\si}(\La^2,\deg(z),d_s(\fs),\delta_1)
SW_{X,\fs}(\vartheta x^{\half(d_s(\fs)-\delta_1)}) \langle
c_1(\fs)-\Lambda,h\rangle^{\delta_2},
\end{aligned}
\end{equation}
with $H_{\chi,\si}$ defined in equation
\eqref{eq:MainJacobiCoefficient}. If $\half(d_s(\fs)-\delta_1)=0$, then
$H_{\chi,\si}(\La^2,\deg(z),d_s(\fs),\delta_1) = 1$.
If $b_2^+(X)=1$, then all invariants in equation \eqref{eq:Main}
are evaluated with respect to the chambers determined by the same
period point in the positive cone of $H^2(X;\RR)$.

\noindent{\em{\bf (c)}} If $\Lambda^2$ and $\delta$ satisfy
$$
2r(\Lambda) < 2\delta 
\leq 
r(\Lambda)+\textstyle{\frac{1}{2}}(r(\Lambda)+i(\Lambda)) - 2,
$$
then equation \eqref{eq:Main} holds with $D^w_X(z) = 0$.
\end{thm}

\begin{rmk}
\label{rmk:H1andH3} 
\begin{enumerate}
\item
When $\delta_1> d_s(\fs)$, the Seiberg-Witten invariant
$SW_{X,\fs}(\vartheta x^{\half(d_s(\fs)-\delta_1)})$ in equation
\eqref{eq:Main} is zero by definition.
\item
If $z=Yz'$ where $Y\in H_3(X;\ZZ)$ and $z'$ cannot be written as $z'=xz''$
for $x\in H_0(X;\ZZ)$, then equations \eqref{eq:Dzero} and \eqref{eq:Main}
hold but the right-hand-side of \eqref{eq:Main} vanishes.
\end{enumerate}
\end{rmk}

The hypothesis in Theorems \ref{thm:DSWSeriesRelation}, \ref{thm:Main},
\ref{thm:FLthm}, and Corollary \ref{cor:FLthm} 
that $X$ is effective can be eliminated if, in the
definition \eqref{eq:SWInTopLevelFunction} of $r(\Lambda)$, we replace
$S(X)$ with the (possibly larger) set of all $\fs\in \Spinc(X)$ for which
the Seiberg-Witten moduli space $M_{\fs}$ (as defined in
\cite[\S 2.3]{FL2a}, with perturbations depending on $\La$) is non-empty. This
additional generality does not seem to be useful in practice, however.

When $b_2^+(X)=1$ and $w\pmod 2$ admits no torsion integral lifts,
Lemma \ref{lem:IdentifyChambers} implies that
the walls defining the chamber structure for the Donaldson invariant
$D^w_X$ are given by the walls for the Seiberg-Witten invariants
appearing in equation \eqref{eq:Main}.  Thus, both sides of the
equation will change when the period point crosses one of these walls.

The function $H_{\chi,\sigma}$ in equation \eqref{eq:Main} is defined as
\begin{equation}
\label{eq:MainJacobiCoefficient}
H_{\chi,\si}(\La^2,\deg(z),d_s(\fs),\delta_1)
=
2^{d}P^{a,b}_d(0),
\end{equation}
where $d$ is a natural number and $a$, $b$ are integers given by
\begin{align*}
d &= \textstyle{\frac{1}{2}}(d_s(\fs) - \delta_1),
\\
a &= \textstyle{\frac{1}{4}}\left( 3r(\La) + i(\La)\right) 
- \textstyle{\frac{1}{2}}\deg(z)-d-2,
\\
b &= \textstyle{\frac{1}{2}}\left(\deg(z)- 2r(\La) - d_s(\fs)\right) 
-\textstyle{\frac{1}{4}}(\chi+\si),
\end{align*}
and $P^{a,b}_d(\xi)$ is a {\em Jacobi polynomial\/} \cite[\S 8.96]{GR},
\begin{equation}
\label{eq:IntroJacobiPolynomial}
P^{a,b}_d(\xi) 
=
\sum_{u=0}^{d}
{\binom{a+d}{d-u}}{\binom {b+d}{u}}(\xi-1)^u(\xi+1)^{d-u},
\quad \xi\in\CC.
\end{equation}
The polynomials $P^{a,b}_d(\xi)$ may in turn be expressed in terms of {\em
hypergeometric functions\/} \cite[\S 9.10]{GR},
\cite{Luke}, as we explain in \S \ref{subsec:Computation}. Since
$c_1(\fs)^2\equiv \sigma\pmod 8$, equation \eqref{eq:SWDimension} for
$d_s$ implies that the expression $-\half d_s-\quarter\left(\chi
+\sigma\right)$ in the definition of $b$ is an integer.

When $\Lambda\in B^\perp\subset H^{2}(X;\ZZ)$ and $X$ has SW-simple type,
the expression \eqref{eq:SWInTopLevelFunction} for $r(\Lambda)$ becomes
\begin{equation}
\label{eq:SWSimpleUDepthDiracIndexParam}
r(\Lambda) =  -\Lambda^{2} + c(X) - (\chi+\sigma),
\end{equation}
and $i(\Lambda)+r(\Lambda)=2c(X)$, by equation
\eqref{eq:PositiveDiracIndexFunction} for $i(\Lambda)$. 
Theorem \ref{thm:Main} then simplifies to:

\begin{thm}
\label{thm:FLthm}
Let $X$ be a four-manifold with odd $b_2^+(X)\geq 3$ and $b_1(X)=0$. Assume
that $X$ is effective and has SW-simple type.  Suppose that $\Lambda\in
B^\perp$ and that $w \in H^{2}(X;\ZZ)$ is a class with $w-\Lambda\equiv
w_{2}(X) \pmod 2$.  Let $\delta \geq 0$ and $0\leq m\leq
[\delta/2]$ be integers. 
\medskip

\noindent{\em{\bf (a)}}
If $\delta< i (\Lambda)$ and $\delta < r(\Lambda)$, then for all $h\in
H_2(X;\RR)$ we have
\begin{equation}
    \label{eq:DVanishingIntro}
D^{w}_{X}(h^{\delta-2m}x^{m})= 0. 
\end{equation}
\noindent{\em{\bf (b)}}
If $\delta<i(\Lambda)$ and $\delta = r(\Lambda)$ we have
\begin{equation}
    \label{eq:DSWrel}
\begin{aligned}
    D^{w}_{X}(h^{\delta-2m}x^{m})
&=
    2^{1-\half (c(X)+\delta)}(-1)^{m+1+\half (\sigma-w^2)}
\\
&\qquad\times \sum_{\fs \in \text{\rm Spin}^c(X)}
    (-1)^{\half(w^{2}+c_{1}(\fs)\cdot w)}SW_X(\fs)
    \langle c_{1}(\fs)-\Lambda,h \rangle^{\delta -2m}.
\end{aligned}
\end{equation}
\noindent{\em{\bf (c)}} 
If $r(\Lambda)<\delta\leq \frac{1}{2}(r(\Lambda)+c(X)-2)$, then equation
    \eqref{eq:DSWrel} holds with $D^{w}_{X}(h^{\delta-2m}x^{m})=0$.
\end{thm}

We recall from \cite[\S 6.1.1]{FrM} that the symmetric algebra
$\Sym(H_2(X;\RR))$ is canonically isomorphic (as a graded algebra) to the
algebra $P(H_2(X;\RR))$ of polynomial functions on $H_2(X;\RR)$. Thus,
given the Donaldson invariants $D^{w}_{X}(h^{\delta-2m}x^{m})$, we can
recover all invariants of the form $D^{w}_{X}(h_1\cdots h_{\delta-2m}x^{m})$.

\begin{cor}
\label{cor:FLthm}
\cite[Theorem 1.1]{FKLM}
Let $X$ be a four-manifold with odd $b_2^+(X)\geq 3$ and $b_1(X)=0$. Assume
that $X$ is effective, abundant, and SW-simple type, with $c(X)\geq 3$. Then
for any $w\in H^2(X;\ZZ)$ with $w\equiv w_2(X)\pmod{2}$ we have
$$
\bS\bW_X^w(h) \equiv 0 \pmod{h^{c(X)-2}}.
$$
\end{cor}

In \S \ref{subsec:ProofsOfMainTheorems} we give a slightly different
and more geometric proof of Corollary \ref{cor:FLthm} than provided in
\cite[\S 2]{FKLM}, using the final case of Theorem \ref{thm:FLthm}. 
This result proves a conjecture of Mari\~no, Moore, and Peradze
\cite{MMPdg}, \cite{MMPhep} for four-manifolds of Seiberg-Witten simple
type, albeit with the additional hypotheses that the four-manifolds are
abundant and effective. The vanishing result for Donaldson invariants in
Theorem \ref{thm:DSWSeriesRelation} can be sharpened: see Theorem 1.2 in
\cite{FKLM}.

\subsection{Remarks on the hypotheses of Theorems \ref{thm:Main},
\ref{thm:FLthm}, and \ref{thm:DSWSeriesRelation}}
To prove Theorem \ref{thm:Main} (and thus Theorems
\ref{thm:DSWSeriesRelation}, \ref{thm:FLthm}, and Corollary
\ref{cor:FLthm}), we employ the moduli space of
$\PU(2)$ monopoles, $\sM_{\ft}/S^1$, as a cobordism between a link of the
moduli space of anti-self-dual connections, $M^w_{\ka}$, and the links of
moduli spaces of Seiberg-Witten monopoles, $M_{\fs}$.  Our application of
the cobordism method in this article requires that 
\begin{enumerate}
\item
The codimension of $M^w_{\ka}$ in $\sM_{\ft}$,
given by twice the complex index of a Dirac operator, is positive
(used in Proposition \ref{prop:LinkOfASD}), and
\item
Only the top level of the Uhlenbeck compactification
$\bar\sM_{\ft}$ of the moduli space of $\PU(2)$ monopoles
contains Seiberg-Witten moduli spaces $M_{\fs}$ with non-trivial invariants
(used in Corollary \ref{cor:CompactReductionFormula}).
\end{enumerate}
In the proof of Assertions (a) and (b) of Theorem \ref{thm:Main}, one has
$2\delta=\deg(z)=\dim M^w_{\ka}$ and the hypotheses $\delta< i(\La)$ and
$\delta\leq r(\La)$ ensure that Conditions (1) and (2) hold, respectively.
In the proof of Assertion (c), one has $\dim M^w_{\ka}=2r(\La)$, which
implies that Condition (2) holds while the inequality in the hypothesis of
(c) implies that Condition (1) is satisfied.

Assertion (a) follows because the hypothesis $\delta<r(\Lambda)$ implies
that there are no Seiberg-Witten moduli spaces contained in $\bar\sM_{\ft}$
with non-trivial invariants, and so $\sM_{\ft}/S^1$ is essentially a
null-cobordism of the link of $M_\kappa^w$. Assertion (c) follows because
the hypothesis on $\deg(z)$ implies that $\deg(z)>M_\kappa^w$, and only the
pairings of Donaldson-type cohomology classes with links of Seiberg-Witten
moduli spaces can be non-trivial.  In the remaining Assertion (b), the
hypotheses imply that the cobordism yields an equality between pairings with
the link of $M_\kappa^w$ and a sum of pairings with the links of
$M_{\fs}$. The same remarks apply to the hypotheses in Assertions (a), (b),
and (c) in Theorem \ref{thm:FLthm}.

The assumption that $\alpha\smile\alpha'=0$ for all $\alpha,\alpha'\in
H^1(X;\ZZ)$ greatly simplifies the calculation of the Chern classes of the
virtual normal bundle of $M_{\fs}$ in $\sM_{\ft}$ (see Corollary 3.30 in
\cite{FL2a}) and hence its Segre classes (see Lemma \ref{lem:SegreOfN}). It
should be possible to remove this condition with more work, but this would
take us a little beyond the scope of this article and \cite{FL2a}. We plan
to address this point elsewhere.

When $b_2^+(X)=1$, we assume that $w\pmod 2$ does not admit a torsion
integral lift in order to avoid complications in defining the chamber in
the positive cone of $H^2(X;\RR)$ with respect to which the Donaldson and
Seiberg-Witten invariants are computed.  See the comments at the end of \S
\ref{subsubsec:DefnDInvariants} and before Lemma \ref{lem:IdentifyChambers}
for further discussion.

The proof of Theorem \ref{thm:DSWSeriesRelation} requires one to choose
classes $\Lambda\in B^\perp$ with optimally prescribed even square in order
to obtain the indicated vanishing results for the Donaldson and
Seiberg-Witten series, as well as compute the first non-vanishing
terms. The hypothesis that $X$ is abundant guarantees that one can find
such classes, though such choices are also possible for some non-abundant
four-manifolds \cite{FKLM}.

\subsection{Remarks and conjectures for formulas for Donaldson invariants
involving Seiberg-Witten strata in arbitrary levels} The following remarks
are intended to convey an outline of the remainder of our work on a proof
of Witten's conjecture in \cite{FL3}, \cite{FL4},
\cite{FL5}, \cite{FLLevelOne}. While some details still remain to be
checked, we are confident that the conclusions stated below are correct
based on our work thus far, despite their conservatively-stated current
status as conjectures rather than firm assertions.

As envisaged in \cite{FLGeorgia}, the $\PU(2)$-monopole program proposed by
Pidstrigatch and Tyurin \cite{PTLocal} for proving Witten's conjecture
\cite{MooreWitten}, \cite{Witten} uses the oriented cobordism
$\sM_{\ft}/S^1$ between 
\begin{itemize}
\item
The links $\bL_{\ft,\kappa}^w$ in $\bar\sM_{\ft}/S^1$ of the anti-self-dual
moduli subspace $M_\kappa^w$ of $\bar\sM_{\ft}/S^1$, and
\item
The links $\bL_{\ft,\fs}$ in $\bar\sM_{\ft}/S^1$ of the Seiberg-Witten
moduli subspaces, $M_{\fs}\times\Sym^\ell(X)$, of the space of ideal
$\PU(2)$ monopoles, $\cup_{\ell=0}^\8(\sM_{\ft_\ell}\times\Sym^\ell(X))$,
containing $\bar\sM_{\ft}$.
\end{itemize}
The program therefore has two principal steps, which we now outline.  The
first step is to define the links $\bL_{\ft,\fs}$ of
$M_{\fs}\times\Sym^\ell(X)$ for arbitrary $\ell\geq 0$ using the gluing
construction of \cite{FL3}, \cite{FL4}, extending the construction in
\cite{FL2a} which just treats the case $\ell=0$. The oriented cobordism
$\sM_{\ft}/S^1$ then yields a formula (with $\deg(z)=2\delta$),
\begin{equation}
\label{eq:RawGeneralCobordismFormula}
D_{X}^w(z) = - 2^{-\delta_c}\sum_{\fs\in\Spinc(X)}
(-1)^{\frac{1}{4}(w-\Lambda+c_1(\fs))^2}
\langle\mu_{p}(z)\smile\mu_{c}^{\delta_c},
[\bL_{\ft,\fs}]\rangle,
\end{equation}
where $\mu_{p}(z)$ and $\mu_{c}$ are Donaldson-type cohomology classes, and 
$\delta_c=\frac{1}{4}(i(\Lambda)-\delta)-1$.  Work in progress \cite{FL5}
then strongly indicates that the pairings on the right-hand side of
equation \eqref{eq:RawGeneralCobordismFormula} have the following general
form, when $b_1(X)=0$ and $z=x^mh^{\delta-2m}$,
\begin{equation}
\label{eq:GeneralFormulaForPairing}
\langle\mu_{p}(z)\smile\mu_{c}^{\delta_c},
[\bL_{\ft,\fs}]\rangle
=
SW_X(\fs)
\sum_{i=0}^{r}
\bdelta_{d,i}\left(\langle c_1(\fs),h\rangle,\langle\Lambda,h\rangle\right)
Q_X^{\ell-i}(h,h),
\end{equation}
where $r = \min(\ell,[\delta/2]-m)$,
$\ell=\frac{1}{4}(\delta-r(\Lambda,c_1(\fs)))$, 
$d=\delta-2(m+\ell-i)$, and $\bdelta_{d,i}$ is a degree-$d$, homogeneous
polynomial (aside from stray powers of $(-1)$ and $2$) in two variables
with coefficients which are degree-$i$ polynomials in $2\chi\pm 3\sigma$,
$(c_1(\fs)-\Lambda)^2$, $\Lambda^2$, and $(c_1(\fs-\Lambda)\cdot c_1(\fs)$.

An interesting feature of the formula
\eqref{eq:GeneralFormulaForPairing} is that one sees a
factorization of the pairings $\langle\mu_{p}(z)\smile\mu_{c}^{\delta_{c}},
[\bL_{\ft,\fs}]\rangle$ into a product of $SW_X(\fs)$ and the term
$\bdelta_d=\sum_{i=1}^r\bdelta_{d,i}Q_X^{\ell-i}$. In particular, the pairing 
\eqref{eq:GeneralFormulaForPairing} vanishes when $SW_X(\fs)=0$, implying
that $X$ is ``effective'' in the sense described in \S
\ref{subsec:Statement}.  The factors $\bdelta_{d,i}$ are similar to those
appearing in the Kotschick-Morgan conjecture \cite{FLKM1},
\cite{KotschickMorgan} for the form of the wall-crossing formula
for the Donaldson invariants of a four-manifold $X$ with
$b_2^+(X)=1$. However, $\bdelta_{d,i}$ is a polynomial in two variables
while the corresponding term in the conjectured wall-crossing formula for
Donaldson invariants is a polynomial in only one variable.

Explicit, direct computations of pairings with the links $\bL_{\ft,\fs}$ of
ideal Seiberg-Witten moduli spaces, $M_{\fs}\times\Sym^\ell(X)$, are 
possible when $\ell=0,1$ or $2$: indeed, Theorem
\ref{thm:DSWSeriesRelation} is proved using the case $\ell=0$ and we prove
an $\ell=1$ analogue of Theorem
\ref{thm:DSWSeriesRelation} in \cite{FLLevelOne}, while the case $\ell=2$
would follow by adapting work of Leness in \cite{LenessWC}. However, direct
computations appear intractable when $\ell$ is large.

The idea underlying the second step of the program is to use the existence
of formulas \eqref{eq:RawGeneralCobordismFormula} and
\eqref{eq:GeneralFormulaForPairing} in conjunction with auxiliary
arguments to prove Witten's conjecture, since more direct calculations of
the link pairings appear difficult.  The work of G\"ottsche
\cite{Goettsche} suggests that such indirect strategies should succeed, as
he was able to compute the wall-crossing formula for the Donaldson
invariants of four-manifolds with $b_2^+(X)=1$ and $b_1(X)=0$, under the
assumption that the Kotschick-Morgan conjecture \cite{KotschickMorgan}
holds for such four-manifolds.  The facts that $\bdelta_{d,i}$ is a function of
two variables and both $\chi$ and $\sigma$ may vary
independently indicate that this second step in the $\PU(2)$-monopole
program is potentially more complicated than that of
\cite{Goettsche}, where the assumption that $b_2^+(X)=1$ implies that
$\sigma=1-b_2^-(X)$ and $\chi=3+b_2^-(X)$ (when $b_1(X)=0$).  However, in
the case of the $\PU(2)$-monopole program there are more sources of
``recursion relations'' of the type used by G\"ottsche, in addition to
those arising from the blow-up formulas of Fintushel-Stern \cite{FSBlowUp},
\cite{FSTurkish} (for Donaldson and Seiberg-Witten invariants). Moreover,
there is a rich supply of examples of four-manifolds where Witten's
conjecture is known to hold.

\subsection{Remarks on Abelian localization}
One of the first observations concerning the moduli space of $\PU(2)$
monopoles is that the Donaldson stratum, $\iota(M_\kappa^w)$, and the
Seiberg-Witten strata, $\iota(M_{\fs})$, are fixed-point sets under the
circle action given by scalar multiplication on the spinor components of
$\PU(2)$-monopole pairs; see \cite[\S 3.1]{FL2a} for a detailed
account. This raises the question of whether the technique
of Abelian localization, as applied to circle actions on compact
manifolds \cite{AtiyahBott} or its generalizations to singular
algebraic varieties (for example, see \cite{GraberPand}), can be usefully
applied 
here to prove Witten's conjecture. As we indicate below, if the
moduli space of $\PU(2)$ monopoles were a compact manifold and
the Donaldson and Seiberg-Witten strata were smooth submanifolds,
then an application of the localization formula would be
equivalent to our construction of links and application of the
$\PU(2)$-monopole cobordism. There is no saving of labor and the
essential point remains, with either equivalent view, to compute
the Chern classes of the normal bundles of the fixed-point sets.
As the moduli space of $\PU(2)$ monopoles is non-compact,
equipped with the highly singular Uhlenbeck compactification or
somewhat less singular but more complicated bubble-tree
compactification, our application of the $\PU(2)$-monopole
cobordism can be thought of as an extension of the localization
method to those differential-geometric, singular settings. 

The technique of Abelian localization does not reduce the information about
neighborhoods of singularities needed to compute intersection pairings
because the localization formula requires a computation of an equivariant
Thom or Euler class of the normal bundle of the fixed point set.  For
example, if $F\subset M$ is the fixed point set and $N\to F$ is its normal
bundle, the equivariant Euler class of $N$ is the Euler class $e(N_{S^1})$
of the bundle
$$
N_{S^1} =\ES^1\times_{S^1}N \to \BS^1\times F,
$$
where $\ES^1\to\BS^1$ is the universal $S^1$ bundle over the classifying
space $\BS^1$.  Let
$$
\pi_{S^1}: \ES^1\times_{S^1}M\to \BS^1\quad\text{and}\quad
\iota:\ES^1\times_{S^1}N \to \ES^1\times_{S^1} M,
$$
be the projection and embedding maps, respectively.  If $\dim M=m$, then
because 
$$
H^m_{S^1}(M\less F;\RR)\cong H^m( (M\less F)/S^1;\RR),
$$ 
any class $x\in H^m_{S^1}(M;\RR)$ has compact support in
$\ES^1\times_{S^1}N$ by dimension counting.  The Abelian localization
formula \cite[Equation (3.8)]{AtiyahBott} states that
\begin{equation}
\label{eq:EquivariantLocalization}
(\pi_{S^1})_* x = \frac{\iota^*x}{e(N_{S^1})}/[F].
\end{equation}
If $h$ is the pullback of the universal first Chern class from $\BS^1$ to
$\BS^1\times F$ and $\pi_F:\BS^1\times F\to F$ is the projection, then the
splitting principle shows that $e(N_{S^1})=\sum_{i=0}^rh^i \pi_F^*
c_{r-i}(N)$, where $r=\rank N$.  A simple algebraic computation (see
\cite{FLLevelOne}) shows that computing $(\pi_{S^1})_* x $ using
formula \eqref{eq:EquivariantLocalization} and finding the inverse of
$e(N_{S^1})$ is equivalent to computing the Segre classes of $N$.  

Thus, if $\sM_{\ft}$ were a compact manifold and the Seiberg-Witten moduli
spaces $M_{\fs}$ were smooth submanifolds of $\sM_{\ft}$, the Abelian
localization method would be equivalent to the one used in this article.

If $X$ is a complex algebraic surface, it should be possible to construct
the Gieseker compactification for the moduli space of $\PU(2)$ monopoles
over $X$, by analogy with the construction of Morgan
\cite{MorganComparison} and Li \cite{LiAlgGeomDonPoly} for the moduli space
of $\PU(2)$ monopoles, and then apply the results of \cite{GraberPand} to
this compactification.  However, one would still need to compute the
equivariant Euler classes of the normal sheaves of strata of ideal,
reducible $\PU(2)$ monopoles, in order to apply
\cite[Equation (1)]{GraberPand}. If $X$ is not algebraic, one would need to
solve the non-trivial problem of defining the normal sheaves of these
strata in gauge-theoretic compactifications.

\subsection{A guide to the article and outline of the proofs of the main
results} The present article is a direct continuation of
\cite{FL2a} and rather than repeat many of the definitions here, we shall
refer to \cite{FL2a}.  The reader may wish to consult the notational index in
\cite{FL2a}, as well as the notational index for the present article which
appears just before \S \ref{sec:Orient}. The first version of this article
was distributed in December 1997 as sections 4--7 of the
preprint \cite{FL2}. 

As in \cite{FL2a}, we let $\fs=(\rho,W)$ denote a \spinc structure on $X$,
where $W$ is a Hermitian, rank-four bundle over $X$ and use $\ft=(\rho,V)$
to denote a ``\spinu structure'' on $X$, where $V$ is a Hermitian,
rank-eight bundle over $X$ \cite[Definition 2.2]{FL2a}. 

One of the main results (Theorem 3.31) of \cite{FL2a} is a calculation the
Chern characters of vector bundles defining tubular neighborhoods of links
of Seiberg-Witten moduli spaces in local, ``thickened'' moduli spaces of
$\PU(2)$ monopoles $\sM_{\ft}$. In \S \ref{sec:Orient} of this article, we
compare the orientations of the moduli spaces and anti-self-dual
connections, Seiberg-Witten monopoles, and their links in the moduli space
of $\PU(2)$ monopoles. In \S \ref{sec:Cohomology} we define cohomology
classes and dual geometric representatives on the moduli space of $\PU(2)$
monopoles and in \S \ref{sec:DegreeZero} we prove Theorem
\ref{thm:Main} by counting the intersection of these
representatives with links of the moduli spaces of anti-self-dual
connections and top-level Seiberg-Witten monopoles.

Compatible choices of orientations for all the moduli spaces appearing in
the stratification \cite[Equation 3.17]{FL2a} of $\sM_{\ft}$
and of the links $\bL^{w}_{\ft,\kappa}$ and $\bL_{\ft,\fs}$
are a basic requirement of the
cobordism method and a discussion of our orientation conventions is taken
up in \S \ref{sec:Orient}. As in \cite[Equation (2.40)]{FL2a}, we let
$\sM_{\ft}^{*,0}\subset\sM_{\ft}$ denote the subspace represented by $\PU(2)$
monopoles which are neither zero-section pairs (corresponding to
anti-self-dual connections) or reducible pairs (corresponding to
Seiberg-Witten monopoles). In \S \ref{subsec:Orientability} we show that
$\sM_{\ft}^{*,0}/S^1$ (the smooth locus or top stratum of $\sM_{\ft}/S^1$)
is orientable, with an orientation determined by a
choice of an orientation for the moduli space of anti-self-dual connections,
$M^w_\ka$, as explained further in \S
\ref{subsec:OrientASD}: this allows us to compute the oriented intersection
of one-manifolds with the link $\bL^{w}_{\ft,\kappa}$, where the one-manifolds
arise as the intersection of geometric representatives of the cohomology
classes on $\sM_{\ft}^{*,0}/S^1$. We also
define an orientation for $\sM_{\ft}^{*,0}/S^1$ induced by an
orientation of $M_{\fs}$, as in
\S \ref{subsec:OrientSW}, and this allows us to
compute the oriented intersection of one-manifolds with the link
$\bL_{\ft,\fs}$.  In \S \ref{subsec:OrientComparison} we compare the two
orientations of $\sM_{\ft}^{*,0}/S^1$ naturally induced by those of
$M^w_{\ka}$ and $M_{\fs}$. Finally, in \S \ref{subsec:OrientLinkReducibles}
we compare the natural orientations of these links with the orientations
obtained by considering them as boundaries of the complement in
$\sM_{\ft}/S^1$ of small open neighborhoods of the anti-self-dual and
Seiberg-Witten strata.

In \S \ref{subsec:CohomologyDef}
and \ref{subsec:GeomRepr} we describe the cohomology classes
$\mu_{p}(\beta)$ on the moduli space $M^*_{\ft}/S^1$ and $\mu_{c}$ on
$\sM_{\ft}^{*,0}/S^1$ and their dual geometric
representatives $\sV(\beta)$ and $\sW$, following the methods of
\cite{DonConn}, \cite{DonPoly}, \cite{DK}, \cite{KMStructure} for the classes
$\mu_{p}(\beta)$. A technical complication not present when
dealing solely with $M^w_{\ka}$ is that the lower strata of $\bar\sM_{\ft}$
have smaller codimension than those of $\barM^w_{\ka}$.  The
unique continuation property for reducible $\PU(2)$ monopoles \cite[Theorem
5.11]{FL1} plays a role here analogous to that of the unique continuation
property for reducible anti-self-dual $\SO(3)$ connections in
\cite{DK}, \cite{KMStructure}. In \S \ref{subsec:ExtendGeomRepr},
we show that the closures $\bar\sV(\beta)$ and $\bar\sW$ of these geometric
representatives in $\bar\sM_{\ft}/S^1$ meet the lower strata of
$\bar\sM_{\ft}^{*,0}/S^1$ transversely, that is, in a subspace of the expected
codimension away from the reducible and zero section solutions appearing in
lower levels.  Thus, the closure of the one-manifolds will have boundaries
in $\barM^w_{\ka}$ or in some stratum
$M_{\fs}\times\Sigma$ of reducible $\PU(2)$ monopoles,
where $\Sigma\subset\Sym^\ell(X)$. The hypotheses
of Theorem \ref{thm:Main} exclude consideration of the
more difficult case $\ell > 0$.
In \S \ref{subsec:CohomOnASDLink} we show that the number of points,
counted with sign, in the boundaries of the one-manifolds defined by an
appropriate choice of geometric representatives in the link
$\bL^{w}_{\ft,\kappa}$ of the stratum $M^w_{\ka}$ is
given by a multiple of the Donaldson invariant, thus completing the proof of
Theorem \ref{thm:CompactReductionFormula}. In the course of proving
this result we also show that $\sM_{\ft}^{*,0}$ is nonempty for sufficiently
negative $p_1(\ft)$ --- see Proposition \ref{prop:PU(2)MonopoleExist} in \S
\ref{subsec:CohomOnASDLink}.

In \S \ref{sec:DegreeZero} we calculate the intersection of these geometric
representatives with the link $\bL_{\ft,\fs}$ of the stratum $M_{\fs}$ and
show that it is given by a multiple of the Seiberg-Witten invariant
associated to the \spinc structure $\fs$ (see Theorem
\ref{thm:DegreeZeroFormula}). The geometric representatives, in general,
intersect the strata $M_{\fs}$ in sets of higher than expected dimension,
so our calculation of the link pairings here may be viewed as a
differential-geometric analogue of the ``excess intersection theory''
calculations discussed in \cite{Fulton}. Combining the link pairing
calculations of \S \ref{sec:Cohomology} and \S
\ref{sec:DegreeZero} and applying the cobordism $\sM_{\ft}^{*,0}/S^1$ then
yields the formulas for Donaldson invariants in terms of Seiberg-Witten
invariants in Theorem \ref{thm:Main}, from which
Theorems \ref{thm:DSWSeriesRelation} and \ref{thm:FLthm}
and Corollary \ref{cor:FLthm} are derived in \S \ref{subsec:ProofsOfMainTheorems}.

\subsubsection*{Acknowledgments}
The authors thank Ron Fintushel, Tom Mrowka, Peter Ozsv\'ath, and Zolt\'an
Szab\'o for helpful conversations.  We are especially grateful to Tom
Mrowka for his many helpful comments during the course of this work, for
bringing a correction to some examples in \cite{FL2} to our attention, as
well pointing out that our Theorem 1.4 in \cite{FL2} could be elegantly
rephrased and specialized to give the version stated as Theorem
\ref{thm:FLthm} here (and as Theorem 2.1 in \cite{FKLM}).  We thank
Dietmar Salamon for generously providing us with a pre-publication version
of his book \cite{SalamonSWBook}.  We thank Victor Pidstrigatch for
inviting us to present our work at the Warwick Symposium on Geometry and
Topology (2000). We thank the Columbia University Mathematics Department,
the Institute for Advanced Study, Princeton, and the Max Planck Institute
f\"ur Mathematik, Bonn for their generous support and hospitality during a
series of visits while this article and its companion \cite{FL2a} were
being prepared.  Finally, we thank the anonymous referee and Simon
Donaldson for their editorial suggestions and comments on the previous
versions of this article.


\twocolumn
\centerline{{\sc Index of Notation}}
\bigskip

$\AAA(X)$ 
\hfill Equation \eqref{eq:AAAXdefn}

$\BB(X)$ 
\hfill Equation \eqref{eq:SWGradedAlgebra}

$\sB_\kappa^w$ 
\hfill \cite[\S 2.1.6]{FL2a}

$C_{\chi,\si}(\cdot)$
\hfill Theorem \ref{thm:DegreeZeroFormula}

$\sC_{\ft}$, $\tsC_{\ft}$, $\sC^{*,0}_{\ft}$, $\sC^*_{\ft}$, $\sC^0_{\ft}$
\hfill \cite[\S 2.1.5]{FL2a}

$D_A$, $D_{A,\vartheta}$, $D_B$, $D_{B,\vartheta}$
\hfill Equation \eqref{eq:PertDirac}

$\bsD_{\ft}$ 
\hfill \S \ref{subsec:Orientability}

$\sD^n$, $\sD^t$ 
\hfill \cite[Equation (3.35)]{FL2a}

$D^w_X$
\hfill Equation \eqref{eq:DefineDInvarBlowUp}

$\bD^w_X$
\hfill Equation \eqref{eq:DSeries}

$\FF^w_{\ka}$ 
\hfill Equation \eqref{eq:ConnectionUniversalSO(3)Bundle}

$\FF_{\ft}$ 
\hfill Equation \eqref{eq:UniversalSO(3)Bundle}

$\sG_{\fs}$
\hfill \cite[\S 2.3]{FL2a}

$\sG_{\ft}$
\hfill \cite[Definition 2.6]{FL2a}

$H_{\chi,\sigma}(\cdot)$
\hfill Equation \eqref{eq:MainJacobiCoefficient}

$\bK_{A,\delta}$
\hfill Equation \eqref{eq:bKAdelta}

$\LL_{\fs}$
\hfill Equation \eqref{eq:UniversalSWDefinition}

$\LL_{\ft}$
\hfill Equation \eqref{eq:DefineC1}

$\bL^w_{\ft,\ka}$, $\bL^{w,\eps}_{\ft,\ka}$
\hfill Equation \eqref{eq:ASDLink}

$\bL_{\ft,\fs}$
\hfill \cite[Definition 3.22]{FL2a}

$M^w_\ka$, $\bar M^w_\ka$  
\hfill \cite[\S 2.1.6 \& 2.2]{FL2a}

$M_{\fs}$, $\tM_{\fs}$
\hfill \cite[\S 2.3]{FL2a}

$\sM_{\ft}$, $\sM^{*,0}_{\ft}$
\hfill \cite[Equations (2.33) \& (2.40)]{FL2a}

$\bar\sM_{\ft}$
\hfill \cite[Equation (2.46)]{FL2a}

$\bar\sM^{\asd}_{\ft}$
\hfill Equation \eqref{eq:ZeroSectionCompactification}

$\bar\sM^{*}_{\ft}$, $\bar\sM^0_{\ft}$, $\bar\sM^{*,0}_{\ft}$
\hfill Equation \eqref{eq:PartClosedMtSubspaces}

$\bar\sM^{\ge \eps}_{\ft}$, $\bar\sM^{*,\ge\eps}_{\ft}$
\hfill Equation \eqref{eq:MtASDChoppedOff}

$\sM_{\ft}(\Xi,\fs)$
\hfill \cite[Definition 3.20]{FL2a}

$N_{\ft}(\Xi,\fs)$, $\tN_{\ft}(\Xi,\fs)$,
\hfill \cite[\S 3.5]{FL2a}

$N^\eps_{\ft}(\Xi,\fs)$, $N^{\le\eps}_{\ft}(\Xi,\fs)$
\hfill \S \ref{subsec:OrientLinkReducibles}

$O^{\asd}(\Omega,w)$
\hfill Definition \ref{defn:ASDOrient}

$O^{\red}(\Omega,\ft,\fs)$
\hfill Definition \ref{defn:ReducibleOrient}

$\sO_A$, $\sO_A'$
\hfill Lemma \ref{lem:DeformingV}

$P^{a,b}_d$
\hfill Equation \eqref{eq:IntroJacobiPolynomial}

$SW_{X}$
\hfill Equation \eqref{eq:RestrictedSWInvtDefn}

$SW_{X,\fs}$
\hfill Equation \eqref{eq:SWFunctionDefn}

$\bS\bW^w_X$
\hfill Equation \eqref{eq:SWSeries}

$\sV(z)$, $\sW$
\hfill \S \ref{subsubsec:GeomReprMuP1} \&
\S \ref{subsubsec:DetClassRepr}

$\bar\sV(z)$, $\bar\sW$
\hfill Definition \ref{defn:GeomReprClosure}

$\tilde X$
\hfill \S \ref{subsubsec:DefnDInvariants}

$\sZ_A$
\hfill Equation \eqref{eq:sZA}

$(a)_n$ 
\hfill Equation \eqref{eq:Comb(a)Defn}

\newpage
\phantom{{\sc Index of Notation}}
\bigskip

$c(X)$
\hfill \S \ref{subsec:Statement}

$c_1(\ft)$, $p_1(\ft)$, $w_2(\ft)$
\hfill \cite[Equation (2.19)]{FL2a}

$d_a(\ft)$ 
\hfill Equation \eqref{eq:AsdDimDiracIndex}

$d_s(\fs)$
\hfill Equation \eqref{eq:SWDimension}

$\deg(z)$
\hfill Equation \eqref{eq:AA(X)Degrees}

$\det\sD_{A,\Phi}$, $\det\bsD_{\ft}$ 
\hfill Equation \eqref{eq:defndetsDt}

$\det\delta_{\hat A}$, $\det\bdelta_E$ 
\hfill Equation \eqref{eq:defndetdeltaE}

$f_A$ 
\hfill Lemma \ref{lem:DeformingV}

$g_A$ 
\hfill Equation \eqref{eq:DefineCPDiffeom}

$i(\La)$
\hfill Equation \eqref{eq:PositiveDiracIndexFunction}

$\bl$
\hfill Equation \eqref{eq:DefineNormFunction}

$n_a(\ft)$
\hfill Equation \eqref{eq:AsdDimDiracIndex}

$n_s'(\ft,\fs)$, $n_s''(\ft,\fs)$
\hfill Equation \eqref{eq:RedefineNormalIndices}

$o(\Omega,w)$
\hfill \S \ref{subsec:OrientASD}

$o_{\ft}(w,\fs)$
\hfill Equation \eqref{eq:OrientationFactor}

$r(\La)$, $r(\La,\fs)$
\hfill Equations
\eqref{eq:SWInTopLevelFunction}
\& \eqref{eq:SWSimpleUDepthDiracIndexParam}

$\fs$ 
\hfill\cite[Definition 2.1]{FL2a}

$\fs^\pm$
\hfill \S \ref{subsec:SWinvariants}

$s_i$, $s_i(N)$
\hfill Lemma \ref{lem:Segre}

$\ft$ 
\hfill\cite[Definition 2.2]{FL2a}

$\tilde \ft$
\hfill Lemma \ref{lem:BlowUpSpinu}

$\ft_\ell$ 
\hfill\cite[Equation (2.44)]{FL2a}

$\Xi$
\hfill \cite[Definition 3.20]{FL2a}

$\bga_A$
\hfill Equation \eqref{eq:ASDKuranishiEmbedding}

$\bdelta_E$
\hfill  \S \ref{subsec:Orientability}

$\delta_{c}$
\hfill Equation \eqref{eq:defndelta_c}

$\delta_i$
\hfill Equation \eqref{eq:AA(X)Degrees}

$\delta_{p}$
\hfill Equation \eqref{eq:defndelta_p}

$\iota$ 
\hfill \cite[Equations (3.4) \& (3.9)]{FL2a}

$\ka$
\hfill \S \ref{subsec:CohomologyDef}

$\mu_{c}$
\hfill Equation \eqref{eq:DefineC1}

$\mu_{p}$
\hfill Equation \eqref{eq:DefineMuMap}

$\mu_{\fs}$
\hfill Equations \eqref{eq:SWMuMap} \& \eqref{eq:SWMuMonomial}

$\nu$
\hfill Definition \ref{defn:nu}

$\pi_{\sB}$
\hfill Equation \eqref{eq:ProjPairsToConns}

$\varrho$ 
\hfill Equation \eqref{eq:DefnGInclusion}

$\varrho_L$ 
\hfill Equation \eqref{eq:DefnCircleMultWL}

$\sigma$, $\chi$, 
\hfill \S \ref{subsec:Statement} 

$\bvarphi_A$
\hfill Equation \eqref{eq:ASDKuranishiMap}

$\om(\bK,o)$ 
\hfill Equation \eqref{eq:ComplexOrientationForK}

$\om(\bL,\rd O)$ 
\hfill Equation \eqref{eq:ASDLinkOrientation1}

$\om(\sZ\cap\bK,\rd O)$ 
\hfill Equation \eqref{eq:BoundaryOrientationKZ}

$\om(Z)$, $\omega(Z)$
\hfill Convention \ref{conv:NormalBundleOrientation}

\onecolumn


\section{Orientations of moduli spaces}
\label{sec:Orient} Following the pattern in \cite{DonOrient} and
\cite[\S 5.4 \& \S 7.1.6]{DK}, we first show that $\sM^{*,0}_{\ft}$ is an
orientable manifold and then show that its orientation is canonically
determined by a choice of homology orientation of our four-manifold $X$ and
an integral lift $w$ of $w_2(\ft)$.  The orientation for $\sM^{*,0}_{\ft}$
will be invariant under the circle action and thus give an orientation for
$\sM^{*,0}_{\ft}/S^1$. We also obtain relations between the orientations of
the  smooth, top stratum $\sM^{*,0}_{\ft}/S^1$, the stratum
$M^w_{\ka}\hookrightarrow \sM_{\ft}/S^1$ defined by the anti-self-dual
moduli space, and the strata $M_{\fs}\hookrightarrow \sM_{\ft}/S^1$
defined by the Seiberg-Witten moduli spaces.

\subsection{Orientability of moduli spaces of PU(2) monopoles}
\label{subsec:Orientability}
In this section we show that $\sM^{*,0}_{\ft}$ is orientable.

As in \cite[\S 2.1.5]{FL2a}, we let $\tsC_{\ft}$ denote the
pre-configuration space of pairs $(A,\Phi)$, where $A$ is a spin connection
on $V=V^+\oplus V^-$ with fixed determinant connection
$A^{\det}=2A_\Lambda$ on $\det(V^+)$ and $\Phi$ is a section of $V^+$. We
defined $\sC_{\ft}=\tsC_{\ft}/\sG_{\ft}$, where $\sG_{\ft}$ is the group
\spinu automorphisms of $V$ \cite[Definition 2.6]{FL2a}.

Recall that $\sD_{A,\Phi} = d^{*,0}_{A,\Phi}+d^1_{A,\Phi}$ is the
``deformation operator'' corresponding to the elliptic
deformation complex \cite[Equation (2.47]{FL2a} for the moduli
spaces $\sM_{\ft}$. Let $\det \bsD_{\ft}$ be the real determinant
line bundle over the pre-configuration space $\tsC_{\ft}$, with
fiber over $(A,\Phi)\in\tsC_{\ft}$ given by
\begin{equation}
\label{eq:defndetsDt}
\det\sD_{A,\Phi} = \Lambda^{\text{max}}(\Ker\sD_{A,\Phi})\otimes
\Lambda^{\text{max}}(\Coker\sD_{A,\Phi})^*.
\end{equation}
(See \cite[\S 5.2.1]{DK} for the construction of determinant line
bundles for families of elliptic operators.)  The kernel and cokernel of
$\sD_{A,\Phi}$ are equivariant with respect to the actions of group
$\sG_{\ft}\times_{\{\pm 1\}}S^1$. The stabilizer, in $\sG_{\ft}\times_{\{\pm
  1\}}S^1$, of the pair $(A,\Phi)$ acts trivially on the fibers of
$\det\bsD_{\ft}$ because this stabilizer is connected and the structure
group of the fiber of $\det\bsD_{\ft}$ (which is a real line) is $\{\pm
1\}$ and thus disconnected.  Hence, the bundle $\det\bsD_{\ft}$ descends to
$\sC_{\ft}/S^1$ and so to $\sC_{\ft}$ as well.  We will show that the
bundle $\det\bsD_{\ft}\to \sC_{\ft}/S^1$ is trivial.

Motivated by the remarks of \cite[p. 330]{FrM}, we say that an orientation
for $\sM_{\ft}$ is a choice of orientation for the real line bundle
$\det\bsD_{\ft}$ (restricted to $\sM_{\ft}$): an orientation of the fibers
of $\det\bsD_{\ft}$ gives orientations for the real lines
$$
\det\sD_{A,\Phi} = \Lambda^{\text{max}}(H^1_{A,\Phi})
\otimes\Lambda^{\text{max}}(H^0_{A,\Phi}\oplus H^2_{A,\Phi})^*,
\qquad [A,\Phi]\in \sM_{\ft}.
$$
If $[A,\Phi]$ is a smooth point of $\sM_{\ft}$, so
$\Coker\sD_{A,\Phi} \cong H^0_{A,\Phi}\oplus H^2_{A,\Phi} = 0$,
then $\Ker\sD_{A,\Phi} = H^1_{A,\Phi} = T_{A,\Phi}\sM_{\ft}$ and
$$
\det\sD_{A,\Phi} = \Lambda^{\text{max}}(H^1_{A,\Phi}) =
\La^{\max}(T_{[A,\Phi]}\sM_{\ft}),
$$
so an orientation for $\det\bsD_{\ft}$ defines an orientation for the
open manifold $\sM^{*,0}_{\ft}$ of smooth points of $\sM_{\ft}$.
Therefore, $\det\bsD_{\ft}$ is an orientation bundle for $\sM_{\ft}$ and
$\sM_{\ft}$ is orientable if $\det\bsD_{\ft}$ is trivial.  As in
\cite{DonOrient}, we show that $\sM_{\ft}$ is orientable because the bundle
$\det\bsD_{\ft}\to \sC_{\ft}/S^1$ has a nowhere vanishing section, that is,
its first Stiefel-Whitney class vanishes.

Suppose $E\to X$ is a rank-two, complex Hermitian bundle with $c_1(E)=w$ and
$p_1(\su(E))=-4\kappa$. Denote the group of determinant-one, unitary
automorphisms of $E$ by $\sG_\kappa^w$ and the space of $\SO(3)$
connections on $\su(E)$ by $\sA_\kappa^w$.  Over the quotient space of
connections $\sB^w_{\ka}=\sA_\kappa^w/\sG_\kappa^w$ there is an orientation
bundle $\det\bdelta_E$ (see \cite[Equation (5.4.2)]{DK}) with fiber over
$[\hat A]\in \sB^w_{\ka}$ given by
\begin{equation}
\label{eq:defndetdeltaE}
\det\delta_{\hat A}
=
\Lambda^{\text{max}}(\Ker\delta_{\hat A})
\otimes\Lambda^{\text{max}}(\Coker\delta_{\hat A}),
\end{equation}
coming from the rolled-up deformation complex for the anti-self-dual
moduli space $M_\kappa^w$,
\begin{equation}
\label{eq:ASDOrientationSequence}
\delta_{\hat A} = d^*_{\hat A}+d^+_{\hat A}:
C^\8(\Lambda^1\otimes\su(E)) \to
C^\8((\Lambda^0\oplus\Lambda^+)\otimes\su(E)).
\end{equation}
Thus, an orientation for $\det\bdelta_E$ defines an orientation for
manifold $M_\kappa^w$, since
$$
\det(d^*_{\hat A}+d^+_{\hat A}) 
=
\Lambda^{\text{max}}(H^1_{\hat A})
\otimes\Lambda^{\text{max}}(H^0_{\hat A}\oplus H^2_{\hat A})^*,
$$
We recall the following result of Donaldson:

\begin{prop}
\label{prop:BxiOrientable} \cite[Corollary 3.27]{DonOrient}
The bundle $\det\bdelta_{E}\to\sB_\kappa^w$ is topologically trivial.
\end{prop}

We now show that $\det\bsD_{\ft}$ is trivial using the fact that
$\det\bdelta_{E}$ is trivial, where $\ft=(\rho,V)$, $V=W\otimes E$, we
identify $\fg_{\ft}\cong\su(E)$ and so $w_2(\ft)\equiv w\pmod{2}$ and
$p_1(\ft)=-4\kappa$.  We shall denote
\begin{equation}
  \label{eq:PertDirac}
  D_{A,\vartheta} = D_A + \rho(\vartheta)
\quad\text{and}\quad
  D_{B,\vartheta} = D_B + \rho(\vartheta),
\end{equation}
where $D_A:C^\8(V^+)\to C^\8(V^-)$ and $D_B:C^\8(W^+)\to C^\8(W^-)$ are the
Dirac operators defined by spin connections $A$ on $V$ and $B$ on $W$,
respectively \cite[\S 2.2 \& \S 2.3]{FL2a}.

\begin{lem}
\label{lem:OrientableModuliSpace} The bundle
$\det\bsD_{\ft}\to \sC_{\ft}/S^1$ is topologically trivial.
\end{lem}

\begin{proof}
We recall that the K-theory isomorphism class of an index bundle
over a compact topological space depends only on the homotopy
class of its defining family of Fredholm operators (see, for
example, \cite[p. 69]{BoossBleecker}). Moreover, the isomorphism
class of the determinant line bundle (over a possibly non-compact
topological space) depends only on the homotopy class of the
family of Fredholm operators \cite[Lemma 6.6.1]{MorganSWNotes}.
In particular, the defining family of Fredholm operators
$\sD_{A,\Phi}$, parameterized by $[A,\Phi]\in\sC_{\ft}/S^1$, is
homotopic through $\sD_{A,t\Phi}$, $t\in[0,1]$, to the family of
Fredholm operators $\sD_{A,0} = (d^*_{\hat A}+ d^+_{\hat A})\oplus
D_{A,\vartheta}$ parameterized by $[A,\Phi]\in\sC_{\ft}/S^1$. Thus,
\begin{equation}
\label{eq:ASDOrientIsom}
\det\sD_{A,\Phi}
\cong {\det}(d^*_{\hat A}+ d^+_{\hat A})\otimes \det D_{A,\vartheta}. 
\end{equation}
Let $\det\bD_V$ be the real determinant line bundle over
$\sC_{\ft}/S^1$ associated to the family of perturbed Dirac operators,
$D_{A,\vartheta}$, where $[A,\Phi]\in\sC_{\ft}$. Let
\begin{equation}
  \label{eq:ProjPairsToConns}
 \pi_{\sB}:\sC_{\ft}/S^1\to\sB_\kappa^w, \quad [A,\Phi]\mapsto [\hat A]
\end{equation}
be the projection. Equation \eqref{eq:ASDOrientIsom} implies
there is an isomorphism $\det\bsD_{\ft} \cong \pi^*_{\sB}\det\bdelta_{E}
\otimes\det\bD_{V}$ of real determinant line bundles, so
$$
w_1(\det\bsD_{\ft})
= \pi^*_{\sB}w_1(\det\bdelta_{E}) + w_1(\det\bD_{V}).
$$
Because the Dirac operators $D_{A,\vartheta}$ have complex kernels and
cokernels, the real line bundle $\det\bD_{V}\to \sC_{\ft}/S^1$
is topologically trivial and hence $w_1(\det\bD_{V})=0$.  By
Proposition \ref{prop:BxiOrientable} we have
$w_1(\det\bdelta_{E})=0$. Combining these observations yields
$w_1(\det\bsD_{\ft}) = 0$.
\end{proof}

\subsection{Orientations of moduli spaces of PU(2) monopoles and
anti-self-dual connections}
\label{subsec:OrientASD}
We introduce an orientation for the $\PU(2)$-monopole moduli space
$\sM_{\ft}$ determined by an orientation for the moduli space
$M_\kappa^w\embed \sM_{\ft}$ of anti-self-dual connections.

An orientation for the line bundle $\det\bsD_{\ft}$ determines an
orientation for $\sM_{\ft}$.  The space $\tsC_{\ft}$ is connected, so the
quotients $\sC_{\ft}$ and $\sC_{\ft}/S^1$ are connected and a choice of
orientation for $\det\bsD_{\ft}$ is equivalent to a choice of orientation
for a fiber $\det\sD_{A,\Phi}$ over a point $[A,\Phi]$. The proof of Lemma
\ref{lem:OrientableModuliSpace} provides a method of orienting
$\det\bsD_{\ft}$ from an orientation for $\det\bdelta_{E}$, and thus from
a homology orientation and integral lift $w$ of $w_2(\ft)$, using the
isomorphism \eqref{eq:ASDOrientIsom} of real determinant lines.
Indeed, it suffices to choose an orientation for the line $\det(d^*_{\hat
  A}+ d^+_{\hat A})$ and thus an orientation for $\det\bdelta_{E}$ and
choose the orientation of $\det \bD_V$ induced from the complex
orientations of the complex kernel and cokernel of $D_{A,\vartheta}$.

To fix our conventions and notation, we outline Donaldson's method for
orienting $\det(d^*_{\hat A}+d^+_{\hat A})$, and thus $\det\bdelta_{E}$,
given a homology orientation $\Om$ and an integral lift $w$ of $w_2(\ft)$:
the detailed construction is described in \cite[\S 3]{DonOrient}. Suppose
$E\cong\underline{\CC}\oplus L$ is a Hermitian, rank-two vector bundle over
$X$, where $\underline{\CC} = X\times\CC$ and $L$ is a complex line bundle
with $c_1(L)=w$.  Then $\su(E)\cong i\underline{\RR}\oplus L$ has
$w_2(\su(E)) \equiv c_1(L) \pmod{2}$, where $\underline{\RR} = X\times\RR$.
If $d_\CC\oplus A_L$ is a reducible connection with respect to the
splitting of $E$, where $d_\CC$ is the trivial connection on $\underline{\CC}$,
and $\hat A = d_\RR\oplus A_L$ is the corresponding reducible connection on
$\su(E)$, where $d_\RR$ the trivial connection on $\underline{\RR}$. Then
the induced rolled-up deformation complex for the anti-self-dual equation
\eqref{eq:ASDOrientationSequence} splits as
\begin{equation}
  \label{eq:SplitASDRedCplx}
d_{\hat A}^*+d_{\hat A}^+ = (d^*+d^+)\oplus(d^*_{A_L}+ d^+_{A_L}),
\end{equation}
where,
\begin{gather}
\label{eq:HomologyDefOpr}
d^* + d^+: C^\8(i\Lambda^1)\to C^\8(i\Lambda^0\oplus i\Lambda^+),
\\
\label{eq:CplxLineTwistedHomologyDefOpr}
d^*_{A_L}+ d^+_{A_L}: C^\8(\Lambda^1\otimes_\RR L)\to
C^\8((\Lambda^0\oplus\Lambda^+)\otimes_\RR L).
\end{gather}
The real determinant line,
$$
\det(d^*+d^+) \cong \Lambda^{\text{max}}(H^1(X;\RR))\otimes
\Lambda^{\text{max}}(H^0(X;\RR)\oplus H^+(X;\RR))^*,
$$
is oriented by a choice of ``homology orientation'' $\Omega$
\cite[\S 3]{DonOrient}, that is, an orientation for $H^1(X;\RR)\oplus
H^{2,+}(X;\RR)$, while $H^0(X;\RR)$ is oriented by the choice of
orientation for $X$ \cite[\S 6.6]{MorganSWNotes}. The
operator $d^*_{A_L}+ d^+_{A_L}$ is complex linear, and hence the complex
orientations of its complex kernel and cokernel determine an orientation
for the real line $\det(d^*_{A_L}+d^+_{A_L})$. Thus, an orientation for
$\det(d^*_{\hat A}+d^+_{\hat A})$ is defined by the class $w$ and homology
orientation $\Omega$.

An isomorphism between any two pairs of Hermitian, rank-two complex vector
bundles $E$, $E'$ over $X$ with first Chern class $w$ can be constructed by
splicing in $|c_2(E)-c_2(E')|$ copies of $\SU(2)$ bundles over $S^4$ with
second Chern class one. Given a $\U(2)$ connection on the bundle over
$X$ with smaller second Chern class, we obtain a $\U(2)$ connection on the
other $\U(2)$ bundle by splicing in copies of the one-instanton on $S^4$.
The excision principle \cite[\S 3]{DonOrient}, \cite[\S 7.1]{DK} implies
that an orientation for one of the pair $\det\bdelta_{E}$, $\det\bdelta_{E'}$
determines an orientation for the other.

For the moduli space $M_\kappa^w$ of anti-self-dual $\SO(3)$ connections,
we let $o(\Omega,w)$ denote the orientation determined by the class $w\in
H^2(X;\ZZ)$ and corresponding split $\U(2)$ bundle, $\underline{\CC}\oplus
L$, with $c_1(L)=w$, together with a homology orientation $\Omega$.

\begin{defn}
\label{defn:ASDOrient}
Let $w\in H^2(X;\ZZ)$ be an integral lift of $w_2(\ft)$.  The orientation
$O^{\asd}(\Omega,w)$ for the line bundle $\det\bsD_{\ft}$ over
$\sC_{\ft}/S^1$, and so for the moduli space $\sM_{\ft}$, is defined by:
\begin{itemize}
\item The orientation of a fiber $\det\sD_{A,\Phi}$ over a point
  $[A,\Phi]\in\sC_{\ft}$, via isomorphism \eqref{eq:ASDOrientIsom},
\item The complex orientation for $\det D_{A,\vartheta}$, and
\item The orientation $o(\Omega,w)$ for $\det(d^*_{\hat A}+d^+_{\hat A})$.
\end{itemize}
\end{defn}

For the moduli space of anti-self-dual connections on an $\SO(3)$ bundle,
we shall need to compare orientations defined by different integral lifts
of its second Stiefel-Whitney class:

\begin{lem}
\cite[p. 283]{DK}
\label{lem:DonOrientDiffDet}
Let $X$ be a closed, oriented, Riemannian four-manifold and let $\Omega$ be
a homology orientation. If $w, w'\in H^2(X;\ZZ)$ obey $w\equiv w'\pmod{2}$, 
then
$$
o(\Omega,w') = (-1)^{\frac{1}{4}(w-w')^2}o(\Omega,w).
$$
\end{lem}

\subsection{Orientations of moduli spaces of PU(2) and Seiberg-Witten
  monopoles} 
\label{subsec:OrientSW} 
We introduce an orientation for the $\PU(2)$-monopole moduli space
$\sM_{\ft}$ determined by an orientation for a Seiberg-Witten moduli space
$M_{\fs}\embed \sM_{\ft}$.

Let $(A,\Phi)=\iota(B,\Psi) =(B\oplus B\otimes A_L,\Psi)$ be a
reducible pair in $\tsC_{\ft}$, with respect to a splitting $V=W\oplus
W\otimes L$, where $\fs=(\rho,W)$ and $\ft=(\rho,V)$ and
$\iota:\tsC_{\fs}\hookrightarrow\tsC_{\ft}$ denotes the embedding (see
Lemma 3.11 in \cite{FL2a}). Recall from
\cite[\S 3.4]{FL2a} that the deformation operator
$\sD_{\iota(B,\Psi)}$ admits a splitting $\sD_{\iota(B,\Psi)} =
\sD^t_{\iota(B,\Psi)}\oplus\sD^n_{\iota(B,\Psi)}$ into tangential and
normal components given by \cite[Equations (3.21) \& (3.22)]{FL2a}; the
splitting is $\sG_{\fs}$-equivariant with respect to the
inclusion $\sG_{\fs} \hookrightarrow \sG_{\ft}$ of automorphism groups in
\cite[Equation (3.10)]{FL2a}. Hence, we have an isomorphism of real
determinant lines,
\begin{equation}
\label{eq:ReduceOrientIsom1}
\det\sD_{\iota(B,\Psi)} \cong
\det\sD^t_{\iota(B,\Psi)} \otimes\det\sD^n_{\iota(B,\Psi)}.
\end{equation}
Furthermore, by comparing equations (2.12), (2.17), and (2.19)
with \cite[Equations (3.26) \& (3.32]{FL2a}, we see that the rolled-up
Seiberg-Witten elliptic deformation complex is identified with
the rolled-up tangential deformation complex (3.34) in
\cite{FL2a}. This identifies an orientation for the line
$\det\sD_{B,\Psi}$ with an orientation for
$\det\sD^t_{\iota(B,\Psi)}$.  Combined with the isomorphism
\eqref{eq:ReduceOrientIsom1}, this yields
\begin{equation}
\label{eq:ReduceOrientIsom} \det\sD_{\iota(B,\Psi)} \cong
\det\sD_{B,\Psi} \otimes\det\sD^n_{\iota(B,\Psi)}.
\end{equation}
The pair $\iota(B,\Psi)\in\tsC_{\ft}$ is a fixed point of the $S^1$ action
on $\tsC_{\ft}$ induced by the $S^1$ action on $V=W\oplus W\otimes L$
given by the trivial action on the factor $W$ and the action by scalar
multiplication on $L$ (see \cite[Equation (3.2)]{FL2a}). The operator
$$
\sD^n_{\iota(B,\Psi)}:
C^\8(\Lambda^1\otimes L)\oplus C^\8(W^+\otimes L)
\to
C^\8(L)\oplus C^\8(\Lambda^+\otimes L)\oplus C^\8(W^-\otimes L)
$$
is gauge equivariant and thus, because $\iota(B,\Psi)$ is a fixed point of
this $S^1$ action, is complex linear.  Hence, $\sD^n_{\iota(B,\Psi)}$ is
complex linear and the complex orientations on its complex kernel and
cokernel induce an orientation for $\det\sD^n_{\iota(B,\Psi)}$.

We recall that a homology orientation $\Omega$ defines an
orientation for $M_{\fs}$ \cite[\S 6.6]{MorganSWNotes}. As in
\cite[\S 2.3]{FL2a}, we let $\tsC_{\fs}$ denote the pre-configuration space of
pairs $(B,\Psi)$, where $\fs=(\rho,W)$, $B$ is a spin connection on $W$,
and $\Phi$ is a section of $W^+$; then $\sC_{\fs}=\tsC_{\fs}/\sG_{\fs}$ is
the configuration space, where $\sG_{\fs}\cong \Map(X,S^1)$ is the group of
\spinc automorphisms of $W$.  If $(B,0)$ is a point in $\tsC_{\fs}$ then
from
\cite[Equations (2.61) \& (2.62)]{FL2a}, the rolled-up
Seiberg-Witten elliptic deformation complex is given by
$$
\sD_{B,0}: C^\8(i\Lambda^1)\oplus C^\8(W^+) \to
C^\8(i\Lambda^0\oplus i\Lambda^+)\oplus C^\8(W^-).
$$
According to \cite[Equations (2.59), (2.60), \& (2.62)]{FL2a}, we
have
$$
\sD_{B,0}
= (d^* + d^+)\oplus D_{B,\vartheta},
$$
where $d^* + d^+$ is the operator in \eqref{eq:HomologyDefOpr} and
$D_{B,\vartheta}$ is the Dirac operator in \eqref{eq:PertDirac}. Thus,
\begin{equation}
\label{eq:SWOrientIsom}
\det\sD_{B,0}
\cong
\det(d^* + d^+)\otimes \det D_{B,\vartheta}.
\end{equation}
The determinant line bundle $\det\bsD_{\fs}$ with fibers
$\det\sD_{B,\Psi}$ is topologically trivial over
$\sC_{\fs}$, so $M_{\fs}$ is orientable and, as $\sC_{\fs}$
is connected, an orientation for the real line
$\det\sD_{B,0}$ defines an orientation for
$\det\bsD_{\fs}$. A homology orientation $\Om$
determines an orientation for $\det(d^* + d^+)$. Since the Dirac
operator $D_{B,\vartheta}$ is complex linear, the complex
orientation for its complex kernel and cokernel defines an
orientation for the real line $\det D_{B,\vartheta}$. The
product of these orientations then defines an orientation for
$\det\sD_{B,0}$ and hence for $\det\bsD_{\fs}$ and $M_{\fs}$.

\begin{defn}
\label{defn:ReducibleOrient}
The orientation $O^{\red}(\Omega,\ft,\fs)$ for the real line
$\det\sD_{A,\Phi}$, and so for the line bundle $\det\bsD_{\ft}$ and the
moduli space $\sM_{\ft}$, is defined, through the isomorphism
\eqref{eq:ReduceOrientIsom}, by:
\begin{itemize}
\item The orientation for $\det\sD_{B,\Psi}$, and thus
  $\det\bsD_{\fs}$, given by the homology orientation $\Omega$,
\item The complex orientation for $\det\sD_{\iota(B,\Psi)}^n$.
\end{itemize}
\end{defn}

\subsection{Comparison of orientations of moduli spaces of PU(2) monopoles}
\label{subsec:OrientComparison}
We now compare the different possible orientations for $\sM_{\ft}$ which we
have defined in the preceding sections.

\begin{lem}
\label{lem:OrientComparison}
Let $\ft$ be a \spinu structure on an oriented four-manifold $X$ and let
$\Om$ be a homology orientation. Suppose that $w$ is an integral lift of
$w_2(\ft)$ and that $\ft$ admits a splitting $\ft=\fs\oplus \fs\otimes L$,
for some complex line bundle $L$.  Then,
\begin{align*}
O^{\asd}(\Omega,w)
&=(-1)^{\frac{1}{4}(w-c_1(L))^2}O^{\asd}(\Omega,c_1(L)), 
\\
O^{\asd}(\Omega,c_1(L)) 
&= O^{\red}(\Omega,\ft,\fs).
\end{align*}
\end{lem}

\begin{proof}
  By Definition \ref{defn:ASDOrient}, the difference between
  $O^{\asd}(\Omega,w)$ and $O^{\asd}(\Omega,c_1(L))$ is equal to the
  difference between the orientations $o(\Omega,w)$ and $o(\Omega,c_1(L))$
  for the moduli spaces of anti-self-dual connections on $\SO(3)$ bundles
  with second Stiefel-Whitney classes $w\pmod{2}$ and $c_1(L)\pmod{2}$,
  respectively. Since $\fg_{\ft}\cong i\underline{\RR}\oplus L$ and
  $w_2(\ft)\equiv w\pmod{2}$ by hypothesis, we have $c_1(L)\equiv
  w\pmod{2}$ and so Lemma \ref{lem:DonOrientDiffDet} applies to compute the
  difference in orientations.

To see the second equality, write $\ft=(\rho,V)$ and $\fs=(\rho,W)$ and let
$$
(A,\Phi)=\iota(B,0)=(B\oplus B\otimes A_L,0)
$$
be a pair in $\tsC_{\ft}$ which is reducible with respect to the
splitting $V = W\oplus W\otimes L$ and which has a vanishing spinor
component, with $A_L=A_\La\otimes (B^{\det})^*$.  Recall from
\cite[Lemma 2.9]{FL2a} that $\hat A$ is then reducible with respect to the splitting
$\fg_{\ft}=i\ubarRR\oplus L$ and can be written as $\hat A= d_{\RR}\oplus
A_L$.  The Dirac operator $D_{A,\vartheta}$ also splits,
\begin{equation}
  \label{eq:DiracSplitting}
D_{A,\vartheta} = D_{B,\vartheta}\oplus D_{B\otimes A_L,\vartheta},
\end{equation}
where
$$
D_{B,\vartheta}: C^\8(W^+)\to C^\8(W^-)
\quad\text{and}\quad
D_{B\otimes A_L,\vartheta}: C^\8(W^+\otimes L)\to C^\8(W^-\otimes L).
$$
The isomorphism \eqref{eq:ASDOrientIsom} of determinant lines giving the
orientation $O^{\asd}(\Omega,c_1(L))$ to the line $\det\sD_{A,0}$ and the
decompositions \eqref{eq:SplitASDRedCplx} of $d_{\hat A}^* + d_{\hat A}^+$
and \eqref{eq:DiracSplitting} of $D_{A,\vartheta}$ at a reducible
connection $A$ yield the isomorphisms
\begin{equation}
  \label{eq:Orient1}
\begin{aligned}
\det\sD_{A,0}
&\cong \det(d_{\hat A}^* + d_{\hat A}^+)
\otimes\det D_{A,\vartheta}
\\
&\cong \det(d^* + d^+)
\otimes \det(d_{A_L}^* + d_{A_L}^+)
\otimes \det D_{B,\vartheta} \otimes\det D_{B\otimes A_L,\vartheta}.
\end{aligned}
\end{equation}
The operators $d_{A_L}^* + d_{A_L}^+$, $D_{B,\vartheta}$, and $D_{B\otimes
  A_L,\vartheta}$ are complex linear and thus have complex kernels and
cokernels.  By Definition \ref{defn:ASDOrient}, the orientation
$O^{\asd}(\Omega,c_1(L))$ is defined by choosing a homology orientation
$\Omega$ for $\det(d^*+d^+)$, and the complex orientation on the remaining
factors on the right-hand-side of \eqref{eq:Orient1}.

On the other hand, the isomorphisms \eqref{eq:ReduceOrientIsom}
and \eqref{eq:SWOrientIsom} of determinant lines giving the
orientation $O^{\red}(\Omega,\ft,\fs)$ to the line
$\det\sD_{A,0}$ yield the isomorphisms
\begin{equation}
  \label{eq:Orient2}
\begin{aligned}
\det\sD_{A,0} &\cong \det\sD_{B,0} \otimes
\det\sD^n_{\iota(B,0)}
\\
&\cong \det(d^*+d^+)\otimes \det D_{B,\vartheta}
\otimes
\det(d^*_{A_L}+d^+_{A_L}) \otimes
\det D_{B\otimes A_L,\vartheta}.
\end{aligned}
\end{equation}
By Definition \ref{defn:ReducibleOrient}, the orientation
$O^{\red}(\Omega,\ft,\fs)$ for $\det\sD_{A,0}$ is induced by the
isomorphism \eqref{eq:Orient2}, a choice of homology orientation $\Omega$
for $\det(d^*+d^+)$, and the complex orientation on the remaining factors
on the right-hand side of \eqref{eq:Orient2}.

The isomorphisms \eqref{eq:Orient1} and \eqref{eq:Orient2} thus yield the
same orientation of $\det \sD_{A,0}$, and therefore
$O^{\asd}(\Omega,c_1(L))=O^{\red}(\Omega,\ft,\fs)$.
\end{proof}

\subsection{Orientations of links of strata of reducible PU(2) monopoles}
\label{subsec:OrientLinkReducibles}
We shall need to compute the oriented intersections of codimension-one
submanifolds of $\sM^{*,0}_{\ft}/S^1$ with links $\bL_{\ft,\fs}$ in
$\sM^{*,0}_{\ft}/S^1$ of the strata $\iota(M_{\fs})$.  These computations
(see \S \ref{sec:DegreeZero}) are performed most naturally with a ``complex
orientation'' of the link $\bL_{\ft,\fs}$ induced from the complex
structure on the fibers of the ``virtual normal bundle'' of $M_{\fs}$ . We then
compare this orientation with the ``boundary orientation'' of
$\bL_{\ft,\fs}$ induced from an orientation of $\sM^{*,0}_{\ft}/S^1$ when
the link is oriented as a boundary of an open subspace of
$\sM^{*,0}_{\ft}/S^1$.  Our orientation conventions for the link
$\bL^{w}_{\ft,\kappa}$ of the stratum $\iota(M^w_{\ka})$ are explained in
\S \ref{subsubsec:GeomRepOnLinkASD}. We assume throughout this subsection
that there are no zero-section pairs in $M_{\fs}$.

Suppose that $Y$ is a connected, finite-dimensional, orientable manifold
with a free circle action.  We give $S^1\subset\CC$ its usual orientation.
If $\lambda_{S^1}$ is a vector in $T_yY$ which is tangent to an $S^1$ orbit
through $y\in Y$, then an orientation $\lambda_Y$ for $\det(T_yY)$ and an
orientation $\lambda_{Y/S^1}$ for $\det(T_y(Y/S^1))$ determine one another
through the convention
\begin{equation}
\label{eq:QuotientOrientationConvention}
\lambda_Y = \lambda_{S^1}\wedge\tilde\lambda_{Y/S^1},
\end{equation}
where $\tilde\lambda_{Y/S^1}\in \La^{\dim Y -1}(T_yY)$ satisfies
$\pi_*(\tilde\lambda_{Y/S^1})= \lambda_{Y/S^1}$ and
$\pi: Y\to Y/S^1$ is the projection.
In particular, orientations for $\sM_{\ft}$ and $\sM_{\ft}/S^1$ determine
one another via convention \eqref{eq:QuotientOrientationConvention}.

Recall from \cite[\S 3.5.3]{FL2a} that the ``thickened moduli space''
$\sM_{\ft}(\Xi,\fs)\subset\sC_{\ft}^{*,0}$ is a finite-dimensional
$S^1$-invariant manifold, defined by a choice of finite-rank,
$S^1$-equivariant, trivial ``stabilizing'' or ``obstruction'' bundle $\Xi$
over an open neighborhood of $\iota(M_{\fs})$ in $\sC_{\ft}^0$. Then
$M_{\fs}\embed \sM_{\ft}(\Xi,\fs)$ is a smooth submanifold with
$S^1$-equivariant normal bundle $N_{\ft}(\Xi,\fs)\to M_{\fs}$ and tubular
neighborhood defined by the image of the $S^1$-equivariant smooth
embedding,
$$
\bga:N_{\ft}(\Xi,\fs)\embed \sM_{\ft}(\Xi,\fs).
$$
An open neighborhood of $\iota(M_{\fs})$ in the moduli space $\sM_{\ft}$ is
recovered as the zero locus of an $S^1$-equivariant section
$\bvarphi$ of the $S^1$-equivariant vector bundle $\bga^*\Xi \to
N_{\ft}(\Xi,\fs)$: 
$$
\bga\left(\bvarphi^{-1}(0)\cap N_{\ft}(\Xi,\fs)\right)\subset \sM_{\ft}.
$$
The section $\bvarphi$ vanishes transversely on $N_{\ft}(\Xi,\fs) -
M_{\fs}$. As in \cite[Definition 3.22]{FL2a} we define the {\em
  link\/} of the stratum $\iota(M_{\fs})$ to be
$$
\bL_{\ft,\fs}
=
\bga\left(\bvarphi^{-1}(0)\cap \PP N_{\ft}(\Xi,\fs)\right)
\subset \sM_{\ft}/S^1,
$$
where $\PP N_{\ft}(\Xi,\fs)$ is the projectivization of the complex vector
bundle $N_{\ft}(\Xi,\fs)$. Via the diffeomorphism,
$$
\bL_{\ft,\fs} \cong \bvarphi^{-1}(0)\cap \PP N_{\ft}(\Xi,\fs),
$$
we can take the right-hand side as our model for the link, where $\bvarphi$
is a section of the complex vector bundle $\bga^*\Xi \to
\PP N_{\ft}(\Xi,\fs)$.

We now define the complex orientation for the link $\bL_{\ft,\fs}$.  The
tangent space of $\PP N_{\ft}(\Xi,\fs)$ is oriented by an orientation on
$M_{\fs}$ and the complex structure on the fibers.  To be precise, at a
point $[B,\Psi,\eta]$ in the fiber of $\PP N_{\ft}(\Xi,\fs)$ over $[B,\Psi]\in
M_{\fs}$, the inclusion of the fiber gives an exact sequence of tangent
spaces,
$$
0 \to T_{[\eta]}(\PP N_{\ft}(\Xi,\fs)|_{[B,\Psi]}) \to
T_{[B,\Psi,\eta]}\PP N_{\ft}(\Xi,\fs) \to
T_{[B,\Psi]}M_{\fs} \to 0,
$$
and thus an isomorphism of determinant lines,
\begin{align}
\label{eq:TanProjNormalSplittingIsom}
&\La^{\text{max}}(T_{[\eta]}(\PP N_{\ft}(\Xi,\fs)|_{[B,\Psi]}))
\otimes \La^{\text{max}}(T_{[B,\Psi]}M_{\fs})
\cong
\La^{\text{max}}(T_{[B,\Psi,\eta]}\PP N_{\ft}(\Xi,\fs)).
\end{align}
According to \cite[Lemma 3.23]{FL2a}, the section $\bvarphi$
vanishes transversely at any point $[B,\Psi,\eta]$ in an open
neighborhood of the zero section $M_{\fs}$ of
$N_{\ft}(\Xi,\fs)$, provided $[B,\Psi,\eta]\notin
M_{\fs}$. Thus, at
$[B,\Psi,\eta]\in \bga^{-1}(\bL_{\ft,\fs})
= \bvarphi^{-1}(0)\cap\PP N_{\ft}(\Xi,\fs)$, for $\eta\neq 0$, the differential of
$\bvarphi$ and the diffeomorphism $\bga$ induce an exact sequence,
$$
0 \to T_{\bga[B,\Psi,\eta]}\bL_{\ft,\fs}
\to
T_{[B,\Psi,\eta]}\PP N_{\ft}(\Xi,\fs)
\to
(\bga^*\Xi)_{[B,\Psi,\eta]} \to 0,
$$
since $T_{\bga[B,\Psi,\eta]}\bL_{\ft,\fs}
\cong\Ker(D\bvarphi)_{[B,\Psi,\eta]}$ and
$\Ran(D\bvarphi)_{[B,\Psi,\eta]} = (\bga^*\Xi)_{[B,\Psi,\eta]}$.  This
exact sequence and the isomorphism
\eqref{eq:TanProjNormalSplittingIsom} induce an isomorphism
\begin{equation}
\label{eq:ReducibleOrientIsomorphism}  
\begin{aligned}
\La^{\text{max}}(T_{\bga[B,\Psi,\eta]}\bL_{\ft,\fs})
&\cong 
\La^{\text{max}}(T_{[B,\Psi,\eta]}\PP
N_{\ft}(\Xi,\fs)) \otimes
(\La^{\text{max}}(\Xi_{\bga[B,\Psi,\eta]}))^*
\\
&\cong
\La^{\text{max}}(T_{[B,\Psi]}M_{\fs}) \otimes
\La^{\text{max}} (T_{[\eta]}(\PP N_{\ft}(\Xi,\fs)|_{[B,\Psi,\eta]}))
\\
&\quad\otimes
(\La^{\text{max}}(\Xi_{\bga[B,\Psi,\eta]}))^*.
\end{aligned}
\end{equation}
The fibers of the bundle $\Xi\to\sM_{\ft}(\Xi,\fs)$ are preserved under the
$S^1$ action.  The complex structure defined by this $S^1$ action
gives an orientation for $\La^{\text{max}}(\Xi_{\bga[B,\Psi,\eta]})$.

\begin{defn}
\label{defn:ComplexOrientation}
The {\em complex orientation\/} of the link $\bL_{\ft,\fs}$ is defined
through the isomorphism \eqref{eq:ReducibleOrientIsomorphism} with the
orientations of the terms on the right-hand-side of
\eqref{eq:ReducibleOrientIsomorphism} given by:
\begin{itemize}
\item The orientation of $M_{\fs}$ defined by a choice of homology
  orientation $\Omega$,
\item The complex orientation of the bundle $\bga^*\Xi\to \PP
  N_{\ft}(\Xi,\fs)$,
\item The complex orientation of the tangent space of a fiber of $\PP
  N_{\ft}(\Xi,\fs)$.
\end{itemize}
\end{defn}

Although the complex orientation given by Definition
\ref{defn:ComplexOrientation} is the natural orientation to use when
computing intersection numbers with $\bL_{\ft,\fs}$, we shall need to
orient $\bL_{\ft,\fs}$ as a boundary when using $\sM^{*,0}_{\ft}/S^1$ as a
cobordism.  We describe this procedure next.

Suppose $Z\subset Y$ is a compact submanifold of an oriented, Riemannian
manifold $Y$, with normal bundle $N\to Z$. If $\vec r$ is the
outward-pointing radial vector at a point $y$ on the boundary $\rd N$ of
the tubular neighborhood, also denoted $N$, then an orientation $\lambda_Y$
for $\det(T_yY)$ and an orientation $\lambda_{\rd N}$ for $\det(T_y(\rd
N))$ determine one another through the convention
\begin{equation}
\label{eq:BoundaryOrientation}
\lambda_{Y}=-\vec r\wedge \lambda_{\rd N},
\end{equation}
choosing the sign in equation \eqref{eq:BoundaryOrientation} so that the
link $\rd N$ is the boundary of $Y-N$.

For $[A,\Phi]\in \bL_{\ft,\fs}$, choose an outward-pointing radial vector
with respect to the thickened tubular neighborhood $N_{\ft}^{<
  \eps}(\Xi,\fs)/S^1$,
\begin{align}
  \label{eq:RadialVectorChoice}
\vec r\in T_{[A,\Phi]}(\sM_{\ft}(\Xi,\fs)/S^1)
&\cong
T_{[A,\Phi]}(N_{\ft}(\Xi,\fs)/S^1) 
\\
  \label{eq:DecompNormalRadial}
&\cong
\RR\cdot\vec r
\oplus
T_{[A,\Phi]}(N_{\ft}^\eps(\Xi,\fs)/S^1).
\end{align}
Because the section $\bvarphi$ of
$\bga^*\Xi$ vanishes transversely on both $N_{\ft}(\Xi,\fs)/S^1$ and its
$\eps$-sphere bundle, for generic $\eps$, we have isomorphisms
\begin{align}
\label{eq:DecompThickenedNormal}
T_{[A,\Phi]}(N_{\ft}(\Xi,\fs)/S^1) 
&\cong
T_{[A,\Phi]}(\sM_{\ft}/S^1) \oplus \Xi_{[A,\Phi]},
\\
\label{eq:DecompThickenedLink}
T_{[A,\Phi]}(N_{\ft}^\eps(\Xi,\fs)/S^1) 
&\cong
T_{[A,\Phi]}\bL_{\ft,\fs} \oplus \Xi_{[A,\Phi]}.
\end{align}
Through the isomorphism \eqref{eq:RadialVectorChoice},
let $\pi_{\sM/S^1}\vec r$ be the orthogonal projection of $\vec r$ onto the
subspace $T_{[A,\Phi]}(\sM_{\ft}/S^1)$ in equation
\eqref{eq:DecompThickenedNormal}.  If $\pi_{\sM/S^1}\vec r = 0$, we would
have $\vec r\in \Xi_{[A,\Phi]}$ and thus tangent to
$N_{\ft}^\eps(\Xi,\fs)/S^1$ at $[A,\Phi]$ by equation
\eqref{eq:DecompThickenedLink}, contradicting our choice of $\vec r$. Since
$\pi_{\sM/S^1}\vec r \neq 0$, a comparison of the isomorphisms
\eqref{eq:DecompNormalRadial}, \eqref{eq:DecompThickenedNormal}, and
\eqref{eq:DecompThickenedLink} yields
\begin{equation}
\label{eq:DecompMRadial}
T_{[A,\Phi]}(\sM_{\ft}/S^1)
\cong
(\pi_{\sM/S^1}\vec r)\cdot\RR\oplus T_{[A,\Phi]}\bL_{\ft,\fs}.
\end{equation}
Hence, we make the

\begin{defn}
\label{defn:LinkBoundaryOrientation}
Given an orientation $\lambda_{\sM/S^1}$ of $\sM_{\ft}/S^1$ and an
outward-pointing radial vector $\vec r$ with respect to the tubular
neighborhood $N_{\ft}^{<\eps}(\Xi,\fs)/S^1$, we define the {\em boundary
  orientation\/} $\lambda_{\rd \sM/S^1}$ of $\bL_{\ft,\fs}$ by
\begin{equation}
\label{eq:LinkBoundaryOrientation}
\lambda_{\sM/S^1}=-\pi_{\sM/S^1}\vec r\wedge \lambda_{\rd\sM/S^1}.
\end{equation}
\end{defn}

\begin{lem}
\label{lem:ReducibleOrientation} The complex orientation (Definition
\ref{defn:ComplexOrientation}) of the link $\bL_{\ft,\fs}$ agrees with the
boundary orientation (Definition \ref{defn:LinkBoundaryOrientation}) of
$\bL_{\ft,\fs}$ determined by the orientation $O^{\red}(\Omega,\ft,\fs)$
for $\sM_{\ft}/S^1$.
\end{lem}

\begin{proof}
  The orientation $O^{\red}(\Omega,\ft,\fs)$ of $\sM_{\ft}$ is
  defined through the isomorphism \eqref{eq:ReduceOrientIsom1}, using the
  orientation for $\det\sD^t_{\iota(B,\Psi)}$ (and thus the tangent space
  for $M_{\fs}$) given by the homology orientation $\Omega$, and the complex
  orientation for $\det\sD^n_{\iota(B,\Psi)}$.  {}From \cite[Equation
  (3.55)]{FL2a} we have an isomorphism $\det(\bsD^n) \cong
  \det([N_{\ft}(\Xi,\fs)]-[\Xi])$ and thus, for any $[B,\Psi]\in M_{\fs}$
  an isomorphism,
  \begin{equation}
    \label{eq:DetVirtualNormalBundle}
    \det \sD^n_{\iota(B,\Psi)} \cong
\La^{\text{max}}(N_{\ft}(\Xi,\fs)|_{[B,\Psi]}) \otimes
(\La^{\text{max}}(\Xi_{[B,\Psi]}))^*,
  \end{equation}
which preserves the orientations defined by the complex
structures of both sides.  The orientation $O^{\red}(\Omega,\ft,\fs)$ of
$\sM_{\ft}$ determines one for $\sM_{\ft}/S^1$ through convention 
\eqref{eq:QuotientOrientationConvention} and a boundary orientation for 
$\bL_{\ft,\fs}$ through convention \eqref{eq:LinkBoundaryOrientation}.

On the other hand, the complex orientation for $\bL_{\ft,\fs}$ uses,
through equation \eqref{eq:ReducibleOrientIsomorphism}, the complex
orientation for the complex projective space given by the fiber of $\PP
N_{\ft}(\Xi,\fs)$.  Comparing equation \eqref{eq:DetVirtualNormalBundle}
with equation \eqref{eq:ReducibleOrientIsomorphism} shows that the
difference between the two orientations lies in how the fibers of the
projections $\PP N_{\ft}(\Xi,\fs)\to M_{\fs}$ and $N_{\ft}(\Xi,\fs)\to
M_{\fs}$ are oriented.  The boundary orientation $\lambda_{\rd\sM/S^1}$ for
$\bL_{\ft,\fs}$ induced by $O^{\red}(\Omega,\ft,\fs)$ on $\sM_{\ft}$ begins
with the complex orientation for the fiber of the projection
$N_{\fs}(\Xi,\fs)\to M_{\fs}$, uses convention
\eqref{eq:QuotientOrientationConvention} to define an orientation for the
fiber of $N_{\fs}(\Xi,\fs)/S^1 \to M_{\fs}$, and then uses convention
\eqref{eq:BoundaryOrientation} to define an orientation for the boundary
$N^{\eps}_{\ft}(\Xi,\fs)/S^1 = \PP N_{\ft}(\Xi,\fs)$ of the bundle $N^{\le
  \eps}_{\ft}(\Xi,\fs)/S^1$.  Hence, it is enough to compare these two
methods of orienting the fibers of $\PP N_{\ft}(\Xi,\fs)$.

We denote the fibers of $N_{\ft}(\Xi,\fs)$, $N_{\ft}^\eps(\Xi,\fs)$, and
$\PP N_{\ft}(\Xi,\fs)$ by $\CC^k$, $(\CC^k)^\eps$, and
$\PP^{k-1}=(\CC^k)^\eps/S^1$, respectively.  If $\CC^k$ has a complex basis
$\{\vec r,v_1,\dots,v_{k-1}\}$, then the complex orientation of $\CC^k$ is
defined by
\begin{equation}
\begin{aligned}
  \label{eq:OrientMatch}
\lambda_{\CC^k}
&=
\vec r\wedge i\vec r\wedge v_1\wedge
iv_1\wedge\cdots\wedge v_{k-1}\wedge iv_{k-1}
\\
&= -i\vec r\wedge \vec r\wedge v_1\wedge
iv_1\wedge\cdots\wedge v_{k-1}\wedge iv_{k-1}.
\end{aligned}
\end{equation}
If $\lambda_{\PP^{k-1}}$ is the complex orientation for $\PP^{k-1}$ and
$\pi:\CC^k\less\{0\}\to\PP^{k-1}$ is the projection, then 
$$
\pi_*(v_1\wedge iv_1\wedge\cdots\wedge v_{k-1}\wedge iv_{k-1})
=
\lambda_{\PP^{k-1}}, 
$$
because $\{v_1,\dots,v_{k-1}\}$ is a complex basis for the tangent space to
$\PP^{k-1}$ at $\pi(\vec r)$, so equation \eqref{eq:OrientMatch} yields the
following relation between the complex orientations of $\CC^k$ and
$\PP^{k-1}$:
\begin{equation}
  \label{eq:CCtoBdrytoProjSpace}
  \lambda_{\CC^k}
=
-i\vec r\wedge \vec r\wedge \lambda_{\PP^{k-1}}.
\end{equation}
On the other hand, if $\CC^k/S^1$ is oriented through convention
\eqref{eq:BoundaryOrientation} by $-\vec r\wedge \lambda_{\PP^{k-1}}$, the
boundary orientation of the link $\PP^{k-1}$ is equal to its
complex orientation, $\lambda_{\PP^{k-1}}$.  By
convention \eqref{eq:QuotientOrientationConvention}, the orientation $-\vec
r\wedge \lambda_{\PP^{k-1}}$ for $\CC^k/S^1$ is induced by the
orientation $-i\vec r\wedge \vec r\wedge \pi^*\lambda_{\PP^{k-1}}$ for
$\CC^k$, which is equal to the complex orientation $\lambda_{\CC^k}$ by
equation \eqref{eq:CCtoBdrytoProjSpace}.  Hence, the complex and boundary
orientations of $\PP^{k-1}$ agree.

Therefore, the complex orientation agrees with the boundary orientation for
$\bL_{\ft,\fs}$, induced by the orientation $O^{\red}(\Omega,\ft,\fs)$
through the conventions \eqref{eq:LinkBoundaryOrientation} and
\eqref{eq:QuotientOrientationConvention}.
\end{proof}


\section{Cohomology classes on moduli spaces}
\label{sec:Cohomology} In this section we introduce cohomology
classes on the moduli space $\sM_{\ft}^{*,0}$ (see \S
\ref{subsec:CohomologyDef}) and construct geometric representatives for
these cohomology classes (see \S \ref{subsec:GeomRepr}).  The $\PU(2)$
monopole program uses the moduli space $\sM_{\ft}^{*,0}/S^1$ as a cobordism
between the link $\bL_{\ft,\kappa}^w$ of the anti-self-dual moduli space
stratum, $\iota(M^w_\kappa)\subset\sM_\ft$, and the links $\bL_{\ft,\fs}$
of the Seiberg-Witten strata, $\iota(M_\fs)\subset \sM_\ft$, giving an
equality between the pairings of the cohomology classes with these links.
The following geometric description should help motivate the constructions of
this section.

The intersection of the geometric representatives with $\sM_{\ft}^{*,0}$ is
a family of oriented one-manifolds.  The links $\bL^{w}_{\ft,\kappa}$ and
$\bL_{\ft,\fs}$ of the strata of zero-section and reducible monopoles
described in \cite[Definitions 3.7 \& 3.22]{FL2a} are oriented
hypersurfaces in $\sM_{\ft}^{*,0}/S^1$.  The intersection of these
hypersurfaces with the one-dimensional manifolds given by the intersection
of the geometric representatives is thus an oriented collection of points.
We would like to use the family of oriented one-manifolds to show that the
total signed count of the points in the intersection of the
geometric representatives with the links is zero (being an oriented
boundary).  This would give an equality between the oriented count of
points in the link of the stratum of zero-section monopoles with the
oriented count of points in the links of the strata of the reducible
monopoles.  In \S \ref{subsec:CohomOnASDLink} we show that the intersection
of the geometric representatives with the link $\bL^{w}_{\ft,\kappa}$ is a
multiple of the Donaldson invariant.  In \S \ref{sec:DegreeZero} we show
that the intersection of the geometric representatives with the links
$\bL_{\ft,\fs}$ can be expressed in terms of Seiberg-Witten invariants.
Hence, the cobordism gives a relation between these two invariants.

In practice, the above argument does not work in the simple manner
just described because $\sM_{\ft}^{*,0}/S^1$ is non-compact: the
non-compactness phenomenon responsible for the difficulty is due
to Uhlenbeck bubbling.  Geometrically, this means there can be
one-manifolds in the intersection of the geometric representatives
with one boundary on a link and the other end approaching a reducible in a
lower level of $\bar\sM_\ft/S^1$. Let 
\begin{equation}
\label{eq:NonZeroNonReducibleCompact}
\bar\sM_{\ft}^{*,0}\subset \bar\sM_{\ft}
\end{equation}
be the subspace consisting of points $[A,\Phi,\bx]$ where $(A,\Phi)$ is
neither a zero-section nor a reducible pair.  In \S
\ref{subsec:ExtendGeomRepr} we describe the intersection of the closure of
the geometric representatives with the lower Uhlenbeck levels of
$\bar\sM_{\ft}^{*,0}/S^1$ and show that for appropriate choices of geometric
representatives these intersections are empty.  Therefore, the ends of the
one-manifolds in $\sM_{\ft}^{*,0}/S^1$ do not approach the lower levels of
$\bar\sM_{\ft}^{*,0}/S^1$.

However, there may still be one-manifolds in the intersection of the geometric
representatives with ends approaching reducible monopoles in lower
Uhlenbeck levels of $\bar\sM_{\ft}/S^1$.  Theorem
\ref{thm:CompactReductionFormula} gives a relationship between the
Donaldson and Seiberg-Witten invariants when there are no reducible
monopoles in the lower levels of $\bar\sM_{\ft}$ and thus the ends of the
one-manifolds do not approach the lower levels of $\bar\sM_{\ft}/S^1$.  To
extend this argument to the case where there are reducible monopoles in the
lower levels of the compactification requires a description of
neighborhoods of the lower strata precise enough to allow the definition of
links of these strata of lower-level reducibles.  As we show in \S
\ref{sec:DegreeZero}, the geometric representatives intersect the strata of
reducible monopoles in sets of larger than expected codimension.  Thus, in
the case of the reducibles, we cannot cut down by geometric representatives
as we do with the zero-section monopoles and restrict our attention to a
generic point. Rather, we are forced to describe the entire link.  When the
reducible monopoles lie in a lower level, these links can be extremely
complicated.  A description of neighborhoods the strata of lower-level
reducibles, sufficient to define links, will be given in \cite{FL3},
\cite{FL4}.

We work with geometric representatives rather than cohomology
classes for two reasons.  First, describing the closure of a
geometric representative in a compactification appears to be
simpler than calculating the extension of a cohomology class.
Second, the topology near points in the lower levels of
$\bar\sM_{\ft}$ need not be locally finite (for example, there may
be infinitely many path-connected components).  Hence, it is not
known if $\bL^{w}_{\ft,\kappa}$ is triangulable and thus it may not
have a fundamental class to pair with the cohomology classes
described in \S \ref{subsec:CohomologyDef}.  This problem also
leads us to work in the category of smoothly stratified spaces
rather than that of piecewise-linear spaces.

\subsection{The cohomology classes}
\label{subsec:CohomologyDef}
In this subsection, we define the cohomology classes on the moduli spaces,
following the prescriptions of \cite{DonConn}, \cite{DonPoly}, \cite{DK},
\cite{KMStructure}.  These classes arise from a universal $\SO(3)$ bundle,
just as in Donaldson theory, and a universal line bundle. 

Recall that $\tsC^*_{\ft}$ denotes the subspace of pairs which are not
reducible, $\tsC^0_{\ft}$ denotes the subspace of pairs which are not
zero-section pairs, and $\tsC^{*,0}_{\ft}=\tsC^{*}_{\ft}\cap
\tsC^{0}_{\ft}$ (see \cite[\S 2.1.5]{FL2a}).  We define a universal $\SO(3)$ 
bundle:
\begin{equation}
\label{eq:UniversalSO(3)Bundle}
\FF_\ft
=
\tsC^*_{\ft}/S^1\times_{\sG_{\ft}}\fg_{\ft}
\to
\sC^*_{\ft}/S^1 \times X.
\end{equation}
The action of $\sG_{\ft}$ in 
\eqref{eq:UniversalSO(3)Bundle} is diagonal so, for $u\in\sG_{\ft}$ and
$(A,\Phi,\xi) \in \tsC^{*,0}_{\ft}\times \fg_{\ft}$, one has
$$
((A,\Phi),\xi)\mapsto (u(A,\Phi),u\xi).  
$$
We now define 
\begin{equation}
\label{eq:DefineMuMap}
\mu_{p}: H_\bullet(X;\RR)\rightarrow H^{4-\bullet}(\sC^{*}_{\ft}/S^1;\RR),
\quad \beta\mapsto -\textstyle{\frac{1}{4}}p_1(\FF_\ft)/\beta.
\end{equation}
Following \cite[Definition 5.1.11]{DK} we can also define a universal $\SO(3)$
bundle over the quotient space of $\SO(3)$ connections,
\begin{equation}
\label{eq:ConnectionUniversalSO(3)Bundle}
\FF_\kappa^w
=
\sA_{\kappa}^{w,*}\times_{\sG_\kappa^w}F
\rightarrow
\sB^{w,*}_{\ka}\times X,
\end{equation}
where $F$ is an $\SO(3)$ bundle over $X$ with $\ka=-\quarter p_1(F)$ and
$w$ is an integral lift of $w_2(F)$, and $\sG_\kappa^w$ is the group of
special unitary automorphisms of the $\U(2)$ bundle $E$ with $\su(E)=F$, so
$p_1(F)=p_1(\su(E))$, and $c_1(E)=w$.  As in \cite[Definition 5.1.11]{DK},
we define cohomology classes on $\sB_\kappa^{w,*}$ via
\begin{equation}
  \label{eq:DonaldsonMuMap}
\nu_{p}: 
H_\bullet(X;\RR)\rightarrow H^{4-\bullet}(\sB^{w,*}_{\ka};\RR),
\quad \beta \mapsto -\textstyle{\frac{1}{4}}p_1(\FF_\kappa^w)/\beta.
\end{equation}
Comparing \eqref{eq:DefineMuMap} and \eqref{eq:DonaldsonMuMap}, we see that
there must be a simple relation, which we now describe, between the
cohomology classes defined by these two $\SO(3)$ bundles.

Recall from \cite[\S 2.1.3]{FL2a} that if $F=\su(E)$ and $V=W\otimes E$, then
we have an identification of automorphism groups,
$\sG_\kappa^w\cong\sG_\ft$, and isomorphisms
\begin{equation}
  \label{eq:SO(3)ConnToSpinConn}
 \sA_\kappa^w(X)\cong\sA_\ft, \quad \hat A\mapsto A,
\quad\text{and}\quad
 \sB_\kappa^w(X)\cong\sB_\ft, \quad [\hat A]\mapsto [A].
\end{equation}
Hence, denoting $\fg_\ft=\su(E)$, we have an isomorphism of $\SO(3)$ bundles
\begin{equation}
  \label{eq:IsomUnivSO(3)Bundles}
\FF_\kappa^w
\cong
\sA_{\ft}^{*}\times_{\sG_\ft}\fg_{\ft}
\rightarrow
\sB^{*}_{\ft}\times X.
\end{equation}
Furthermore, there are natural embeddings
\begin{equation}
  \label{eq:SpinConnsIntoPairs}
  \iota:\sA_\ft\embed\tsC_\ft, \quad A\mapsto(A,0),
\quad\text{and}\quad
\iota:\sB_\ft\embed\sC_\ft, \quad[A]\mapsto[A,0]. 
\end{equation}
Using \eqref{eq:SpinConnsIntoPairs} together with the isomorphism
\eqref{eq:IsomUnivSO(3)Bundles} and the definition
\eqref{eq:UniversalSO(3)Bundle} of $\FF_{\ft}$, we see that
\begin{equation}
  \label{eq:PullBackUnivSO(3)Bundles}
(\pi_{\sB}\times \id_X)^*\FF_\kappa^w = \FF_{\ft}
\quad\text{and}\quad
(\iota\times \id_X)^*\FF_{\ft} = \FF_\kappa^w,
\end{equation}
where $\pi_{\sB}:\sC^{*}_{\ft}/S^1\to \sB^{w,*}_{\ka}$ is the restriction
of the map \eqref{eq:ProjPairsToConns} to $\sC^{*}_{\ft}/S^1$. Since
$\pi_{\sB}$ is a deformation retract, we obtain the following relation
between the cohomology classes on $\sC_{\ft}^*/S^1$ and
$\sB_{\kappa}^{w,*}$.

\begin{lem}
\label{lem:PullbackUniversalSO3Bundle}
If $\beta\in H_\bullet(X;\RR)$, then
$\pi^*_{\sB}\nu_{p}(\beta)=\mu_{p}(\beta)$ or, equivalently,
$\iota^*\mu_{p}(\beta)=\nu_{p}(\beta)$.
\end{lem}

Because $(\iota\times\id_X)^*\mu_{p} = \nu_{p}$, we shall henceforth
write $\mu_{p}$ for both $\mu_{p}$ and $\nu_{p}$.

Lastly, we define a universal complex line bundle,
\begin{equation}
\label{eq:DefineC1LineBundle} \LL_{\ft} =
\sC^{*,0}_{\ft}\times_{(S^1,\times -2)}\CC \to \sC^{*,0}_{\ft}/S^1,
\end{equation}
where the $S^1$ action defining $\LL_{\ft}$ is given, for $[A,\Phi]\in
\sC^{*,0}_{\ft}$, $e^{i\theta}\in S^1$, and $z\in\CC$, by
\begin{equation}
\label{eq:C1ClassS1Action}
\left([A,\Phi],z\right) 
\mapsto 
\left([A,e^{i\theta}\Phi], e^{2i\theta}z\right).
\end{equation}
The factor of $2$ is necessary in the action \eqref{eq:C1ClassS1Action}
because $-1\in S^1$ acts on $\tsC_{\ft}$ as $-\id_V\in\sG_t$ and thus
$-1\in S^1$ acts trivially on $\sC^{*,0}_{\ft}$.  The negative sign in the
quotient \eqref{eq:DefineC1LineBundle} indicates that the $S^1$ action is
diagonal, and is chosen to give a more convenient sign in Lemma
\ref{lem:CohomOnASDNormalSlice}.
We then define an additional cohomology class in $H^2(\sC_\ft^{*,0}/S^1;\RR)$,
\begin{equation}
\label{eq:DefineC1}
\mu_{c}
= 
c_1(\LL_{\ft}).
\end{equation}
In Lemma \ref{lem:CohomOnASDNormalSlice} we will show that the class
$\mu_{c}$ is non-trivial on the link of the subspace
$\iota(M_\kappa^w)\subset\sM_{\ft}/S^1$.  Thus, $\mu_{c}$ does not extend
over $\iota(\sB_{\ft}^*)\subset\sC_{\ft}^*/S^1$ as the restriction of an
extension to a contractible neighborhood of point
$[A,0]\in\iota(\sB_{\ft}^*)$ would have to be trivial, contradicting Lemma
\ref{lem:CohomOnASDNormalSlice}.

\subsection{Geometric representatives}
\label{subsec:GeomRepr}
To avoid having to define the link of $M_\kappa^w$ in $\sM_{\ft}$ as a
homology class, we work with geometric representatives of these cohomology
classes.  We define geometric representatives $\sV(\beta)$ and $\sW$ to
represent $\mu_{p}(\beta)$ and $\mu_{c}$ respectively.  To facilitate the
description of the intersection of the closures of $\sV(\beta)$ and $\sW$
in $\bar\sM_{\ft}$ with the lower strata, we construct geometric
representatives with certain localization properties---they are pulled back
from configuration spaces over proper open subsets $U\subset X$.  We let
$\sB^w_{\ka}(U)$ and $\sC_{\ft}(U)$ denote the quotient spaces of
connections and pairs respectively on $U\subseteq X$, where $\ka=-\quarter
p_1(\ft)$ and $w$ is an integral lift of $w_2(\ft)$.

\subsubsection{Stratified spaces}
We begin by recalling a definition of a stratified space (see
\cite[Definition 11.0.1]{MMR}) that will be sufficient for the
purposes of defining our intersection pairings.

\label{subsubsec:StratSpaces}
\begin{defn}
\label{defn:Stratifications}
\cite{GorMacPh}, \cite{Mather}, \cite{MMR} A {\em smoothly stratified
  space\/} $Z$ is a topological space with a {\em smooth stratification\/}
given by a disjoint union, $Z=Z_0\cup Z_1\cup\cdots\cup Z_n$, where the
{\em strata\/} $Z_i$ are smooth manifolds.  There is a partial ordering
among the strata, given by $Z_i<Z_j$ if $Z_i\subset \barZ_j$.  There is a
unique stratum of highest dimension, $Z_0$, such that $\barZ_0=Z$, called
the {\em top stratum}.  If $Y$, $Z$ are smoothly stratified spaces, a map
$f:Y\rightarrow Z$ is {\em smoothly stratified\/} if $f$ is a continuous
map, there are smooth stratifications of $Z$ and $Y$ such that $f$
preserves strata, and restricted to each stratum $f$ is a smooth map.  A
subspace $Y\subset Z$ is {\em smoothly stratified\/} if the inclusion is a
smoothly stratified map.
\end{defn}

\begin{rmk}
\label{rmk:StratSubspace}
If $Z$ is a smoothly stratified space and $f:Z\rightarrow \RR$
is a smoothly stratified map, that is, $f$ is a continuous
map which is smooth on each
stratum, then for generic values of $\varepsilon$, the preimage
$f^{-1}(\varepsilon)$ is a smoothly stratified subspace of $Z$.
\end{rmk}

We shall use the following definition of a geometric
representative:

\begin{defn}
\label{defn:GeomRepresentative}
\cite[p 588]{KMStructure}.
Let $Z$ be a smoothly stratified space.  A {\em geometric
representative\/} for a rational cohomology class $\mu$ of
dimension $c$ on $Z$ is a closed, smoothly stratified subspace
$\sV$ of $Z$ together with a rational coefficient $q$,
the {\em multiplicity\/}, satisfying
\begin{enumerate}
\item
The intersection $Z_0\cap \sV$ of $\sV$ with the top stratum $Z_0$ of
$Z$ has codimension $c$ in $Z_0$ and has an oriented normal
bundle.
\item
The intersection of $\sV$ with all strata of $Z$ other than the top
stratum has codimension $2$ or more in $\sV$.
\item
The pairing of $\mu$ with a homology class $h$ of dimension $c$
is obtained by choosing a smooth singular cycle representing $h$
whose intersection with all strata of $\sV$ has the codimension
$\dim Z_0-c$ in that stratum of $\sV$, and then taking $q$ times
the count (with signs) of the intersection points between the
cycle and the top stratum of $\sV$.
\end{enumerate}
\end{defn}

\begin{defn}
\label{defn:IntersectionNumber}
Let $\sV_1,\dots,\sV_n$ be geometric representatives on a compact,
smoothly stratified space $Z$ with multiplicities $q_1,\dots,q_n$.
Assume
\begin{enumerate}
\item
The sum of the
codimensions of the $\sV_i$ is equal to the dimension of the top
stratum $Z_0$ of $Z$.
\item
For every smooth stratum $Z_s$ of $Z$, the smooth
submanifolds $\sV_i\cap Z_s$ intersect transversely.
\end{enumerate}
Then dimension-counting and the definition of a geometric
representative imply that the intersection $\cap_i \sV_i$
is a finite collection of points in the top stratum $Z_0$:
$$
\sV_1\cap\dots \cap \sV_n=\{v_1,\dots,v_N\}\subset Z_0.
$$
Let $\eps_j=\pm 1$ be the sign of this intersection at $v_j$.
Then we define the {\em intersection number
of the $\sV_i$ in $Z$\/} by setting
$$
\#\left(\sV_1\cap\cdots\cap \sV_n\cap Z\right) 
=
\left(\prod_{i=1}^n q_i\right)\sum_{j=1}^N\eps_j.
$$
A cobordism between
two geometric representatives $\sV$ and $\sV'$ in $Z$ with
the same multiplicity
is a geometric representative $\sW\subset Z\times[0,1]$
which is transverse to the boundary and with
$\sW\cap Z\times\{0\}=\sV$ and $\sW\cap Z\times\{1\}=\sV'$,
with the obvious orientations of normal bundles.

The definition of intersection number does not change
if $\sV_i$ is replaced by $\sV_i'$ and there is a cobordism
between $\sV_i$ and $\sV_i'$ whose intersection with
the other geometric representatives is transverse in each
stratum.  One can see this by observing that the intersection
of the cobordism $\sW$ with the other geometric representatives will
be a collection of one-manifolds contained in $Z_0$
because the lower strata of $Z$ have codimension two.
The boundaries of these one manifolds are
the points in the two intersections
$$
\sV_1\cap \dots \cap \sV_n\quad\text{and}\quad
\sV_1\cap\dots \sV_{i-1}\cap \sV_i'\cap\dots\cap \sV_n,
$$
giving the equality of oriented intersection numbers.
\end{defn}

\begin{rmk}
In Definition \ref{defn:IntersectionNumber},
it is necessary to assume that the geometric representatives
have transverse intersection in each stratum because
we cannot assume there are perturbations of the geometric
representatives which do intersect transversely.
The definition of a smoothly-stratified space
in Definition \ref{defn:Stratifications} does not
control the topology of one strata near another.
If there is ``control data'' on a neighborhood of
one strata in another (see \cite[p. 42]{GorMacPh}, as is true
for Whitney-stratified spaces,
then such perturbations are constructed in \cite[\S 1.3]{GorMacPh}.
Instead, we will construct our geometric representatives pulling them
back from a smooth manifold where one can assume that generic choices
of the geometric representatives intersect transversely.
\end{rmk}

\subsubsection{Preliminaries for localization}
\label{subsubsec:LocalizPrelim}
By construction, our geometric representatives will be ``determined by
restriction to submanifolds'' of $X$, in the sense that they have the
following localization property:

\begin{defn}
\label{defn:DeterminedByRestriction}
Let $U\subset X$ be a submanifold.  A geometric representative $\sV$ in
$\sB^{w,*}_{\ka}$ or $\sC^{*,0}_{\ft}/S^1$ is {\em determined by
  restriction to $U\subset X$\/} if there is a geometric representative
$\sV_U$ in $\sB^{w,*}_{\ka}(U)$ or $\sC^{*,0}_{\ft}(U)/S^1$ such that
$\sV=r^{-1}_U(\sV_U)$, where $r_U$ is the map given by restricting connections
or pairs to the submanifold $U$.
\end{defn}

This localization property will allow a partial description of the
intersection of the closures of the geometric representatives in the
subspace $\bar\sM_{\ft}^{*,0}/S^1$ with the lower strata in this
compactification.  The technical issue which has complicated this
localization technique since its introduction in
\cite{DonPoly}, \cite{DK} (see \cite[p. 192]{DK}) is that there can be
pairs (connections) which are irreducible on $X$ but are reducible when
restricted to a submanifold $Y\subset X$.  The bundles over
$\sB^{w,*}_{\ka}(Y)$, whose sections define the geometric representative,
do not extend over $\sB^w_{\ka}(Y)$.  Therefore, the pullback of these sections
do not have good properties (transversality, for example) over the subspace
of connections in $\sB^{w,*}_{\ka}(X)$ which are reducible when restricted
to $Y$. When working with the moduli space $M^{w}_{\ka}$ of anti-self-dual
$\SO(3)$ connections, this problem can be overcome, if $X$ is simply
connected, by working with a tubular neighborhood $\nu(Y)$ of $Y\subset X$.
The local-to-global reducibility result of \cite[Lemma 4.3.21]{DK} implies
that any anti-self-dual connection which is irreducible on $X$ must be
irreducible on $\nu(Y)$ if $X$ is simply connected.  If $X$ is not simply
connected, there can be ``twisted reducible'' connections, (see \cite[Lemma
3.5]{FL2a} or \cite[Lemma 2.4]{KMStructure}) which are irreducible on $X$
but reducible when restricted to a tubular neighborhood.  The notion of a
``suitable open neighborhood'' of $Y$ (see Definition
\ref{defn:SuitableNeighborhood}) was introduced in \cite{KMStructure} to
deal with the problem of twisted reducibles.  Any anti-self-dual connection
which is irreducible on $X$ must be irreducible on a suitable open
neighborhood.  We can then define geometric representatives by pulling back
sections of bundles over $\sB^{w,*}_{\ka}(U(Y))$, where $U(Y)$ is a
suitable open neighborhood of $Y$.  When points $[\hat A_\alpha]$ in
this geometric representative approach a point $[\hat A_0,\bx]$ in a lower
Uhlenbeck level where the support of $\bx$ is disjoint from $U(Y)$, the
point $[\hat A_0]$ will also be in this geometric representative.  Thus,
either $[\hat A_0]$ is in the geometric representative or the support of
$\bx$ meets $U(Y)$.  Both of these conditions have high enough codimension
in the lower Uhlenbeck levels of the compactification to ensure, via
dimension-counting arguments (see
\cite[pp.  592--593]{KMStructure}), that the intersection of the geometric
representatives is compactly supported in the top level $M^{w}_{\ka}$ of
the Uhlenbeck compactification $\barM^{w}_{\ka}$.

We begin by recalling the following definition of Kronheimer and Mrowka:

\begin{defn}
\label{defn:SuitableNeighborhood}
\cite[p. 589]{KMStructure}
A smooth submanifold-with-boundary or
open set $U\subsetneqq X$ is {\em suitable\/} if
the induced map $H_1(U;\ZZ/2\ZZ)\rightarrow H_1(X;\ZZ/2\ZZ)$ is surjective.
\end{defn}

Let $Y\subset X$ be a submanifold with tubular neighborhood $\nu(Y)$. If
$D$ is a set of embedded loops generating $H_1(X;\ZZ/2\ZZ)$, which are
mutually disjoint and transversal to $Y$, then a tubular neighborhood
$\nu(Y\cup D)$ of $Y\cup D \subset X$ is a suitable neighborhood $Y$.  By
tubular neighborhood of the possibly singular space $Y\cup D$, we mean a
smoothing of the union of the tubular neighborhoods of $Y$ and of $D$.

\begin{rmk}
\label{rmk:SuitableNbhdAsdRestriction}
If $[\hat A]\in M^{w,*}_{\ka}$ then the restriction of $\hat A$ to any
suitable open neighborhood is irreducible by unique continuation
\cite[Lemma 4.3.21]{DK}. The only irreducible anti-self-dual
connections which could be reducible when restricted to an open set 
(which is not suitable) are the
twisted reducibles (see \cite[\S 3.2]{FL2a} or \cite[p. 586]{KMStructure}).
\end{rmk}

The corresponding local-to-global reducibility result for $\PU(2)$
monopoles which are not zero-section pairs \cite[Theorem 5.11]{FL1} is
stronger than that for anti-self-dual connections: 

\begin{thm}
\label{thm:LocalToGlobalReducible}
\cite[Theorem 5.10]{FL1} Suppose $(A,\Phi)$ is a solution to the perturbed
$\PU(2)$ monopole equations \cite[Equation (2.32)]{FL2a} over a connected,
oriented, smooth four-manifold $X$ with smooth Riemannian metric such that
$(A,\Phi)$ is reducible on a non-empty open subset $U\subset X$. Then
$(A,\Phi)$ is reducible on $X$ if
\begin{itemize}
\item $\Phi\not\equiv 0$ on $X$, or
\item $\Phi\equiv 0$, and $M^{w}_{\ka}$ contains no twisted reducibles
or $U$ is suitable. 
\end{itemize}
\end{thm}

Both suitable and tubular neighborhoods of submanifolds are open subsets of
$X$ and thus have codimension zero.  However, a tubular neighborhood admits
a retraction onto the submanifold while a suitable neighborhood admits a
retraction onto the union of the submanifold and a collection of loops.
Hence, these neighborhoods may, for the purposes of counting intersection
points, be thought of as having codimension equal to that of the
submanifold or to that of the union of the submanifold and some loops in
$X$.  In this sense, the suitable neighborhood of a point has smaller
codimension than a tubular neighborhood of a point.  Because the lower
strata of $\bar\sM_{\ft}$ do not have codimension as large as those of
$\barM^{w}_{\ka}$ and because the suitable neighborhoods do not have
codimension as large as tubular neighborhoods, we will require an
additional technical condition on the elements of $H_\bullet(X;\RR)$ to
ensure the intersection of the geometric representatives does not intersect
the lower levels away from the reducible pairs.

Let $Y\subset X$ be a smooth submanifold.  If $Y$ is a manifold with
boundary, the manifold structure of the configuration space
$\sB^w_{\ka}(Y)$ is described in, say, \cite[p. 262]{DonPoly}, \cite[p.
192]{DK}, \cite[\S 2(a)]{TauL2}; the corresponding slice result for
$\sC_{\ft}(Y)/S^1$ can be obtained from the slice result for
manifolds without boundary \cite[Proposition 2.12]{FL1} by taking into
account the Neumann boundary conditions as in \cite{TauL2}. Let
$r_Y:\sC_{\ft}/S^1\rightarrow \sC_{\ft}(Y)/S^1$ and
$r_Y:\sB^w_{\ka}\rightarrow \sB^w_{\ka}(Y)$ denote the restriction maps
defined by $[A,\Phi]\mapsto [A|_Y,\Phi|_Y]$ and $[\hat A]\mapsto [\hat
A|_Y]$, respectively.  We will use the same notation for the restriction
map on any domain.

We define $\sC^*_{\ft}(X,U)$ to be the quotient space of pairs on $X$ which
are irreducible when restricted to $U$, let $\sC^0_{\ft}(X,U)$ denote the
quotient space of pairs on $X$ which are not zero-section pairs when
restricted to $U$, and let $\sC^{*,0}_{\ft}(X,U) = \sC^{*}_{\ft}(X,U)\cap
\sC^{0}_{\ft}(X,U)$. The space
$\sB^{w,*}_{\ka}(X,U)$ is similarly defined.

If $[A,\Phi]\in \sM^*_{\ft}$, then Theorem \ref{thm:LocalToGlobalReducible}
implies that the restriction of the $\SO(3)$ connection $\hat A$ to
the suitable neighborhood
$\nu(Y\cup D)$ cannot be reducible; if $[A,\Phi]\in \sM_{\ft}^{*,0}$, so
we further assume $\Phi\not\equiv 0$, then Theorem
\ref{thm:LocalToGlobalReducible} implies that the restriction of the
connection $\hat A$ to $\nu(Y)$ cannot be reducible.
There is a disjoint decomposition
$$
\sM_{\ft}^* = \sM_{\ft}^{*,0} \cup \iota(M^{w,*}_{\ka}).
$$
The unique continuation result for reducible anti-self-dual $\SO(3)$
connections \cite[Lemma 4.3.21]{DK} and $\PU(2)$ monopoles (Theorem
\ref{thm:LocalToGlobalReducible}), and the preceding decomposition
and remarks yield

\begin{lem}
\label{lem:ContainmentModuliSpace}
Let $U\subset X$ be an open subset and let $Y\subset X$ be a submanifold.
Then the following inclusions hold:
$$
\pi_{\sB}(\sM_{\ft}^{*,0})\subset \sB^{w,*}_{\ka}(X,U)
\quad\text{and}\quad M^{w,*}_{\ka}\subset\sB^{w,*}_{\ka}(X,\nu(Y\cup D)),
$$
where, as in equation \eqref{eq:ProjPairsToConns},
$\pi_{\sB}:\sC_{\ft}\to\sB^w_{\ka}$  is the
projection $[A,\Phi]\mapsto [\hat A]$.
\end{lem}

We can now proceed to construct geometric representatives for the classes
$\mu_{p}(\beta)\in H^\bullet(\sM_{\ft}^{*}/S^1;\RR)$ and $\mu_{c}\in
H^\bullet(\sM_{\ft}^{*,0}/S^1;\RR)$.

\subsubsection{The geometric representatives for $\mu_{p}$}
\label{subsubsec:GeomReprMuP1}
Let $Y\subset X$ be a smooth submanifold and let $\beta=[Y]\in
H_\bullet(X;\RR)$. Let $\nu(Y\cup D)$ be a suitable open neighborhood of $Y$.
In \cite[p.  588--595]{KMStructure}, geometric representatives
for the classes $\mu_{p}(\beta)\in H^\bullet(M_\kappa^{w,*};\RR)$,
$$
r_{\nu(Y\cup D)}^{-1}(\sV(\beta))\subset \sB^{w,*}_{\ka}(X,\nu(Y\cup D)),
$$
are defined which have the property that they are determined by
$$
\sV(\beta) \subset \sB^{w,*}_{\ka}(\nu(Y\cup D)).
$$
Let $(r_{\nu(Y\cup
  D)}\pi_\sB)^{-1}(\sV(\beta))$ be the preimage of this geometric
representative with respect to the map
$$
\pi_\sB:\sC^*_{\ft}(X,\nu(Y\cup D))/S^1\to\sB^{w,*}_{\ka}(X,\nu(Y\cup D)). 
$$
The following
result is a clear consequence of the definitions
and Lemma \ref{lem:ContainmentModuliSpace}.

\begin{lem}
\label{lem:MuP1GeomRepr}
If $Y\subset X$ is a smooth submanifold, representing a class
$\beta\in H_\bullet(X;\RR)$, then
$(r_{\nu(Y\cup  D)}\pi_\sB)^{-1}(\sV(\beta))\subset \sM^*_{\ft}/S^1$
is a geometric representative for
$\mu_{p}(\beta) \in H^{4-\bullet}(\sM_{\ft}^{*}/S^1;\RR)$ and is determined
by restriction to $\nu(Y\cup D)\subset X$.
\end{lem}

Henceforth, we shall abuse notation slightly and write $\sV(\beta)$ for
$\sV(\beta)\subset \sB^{w,*}_{\ka}(\nu(Y\cup D))$, for its preimage
$r_{\nu(Y\cup D)}^{-1}(\sV(\beta))\subset \sB^{w,*}_{\ka}(X,\nu(Y\cup D))$,
and for $(r_{\nu(Y\cup D)}\pi_\sB)^{-1}(\sV(\beta)) \subset
\sC^*_{\ft}(X,\nu(Y\cup D))/S^1$.

\subsubsection{A representative for the determinant class}
\label{subsubsec:DetClassRepr}
Recall that $\sC^{*,0}_{\ft}$ is an $S^1$ bundle over
$\sC^{*,0}_{\ft}/S^1$.  Let $\nu(x)$ be a tubular neighborhood of $x$ and
let $s$ be a generic, smooth, $C^0$ bounded section of the line bundle
\begin{equation}
  \label{eq:OpenSubsetUnivLineBun}
\LL_{\ft}(\nu(x))
=
\sC^{*,0}_{\ft}(\nu(x))\times_{(S^1,\times -2)}\CC
\rightarrow \sC^{*,0}_{\ft}(\nu(x)).
\end{equation}
The action of $S^1$ in the above is the same as that in
\eqref{eq:C1ClassS1Action}, for the definition
\eqref{eq:DefineC1LineBundle} of the universal line bundle
$\LL_{\ft}\to\sC_{\ft}^{*,0}$.  The pullback
$r_{\nu(x)}^*\LL_{\ft}(\nu(x))$ is thus isomorphic to the restriction of
$\LL_{\ft}$ to $\sC^{*,0}_{\ft}(X,\nu(x))$.  We define a geometric
representative by
\begin{equation}
\label{eq:SectionMuC1}
\sW = (r_{\nu(x)}^*s)^{-1}(0) \subset \sC^{*,0}_{\ft}(X,\nu(x)).
\end{equation}
Since $\sM_{\ft}^{*,0}\subset \sC^{*,0}_{\ft}(X,\nu(x))$
by Theorem \ref{thm:LocalToGlobalReducible} and the unique
continuation theorem for the Dirac operator,
\cite[Lemma 5.12]{FL1},
the proof of the next lemma is then clear.

\begin{lem}
\label{lem:DefineC1GeomRepr}
  For a generic choice of section $s$, the zero locus $\sW$ of the section
  $r_{\nu(x)}^*s$ is a geometric representative for $\mu_{c}\in
  H^2(\sM_{\ft}^{*,0}/S^1;\RR)$ and is determined by restriction to $\nu(x)$.
\end{lem}

Let $1\in\AAA(X)$ be the element of degree zero.  
If $z=\beta_1\cdots\beta_r\in\AAA(X)$, we write
\begin{equation}
  \label{eq:AA(X)Degrees}
\delta_i = \sum_{\{p|\beta_p\in H_i(X;\RR)\}} 1
\quad\text{and}\quad
\deg(z) = \sum_{i=0}^4 (4-i)\delta_i.
\end{equation}
For monomials $z=\beta_1\cdots\beta_r$, we set
\begin{equation}
  \label{eq:DonMuMapValues}
\begin{aligned}
\mu_{p}(z) 
&= 
\mu_{p}(\beta_1)\smile\cdots\smile\mu_{p}(\beta_r), 
\\
\sV(z) 
&= 
\sV(\beta_1)\cap\cdots\cap \sV(\beta_r),
\end{aligned}
\end{equation}
and define $\mu_{p}(z)$ and $\sV(z)$ for arbitrary elements $z\in\AAA(X)$ by
$\RR$-linearity \cite[p. 595]{KMStructure}. We write
\begin{equation}
  \label{eq:C1ClassProducts}
\mu_{c}^m
=
\underbrace{\mu_{c}\smile\cdots\smile\mu_{c}}_{\text{$m$ times}}
\quad\text{and}\quad
\sW^m
=
\underbrace{\sW\cap\cdots\cap \sW}_{\text{$m$ times}},
\end{equation}
for products of the class $\mu_{c}$ and its dual $\sW$
(with the understanding that the copies of $\sW$ in the above
representative are defined with different points $x$ and different
transversely intersecting sections $s$
in Lemma \ref{lem:DefineC1GeomRepr}).

\subsection{Extension of the geometric representatives}
\label{subsec:ExtendGeomRepr}
The Uhlenbeck closure $\bar\sM_{\ft}$ of the $\PU(2)$
monopole moduli space $\sM_{\ft}$ is described in \cite[\S 2.2]{FL2a}.
The space $\bar\sM_{\ft}$ is compact \cite[Theorem 2.12]{FL2a},
\cite[Theorem 1.2]{FL1}.
We shall need to consider the following subspaces of $\bar\sM_{\ft}$:
\begin{equation}
\label{eq:PartClosedMtSubspaces}
\begin{aligned}
\bar\sM^{*}_{\ft}
&= 
\{[A,\Phi,\bx] \in \bar\sM_{\ft}: \text{$A$ is irreducible}\},
\\
\bar\sM^0_{\ft}
&=
\{[A,\Phi,\bx] \in \bar\sM_{\ft}: \Phi\not\equiv 0\},
\end{aligned}
\end{equation}
and so, as defined
in \eqref{eq:NonZeroNonReducibleCompact},
$\bar\sM_{\ft}^{*,0} = \bar\sM^{*}_{\ft}\cap \bar\sM^0_{\ft}$. We
also define
\begin{equation}
  \label{eq:MtASDChoppedOff}
\sM^{\ge\eps}_{\ft}
=
\{[A,\Phi]\in \sM_{\ft}: \|\Phi\|^2_{L^2}\ge\eps\}
\subset
\sM_{\ft}^{0},
\end{equation}
and the subspace $\sM^{*,\ge\eps}_{\ft}\subset\sM_{\ft}^{*,0}$ is defined
analogously. By \cite[Theorem 2.13]{FL2a}, the dimension of the highest stratum
$\sM_{\ft}^{*,0}$ of $\bar\sM_{\ft}$ is given by
\begin{equation}
  \label{eq:DimHighestStratum}
\dim \sM_{\ft}^{*,0}
=
d_a(\ft)+2n_a(\ft),
\end{equation}
where
\begin{equation}
  \label{eq:AsdDimDiracIndex}
\begin{aligned}
d_a(\ft) & =-2p_1(\ft)-\textstyle{\frac{3}{2}}(\chi+2\sigma),
\\
n_a(\ft) & = \textstyle{\frac{1}{4}}(p_1(\ft)+c_1(\ft)^2-\sigma).
\end{aligned}
\end{equation}
Recall that the \spinu structure $\ft_\ell$ defined in \cite[Equation
(2.44)]{FL2a} has $p_1(\ft_\ell)=p_1(\ft)+4\ell$ and so equation
\eqref{eq:DimHighestStratum} implies $\dim\sM_{\ft_\ell} =
\dim\sM_{\ft}-6\ell$. The strata of $\sM_{\ft_\ell}/S^1\times\Sym^\ell(X)$
then have codimension at least $2\ell$ relative to the top stratum. Thus we
can calculate the intersection of geometric representatives whose
intersections with the lower strata have the expected codimensions because
this ensures that (by the usual dimension-counting argument) the
intersection of the geometric representatives will be in the top stratum.

\begin{defn}
\label{defn:GeomReprClosure}
The closures of the geometric representatives, $\sV(\beta)$, $\sW$, in
$\bar\sM_{\ft}/S^1$ are denoted by $\bar\sV(\beta)$, $\bar\sW$, respectively.
For $z=\beta_1\cdots\beta_r\in\AAA(X)$, the positive generator $x\in
H_0(X;\ZZ)$, and an integer $m\ge 0$, we denote
\begin{equation}
  \label{eq:IntersectionClosureReps}
\bar\sV(z) = \bar\sV(\beta_1)\cap\cdots\cap\bar\sV(\beta_r)
\quad\text{and}\quad
\bar\sW^m = \underbrace{\bar\sW\cap\cdots\cap\bar\sW}_{\text{$m$ times}}.
\end{equation}
\end{defn}

We shall see in Lemma \ref{lem:CyclesExtension}
that these closures intersect the lower strata of
$\bar\sM^{*}_{\ft}/S^1$ in sets of the appropriate codimension,
except for $\bar\sV(x)$ where $x\in H_0(X;\ZZ)$ (see the remarks
following the proof of Lemma \ref{lem:CyclesExtension}),
and thus are geometric representatives on the compactification, away from
the zero-section and reducible monopoles.  The description of the
intersection of $\bar\sV(\beta)$ and $\bar\sW$ with the lower strata given in
this section is incomplete, as it does not give the multiplicities of all
components of these intersections.  A more complete description will be
given in \cite{FL5} using information about neighborhoods of the lower
strata in $\bar\sM^{*}_{\ft}$ obtained from gluing maps.

Recall that $\Sym^{\ell}(X)$ is a smoothly stratified space, the
strata being enumerated by partitions of $\ell\in\NN$.  For
$i=1,\dots,\ell$, let $\pi_i:X\times\cdots\times X\rightarrow X$ be
projection onto the $i$th factor.  For any subset $Y\subset X$, let
$S^\ell(Y)$ be the image of $\pi_1^{-1}(Y)\cup\cdots\cup\pi_\ell^{-1}(Y)$
in $\Sym^{\ell}(X)$ under the projection $X^\ell\rightarrow
\Sym^\ell(Y)$.  If $\Sigma\subset\Sym^{\ell}(X)$ is a smooth stratum,
we define $S_\Sigma(Y)=S^\ell(Y)\cap \Sigma$.  Let
$\pi_\Sigma: \sM_{\ft_\ell}\times\Sigma\to\Sigma$ be the projection.

On each space $\sM^*_{\ft_\ell}/S^1$ and
$\sM^{*,0}_{\ft_\ell}/S^1$ there are geometric representatives
$\sV_\ell(\beta)$ and $\sW_\ell$ defined by the same construction as
$\sV(\beta)$ and $\sW$, except using the bundles $\fg_{\ft_\ell}$ instead of
$\fg_\ft$. We write $\sV_\ell(\beta)$ and $\sW_\ell$ for the pullbacks of
these geometric representatives to $\sM_{\ft_\ell}\times\Sym^\ell(X)$.

\begin{lem}
\label{lem:CyclesExtension}
Let $\ell \ge 0$ be an integer, let
$\Sigma\subset \Sym^{\ell}(X)$ be a smooth stratum, and let $\beta\in
H_\bullet(X;\RR)$.
\begin{enumerate}
\item If $\beta$ has a smoothly embedded
representative $Y\subset X$ with a suitable neighborhood
$\nu(D\cup Y)$ and $x\in X$ has a tubular neighborhood
$\nu(x)$, then
\begin{enumerate}
\item
$\bar\sV(\beta)\cap (M^{*}_{\ft_\ell}/S^1\times \Sigma)
        \subset \sV_\ell(\beta)\cup \pi^{-1}_\Sigma(S_\Sigma(\nu(Y\cup D)))$,
\item
$\bar\sW\cap (M^{*,0}_{\ft_\ell}/S^1\times\Sigma)
        \subset \sW_\ell\cup \pi^{-1}_\Sigma(S_\Sigma(\nu(x)))$.
\end{enumerate}
\item If $\iota(M_{\fs})\subset \sM_{\ft}$ and $\beta\in H_2(X;\RR)$ is a
  two-dimensional class with $\langle c_1(\ft)-c_1(\fs),\beta\rangle\neq
  0$, and $\gamma\in H_1(X;\RR)$, then we have the following reverse
  inclusions:
\begin{enumerate}
\item
   $\iota(M_{\fs})\subset \bar\sV(\beta)$,
\item
   $\iota(M_{\fs})\subset \bar\sV(\gamma)$,
\item
   $\iota(M_{\fs})\subset\bar\sV(x)$,
\item
   $\iota(\barM^{w}_{\ka})\cup \iota(M_{\fs})\subset\bar\sW$.
\end{enumerate}
\end{enumerate}
\end{lem}

\begin{rmk}
{}From the expression for $\mu_p(\beta)$ in Corollary
\ref{cor:CohomologyOnReducibleLink}, when $\beta\in H_3(X;\RR)$, one can
see that $\mu_p(\beta)$ extends across $\iota(M_{\fs})$. Thus,
$\bar\sV(\beta)$ should be transverse to $\iota(M_{\fs})$ and so 
$\bar\sV(\beta)\cap\iota(M_{\fs})$ would be a codimension-one subset of
$\iota(M_{\fs})$ in this case. 
\end{rmk}

\begin{proof}
Here we only prove Assertion (1) .  Assertions (2a), (2b),
and (2c) will be shown in Corollary
\ref{cor:CohomologyOnReducibleLink}, while Assertion (2d) will follow
from Lemma \ref{lem:CohomOnASDNormalSlice}.

We prove Assertion (1a) about $\bar\sV(\beta)$ by
restricting pairs to the complement of the set
$\pi^{-1}_\Sigma(S_\Sigma(\nu(Y\cup D)))$.  We
assume that $[A_\alpha,\Phi_\alpha]\in \sV(\beta)$ is a sequence of
points in $\sM_{\ft}^{*}/S^1$  converging to the point
$$
[A_\8,\Phi_\8,\by]
\in
(M^{*}_{\ft_\ell}/S^1\times\Sigma)
-
\pi^{-1}_\Sigma(S_\Sigma(\nu(Y\cup D))).
$$
Given a suitable neighborhood $U=\nu(Y\cup D)$ of
$Y\subset X-\by$, we may choose a positive
constant $r$ such that
$$
U \subset X - \cup_{y\in\by}B(y,r).
$$
By the definition of Uhlenbeck convergence, $(A_\alpha,\Phi_\alpha)$
converges in the $C^\8$ topology to $(A_\8,\Phi_\8)$ on
$X-\cup_{y\in\by}B(y,r)$, modulo gauge transformations, and thus
\begin{equation}
\lim_{\alpha\to\8}[A_\alpha|_U,\Phi_\alpha|_U]
=
[A_\8|_U,\Phi_\8|_U].
\label{eq:ConvergenceLocalization}
\end{equation}
Let $\sV_Y(\beta)\subset \sB^{w,*}_{\ka}(\nu(Y\cup D))$ be the
geometric representative whose pullback defines $\sV(\beta)$.
By Lemma
\ref{lem:MuP1GeomRepr}, if $[A_\alpha,\Phi_\alpha]\in \sV(\beta)$, then
$[(A_\alpha,\Phi_\alpha)|_{U}]\in \sV_Y(\beta)$.  Because $\sV_Y(\beta)$ is
closed (see the definition in \cite[pp. 588--592]{KMStructure}),
equation \eqref{eq:ConvergenceLocalization} implies that
$[(A_\8,\Phi_\8)|_U]\in \sV_Y(\beta)$ and thus $[A_\8,\Phi_\8]\in
\sV_\ell(\beta)$.

The same argument proves Assertion (1b) concerning $\bar\sW$,
except one observes that one can replace a suitable neighborhood
$U$ of $x$ with a tubular neighborhood $\nu(x)$.
\end{proof}

Lemma \ref{lem:CyclesExtension} shows that the intersection
of $\bar\sV(\beta)$ with the lower levels of $\bar\sM_{\ft}^*/S^1$
has the same codimension as that of $\sV(\beta)$ in $\sM_{\ft}^*/S^1$,
unless $\beta\in H_0(X;\ZZ)$. This is only significant if
$z=\beta_1\cdots\beta_r$ contains both a three-dimensional
and a four-dimensional homology class.  Then the loops in the
suitable neighborhood $\nu(D\cup x)$ may intersect the three-manifold
$Y$.
In general, if $U_i$ is a suitable neighborhood of a smooth representative
of $\beta_i$ , then for any $x\in X$, the inequality
\begin{equation}
\label{eq:CoDimBound1}
\sum_{\{i:\ x\in U_i\}}(4-\dim\beta_i)\leq 5,
\end{equation}
holds \cite[Equation (2.7)]{KMStructure}.  If there are either no
three-dimensional classes among the $\beta_i$ or no
four-dimensional classes among the $\beta_i$,
the inequality \eqref{eq:CoDimBound1} can be improved to
\begin{equation}
\label{eq:CoDimBound2}
\sum_{\{i:\ x\in U_i\}}(4-\dim\beta_i)\leq 4.
\end{equation}
If $z=\beta_1\cdots\beta_r\in\AAA(X)$ and there is
a collection of suitable
neighborhoods $U_i$ of smooth representatives of $\beta_i$
satisfying \eqref{eq:CoDimBound2}, then we call $z$
{\em intersection-suitable\/}.

\begin{lem}
\label{lem:IntersectionSuitable}
If $z=\beta_1\cdots\beta_r\in\AAA(X)$ and either $\beta_i\notin H_0(X;\ZZ)$
for $i=1,\dots,r$ or $\beta_i\notin H_3(X;\RR)$ for $i=1,\dots,r$ then
$z$ is intersection-suitable.
\end{lem}

Let $z\in\AAA(X)$ and $\delta_{c}$ be a non-negative integer which
satisfy
\begin{equation}
\label{eq:defndelta_c}
\deg(z)+2\delta_{c}=d_a+2n_a-2,
\end{equation}
so for generic choices of the geometric representatives,
$\sV(z)\cap \sW^{\delta_{c}}\cap \sM_{\ft}^{*,0}/S^1$ is a collection
of one-manifolds. If
$$
[A_0,\Phi_0,\by]\in \bar\sV(z)\cap
\bar\sW^{\delta_{c}}\cap\bar\sM_{\ft}^{*,0}/S^1,
$$
where $\by\in\Sym^\ell(X)$, then by equation \eqref{eq:CoDimBound1}
and Lemma \ref{lem:CyclesExtension} we have
$$
[A_0,\Phi_0] \in \sV_{\ell}(z')\cap \sW_{\ell}^{\delta_c-j},
$$
for some $z'=\beta_{i_1}\cdots\beta_{i_q}$ and $1\leq i_1<\cdots<i_q\leq
r$, where $0\leq j\leq \delta_c$. The preceding intersection has codimension
greater than or equal to
$$
\deg(z)+2\delta_{c}-5\ell = d_a+2n_a-2-5\ell.
$$
Then, because $\dim(\sM_{\ft_\ell}/S^1)=d_a+2n_a-1-6\ell$, the
intersection 
\begin{equation}
\label{eq:LowerLevelIntersection}
\sV_{\ell}(z')\cap \sW_{\ell}^{\delta_c-j}\cap\sM_{\ft_{\ell}}/S^1
\end{equation}
has dimension less than or equal to
$$
\dim(\sM_{\ft_\ell}/S^1)-(\deg(z)+2\delta_{c}-5\ell)
=1-\ell.
$$
Hence, there could be a point in $ \bar\sV(z)\cap
\bar\sW^{\delta_{c}}\cap\bar\sM_{\ft}^{*,0}/S^1$ contained in the level
$X\times \sM^{*,0}_{\ft_1}/S^1$, where $\ell=1$. However, if $z$ is
intersection-suitable, the dimension of the intersection
\eqref{eq:LowerLevelIntersection} is $1-2\ell$,
using equation \eqref{eq:CoDimBound2}, so intersection will be empty if
$\ell>0$.  Thus, we have the following corollary to Lemma
\ref{lem:CyclesExtension}.

\begin{cor}
\label{cor:IntersectionOfSmoothLowerStrata}
Let $z\in\AAA(X)$ be intersection-suitable and let
$\delta_{c}$ be a non-negative
integer satisfying
$$
\deg(z)+2\delta_{c}=\dim(\sM_{\ft}^{*,0}/S^1)-1=d_a+2n_a-2.
$$
Then for generic choices of geometric representatives,
the intersection
$$
\bar\sV(z) \cap \bar\sW^{\delta_{c}} \cap \bar\sM^{*,0}_{\ft}/S^1,
$$
is a collection of one-dimensional manifolds,
disjoint from the lower strata of $\bar\sM^{*}_{\ft}/S^1$.
\end{cor}

\begin{rmk}
The restriction that $z\in\AAA(X)$ be intersection-suitable
is a technical one and it should be possible to remove it; we plan to
address this point in a subsequent paper.
\end{rmk}

\subsection{Geometric representatives and zero-section monopoles}
\label{subsec:CohomOnASDLink}
Our goal in this subsection is to show that the signed count of the points
in the intersection of the geometric representatives with the link of
Donaldson moduli space,
$\#(\bar\sV(z)\cap\bar\sW^{n_a-1}\cap\bL_{\ft,\kappa}^w)$, 
can be expressed in terms of the intersection number 
$\#(\bar\sV(z)\cap \bar M_{\kappa}^w)$, which defines a Donaldson invariant (at
least when $w$ is chosen so that no $\SO(3)$ bundle with second
Stiefel-Whitney class $w\pmod{2}$ admits a flat connection); the conclusion
is stated as Proposition \ref{prop:LinkOfASD}.

\subsubsection{Geometric representatives on the stratum of zero-section
  monopoles} {}From the construction of the Uhlenbeck compactifications for
$M^w_{\ka}$ and $\sM_{\ft}$, the smoothly stratified embedding
\cite[Equation (3.5)]{FL2a},
$$
\iota:M^w_{\ka}\embed \sM_{\ft}, \quad [\hat A]\mapsto
[A,0],
$$
extends (see \cite[\S 2.2]{FL2a}, \cite[\S 4]{FL1}, \cite[\S
4.4]{DK}) to smoothly stratified embedding
$$
\iota:\barM^w_{\ka}\embed \bar\sM_{\ft}, \quad [\hat A,\bx]\mapsto [A,0,\bx],
$$
where $\ka=-\quarter p_1(\ft)$ and $w$ is any integral lift of
$w_2(\ft)$.  Define
\begin{equation}
\label{eq:ZeroSectionCompactification}
\barM^{\asd}_{\ft}
=
\{[A,\Phi,\bx]\in \bar\sM_{\ft}: \Phi = 0\}
\subset
\bar\sM_{\ft}.
\end{equation}
The space $\bL^{w}_{\ft,\kappa}$ defined in
\cite[Definition 3.7]{FL2a} serves
as a link for $\barM^{\asd}_{\ft}$.  Although the definitions
imply that
$$
\iota(\barM^{w}_{\ka})\subset \barM^{\asd}_{\ft},
$$
the reverse inclusion might not be true. For example, suppose $[\hat
A_0]$ is the gauge-equivalence class of a flat connection on an $\SO(3)$
bundle $F$ over $X$ with $w_2(F)=w\pmod{2}$ and $-\frac{1}{4}p_1(F)=0$.
The Uhlenbeck compactification $\barM^{w}_{\ka}$ might not contain all
points $[\hat A_0,\bx]\in M_0^w\times\Sym^\ell(X)$ because there can be
obstructions to gluing \cite{TauIndef} onto flat connections, as the
Freed-Uhlenbeck generic metrics theorem does not guarantee they are smooth
points of their moduli spaces \cite{DK}, \cite{FU}.  However, there could be
a sequence of points $[A_\alpha,\Phi_\alpha]\in \sM_{\ft}$ converging to
$[A_0,0,\bx]$ in the Uhlenbeck topology but no sequence $[\hat A_\alpha']\in
M^{w}_{\ka}$ also converging to $[\hat A_0,\bx]$.  

\begin{defn}
\label{defn:Good}
A class $v\in H^2(X;\ZZ/2\ZZ)$ is {\em good\/} if no integral lift of $v$
is torsion.
\end{defn}

Observe that a class $v\in H^2(X;\ZZ/2\ZZ)$ is good if and only if no line
bundle $L$ over $X$ with first Chern class $c_1(L)\equiv v\pmod{2}$ admits
a flat connection or if and only if no $\SO(3)$ bundle over $X$ with
second Stiefel-Whitney class $v$ admits a flat connection. Thus,
Lemma 3.2 in \cite{FL2a} gives a criterion for $v$ to be good.

Hence, if $w_2(\ft)$ is good, there are no flat $\SO(3)$ connections in
$\bar\sM_{\ft}$ and, when there no obstructions to gluing (for example,
when the metric on $X$ is generic in the sense of \cite{DK,FU}),
it follows from Taubes' gluing theorem for anti-self-dual $\SO(3)$
connections \cite{TauSelfDual}, \cite{TauStable} that
$$
\barM^{\asd}_{\ft}\subset \iota(\barM^{w}_{\ka}). 
$$
The preceding discussion yields

\begin{lem}
\label{lem:ASDClosureEquality}
Let $\ft$ be a \spinu structure on a closed, oriented four-manifold $X$
with generic metric, $b_2^+(X)>0$ and $w_2(\ft)\equiv
w\pmod{2}$, for $w\in H^2(X;\ZZ)$. If $w\pmod{2}$ is good, then
\begin{equation}
  \label{eq:ASDClosureEquality}
\barM^{\asd}_{\ft} = \iota(\barM^{w}_{\ka}).  
\end{equation}
\end{lem}

The constraint that $w\pmod{2}$ is good is also used to separate the strata
of zero-section $\PU(2)$ monopoles from the strata of reducible monopoles,
so that the moduli space of $\PU(2)$ monopoles gives a smooth cobordism
between their links. Therefore, when $w\pmod{2}$ is good, equation
\eqref{eq:ASDClosureEquality} holds, we have a disjoint union
$$
\bar\sM_{\ft} = \bar\sM_{\ft}^0\cup \iota(\bar M_\kappa^w),
$$
and $\bL_{\ft,\kappa}^w$ is a link of $\iota(\bar
M_\kappa^w)\subset\bar\sM_{\ft}/S^1$.

\subsubsection{A definition of the Donaldson invariants}
\label{subsubsec:DefnDInvariants}
Fix $w\in H^2(X;\ZZ)$ and let $z\in\AAA(X)$ be a monomial whose degree
$\deg(z)$ satisfies equation \eqref{eq:Mod8}. Then let
$\ka\in\frac{1}{4}\ZZ$ be defined by
\begin{equation}
  \label{eq:Kappa}
\deg(z)=8\ka -\textstyle{\frac{3}{2}}\left(\chi+\sigma\right).
\end{equation}
Let $\tilde X = X\#\overline{\CC\PP}^2$ denote the blow-up of $X$ and let
$e\in H_2(\tilde X;\ZZ)$ be the exceptional class and let $\PD[e]$ be its
Poincar\'e dual. Since $-4\kappa=w^2\pmod{4}$ by equations \eqref{eq:Mod8}
and \eqref{eq:Kappa} and thus $-4(\kappa+1/4)=(w+\PD[e])^2\pmod{4}$, we can
find an $\SO(3)$ bundle $\tilde F\to \tilde X$ with $p_1(\tilde
F)=-4(\kappa+1/4)$ and $w_2(\tilde F)=(w+\PD[e])\pmod{2}$
\cite[Theorem 1.4.20]{GompfStipsicz}. We can therefore define
$M_{\kappa+1/4}^{w+\PD[e]}(\tilde X)$, the moduli space of anti-self-dual
$\SO(3)$ connections on $\tilde F$.  Then the {\em Donaldson invariant\/} is
defined by \cite[p.  594]{KMStructure}
\begin{equation}
D_{X}^{w}(z)
=
\#\left(\bar\sV(ze)\cap \barM^{w+\PD[e]}_{\ka+1/4}(\tilde X)\right),
\label{eq:DefineDInvarBlowUp}
\end{equation}
where the moduli space $M^{w+\PD[e]}_{\ka+1/4}(\tilde X)$ is given the
orientation $o(\Omega,w+\PD[e])$  and on the right-hand side, we consider
$z$ to be a monomial in $\AAA(\tilde X)$ via the inclusion 
$$
H_2(X;\RR)\subset H_2(\tilde X;\RR)=H_2(X;\RR)\oplus\RR[e]
\quad\text{and}\quad
\AAA(X)\subset \AAA(\tilde X).  
$$
The Donaldson invariant is independent of the choice of
generic geometric representatives and, when $b_2^+(X)>1$, independent of the
metric. 

If $M^{w}_{\ka}(X)$ is given the orientation $o(\Omega,w)$,
a well-known special case of the blow-up formula \cite[Lemma
3.13]{FSStructure},  \cite[Theorem 6.9]{KotschickBPlus1} implies that
\begin{equation}
  \label{eq:KotschickBlowUpFormula}
\#\left(\bar\sV(z)\cap \barM^{w}_{\ka}(X)\right)
=
\#\left(\bar\sV(ze)\cap \barM^{w+\PD[e]}_{\ka+1/4}(\tilde X)\right),
\end{equation}
when the intersection number on the left is well-defined, for example, when
$w\pmod{2}$ is good in the sense of Definition \ref{defn:Good}. However,
the blow-up 
trick \cite{MorganMrowkaPoly} ensures that the intersection number on the
right is well-defined for arbitrary $w\in H^2(X;\ZZ)$, since $w+\PD[e]\in
H^2(\tilde X;\ZZ)$ is good.

In \cite[p. 585]{KMStructure}, Kronheimer and Mrowka require that
$b_2^+(X)-b_1(X)$ be odd. However, for the purposes of defining the Donaldson
invariants for a closed four-manifold, with $b_1(X)$ possibly non-zero, one
can have non-trivial Donaldson invariants when $b_2^+(X)-b_1(X)$ is even as
they point out in \cite[p. 595]{KMStructure}. The reason for the constraint
is that their structure theorem is only stated for the case $b_1(X) = 0$;
the invariants become more difficult to compute when $b_1(X)>0$. If $b_1(X)
= 0$, then the Donaldson invariants are necessarily trivial unless $b_2^+(X)$
is odd (and the moduli spaces of anti-self-dual connections are even
dimensional).

When $b_2^+(X)=1$, the Donaldson invariant \eqref{eq:DefineDInvarBlowUp}
depends on the ``chamber'' in $H^2(\tilde X;\RR)$ defined by metric $\tilde
g$ on $\tilde X$.  Specifically, if $\omega(\tilde g)$ is the unique
unit-length harmonic two-form which is self-dual with respect to $\tilde
g$---the {\em period point\/} for $\tilde g$---and lies in the positive cone
(determined by the homology orientation) of $H^2(\tilde X;\RR)$, then the
intersection number on the right-hand-side of equation
\eqref{eq:DefineDInvarBlowUp} changes whenever
the sign of $\omega(\tilde g)\smile \alpha$ changes for some $\alpha\in
H^2(\tilde X;\ZZ)$ satisfying
\begin{equation}
\label{eq:WallCondition}
\alpha\equiv w+\PD[e]\pmod 2
\quad
\text{and}\quad
\alpha^2=-4(\ka+1/4)+4\ell,
\quad\text{for some $\ell\in\NN$}.
\end{equation}
(The classes $\alpha$ correspond to split $\SO(3)$ bundles over $\tilde X$,
namely $\underline{\RR}\oplus L$ with $c_1(L)=\alpha$, so they have first
Pontrjagin number $\alpha^2$ and second Stiefel-Whitney class
$\alpha\pmod{2}$; see \cite{KotschickBPlus1}, \cite{KotschickMorgan} for
further explanation.)  For any $\alpha\in H^2(\tilde X;\RR)$ which is {\em
not torsion\/}, the subset
$$
\{h\in H^2(\tilde X;\RR): h\smile h>0
\text{ and }h\smile\alpha = 0\}
$$
of the positive cone of $H^2(\tilde X;\RR)$
is an {\em $\alpha$-wall\/}.  If $\alpha$ obeys condition
\eqref{eq:WallCondition}, then $\alpha$ is non-torsion since $w+\PD[e]\pmod
2$ is good and the resulting $\alpha$-wall is called a {\em
$(w+\PD[e],-4\ka-1)$-wall\/}. The connected components of the complement in
the positive cone of $H^2(\tilde X;\RR)$ of the union of
$(w+\PD[e],-4\ka-1)$-walls are called {\em
$(w+\PD[e],-4\ka-1)$-chambers\/}. Hence, the intersection pairing in
definition \eqref{eq:DefineDInvarBlowUp} changes if the period point
$\omega(\tilde g)$ moves from one $(w+\PD[e],-4\ka-1)$-chamber to another.

We now discuss how a choice of a metric $g$ on $X$
determines a chamber in the positive cone of $H^2(\tilde X;\RR) \cong
\RR[e]\oplus H^2(X;\RR)$.
Assume first that $w\pmod 2$ is good in the sense of Definition
\ref{defn:Good}. Therefore, if $\beta\in H^2(X;\ZZ)$ satisfies
\begin{equation}
\label{eq:UnblownUpWallCondition}
\beta\equiv w\pmod 2
\quad
\text{and}\quad
\beta^2=-4\kappa+4\ell,
\quad\text{for some $\ell\in\NN$},
\end{equation}
then $\beta$ is non-torsion and thus defines a $(w,-4\ka)$-wall in
$H^2(X;\RR)$. Moreover, $\alpha=\beta+\PD[e]$ satisfies condition
\eqref{eq:WallCondition} and defines a $(w+\PD[e],-4\ka-1)$-wall in
$H^2(\tilde X;\RR)$. This establishes an inclusion of $(w,-4\ka)$-chambers in
$H^2(X;\RR)$ into {\em related $(w+\PD[e],-4\ka-1)$-chambers\/} in
$H^2(\tilde X;\RR)$.

If $g$ is a generic metric on $X$, then the period point $\omega(g)\in
H^2(X;\RR)$ does not lie on any wall and there is a unique chamber in the
positive cone of $H^2(X;\RR)$ which contains $\omega(g)$ \cite[Corollary
4.3.15]{DK}. 

If the metric $\tilde g$ on $\tilde X$ is constructed by splicing together
a generic metric $g$ on $X$ and the Fubini-Study metric on
$\overline{\CC\PP}^2$ along a ``long neck'', then $\omega(\tilde g)$
converges to $\omega(g)$ in $L^2$ as the length of the neck tends to
infinity, viewing both $\omega(\tilde g)$ and $\omega(g)$ as elements of
$H^2(\tilde X;\RR)$. Thus, $\om(\tilde g)$ will lie in the chamber in the
positive cone of $H^2(\tilde X;\RR)$ related to the chamber in the positive
cone of $H^2(X;\RR)$ containing $\om(g)$.  Thus, when $w\pmod 2$ is good,
the intersection pairing \eqref{eq:DefineDInvarBlowUp} defining the
invariant $D_X^w(z)$ is defined with respect to the chamber in $H^2(\tilde
X;\RR)$ related to the chamber in $H^2(X;\RR)$ determined by the period
point $\omega(g)$.

If $w\pmod 2$ is not good, one can have torsion classes $\beta \in
H^2(X;\ZZ)$ satisfying condition \eqref{eq:UnblownUpWallCondition} and thus
$\alpha=\beta+\PD[e]\in H^2(\tilde X;\ZZ)$ satisfying
\eqref{eq:WallCondition}. The corresponding $(w+\PD[e],-4\ka-1)$-wall is
given by
$$
\{h\in H^2(\tilde X;\RR): h\smile h>0
\text{ and }h\smile\PD[e] = 0\}.
$$
Hence, the period point $\omega(g)\in H^2(X;\RR)$ for any metric $g$ on $X$
lies in this $\PD[e]$-wall, since $H^2(\tilde X;\RR)\cong \RR[e]\oplus
H^2(X;\RR)$, and the fact that $\omega(\tilde g)$ converges to $\omega(g)$
as the length of the neck tends to infinity does not determine the chamber
of $\omega(\tilde g)$ without a delicate analysis of the sign of
$\omega(\tilde g)\smile\PD[e]$ (see \cite{Yang}).  We plan to address the
case $b_2^+(X)=1$ elsewhere and so in the present article, if $b_2^+(X)=1$,
we only consider the dependence of the invariants $D_X^w(z)$ on the chamber
in $H^2(X;\RR)$ when $w$ is good.

\subsubsection{Geometric representatives on the link of the stratum of
zero-section monopoles}
\label{subsubsec:GeomRepOnLinkASD}
We now turn to the arguments leading to a proof of Proposition
\ref{prop:LinkOfASD}, which expresses the intersection number
\begin{equation}
\label{eq:ASDLinkIntersection} 
\#\left( \bar\sV(z)\cap\bar\sW^{n_a-1}\cap \bL^{w}_{\ft,\kappa}\right)
\end{equation}
in terms of
$$
\#\left(\bar\sV(z)\cap \iota(\bar M_{\kappa}^w)\right)
=
\#(\bar\sV(z)\cap \bar M_{\kappa}^w),
$$
which is equal to the Donaldson invariant $D_X^w(z)$ when $w$ is good.

Note that by the construction of the geometric representatives and
definition of the Donaldson invariants \cite{DK}, \cite[\S
2]{KMStructure}, one has
$$
\bar\sV(z)\cap \barM^w_{\ka}
=
\sV(z)\cap M^w_{\ka}.
$$
That is, the intersection is contained in the top stratum $M_{\kappa}^w$ of
the compactification $\bar M_\kappa^w$. Therefore,
to calculate the intersection number \eqref{eq:ASDLinkIntersection} it will
be enough to examine small neighborhoods of the points in the intersection
$\sV(z)\cap \iota(M^w_{\ka})$.  Such neighborhoods are described by
Kuranishi models, which we now describe.

Suppose $[A,0]\in \sV(z)\cap \iota(M^w_{\ka})$.  Applying the Kuranishi
method to describe the zero locus of the $\PU(2)$ monopole equations
using \cite[Corollary 3.6]{FL2a}, we obtain a
smooth $S^1$-equivariant embedding
\begin{equation}
\label{eq:ASDKuranishiEmbedding}
\bga_A:\sO_A\subset T_{[\hat A]}M^w_{\ka}\oplus \Ker D_{A,\vartheta}
\to
(A,0) + \Ker d_{A,0}^*\subset\tsC_{\ft},
\end{equation}
of a precompact, open $S^1$-invariant neighborhood $\sO_A$ of the origin
with image $\bga_A(\sO_A/S^1)\subset \sC_{\ft}/S^1$, where $S^1$ acts on the
domain by scalar multiplication on $\Ker D_{A,\vartheta}$, and a smooth
$S^1$-equivariant map
\begin{equation}
\label{eq:ASDKuranishiMap}
\bvarphi_A: \sO_A\subset T_{[\hat A]}M^w_{\ka}\oplus \Ker D_{A,\vartheta}
\to
\Coker D_{A,\vartheta},
\end{equation}
such that
\begin{align}
\label{eq:ImageKurZeroLocusToPU2Mod}
\bga_A\left((\bvarphi_A^{-1}(0)\cap\sO_A)/S^1\right)
&=
\sM_{\ft}/S^1\cap\bga_A(\sO_A/S^1),
\\
\label{eq:ImageKurZeroLocusToASDMod}
\bga_A\left((T_{[\hat A]}M^w_{\ka}\oplus\{0\})\cap\sO_A\right)
&=
\iota(M^w_{\ka})\cap\bga_A(\sO_A),
\end{align}
are open neighborhoods of the point $[A,0]$ in $\sM_{\ft}/S^1$ and
$\iota(M^w_{\ka})$, respectively. Because points in
$\sM_{\ft}^{*,0}$ are regular, the map $\bvarphi_A$ vanishes
transversely on
$$
\sO_A-(T_{[\hat A]}M^w_{\ka}\oplus\{0\}).
$$
(Compare the proof of Assertion (4) in \cite[Theorem 3.21]{FL2a}.) For
convenience, we set
\begin{equation}
  \label{eq:sZA}
  \sZ_A = \bvarphi_A^{-1}(0)\cap\sO_A.
\end{equation}
Equation \eqref{eq:ImageKurZeroLocusToPU2Mod} implies
that the $S^1$-equivariant embedding $\bga_A$ descends to a homeomorphism
from a neighborhood of the origin onto a neighborhood of $[A,0]$,
\begin{equation}
  \label{eq:KurHomeoZeroLocusToPU2Mod}
\bga_A:(\bvarphi_A^{-1}(0)\cap\sO_A)/S^1 
\cong
\sM_{\ft}/S^1\cap\bga_A(\sO_A/S^1),
\end{equation}
which restricts to a diffeomorphism on a smooth stratum.

In \cite[\S 3.2]{FL2a}, we constructed the link $\bL^w_{\ft,\ka}$
using the $S^1$-invariant ``distance function'',
\begin{equation}
\label{eq:DefineNormFunction}
\bl: \sC_{\ft}\to [0,\infty),\quad
[A,\Phi]\mapsto \|\Phi\|_{L^2}^2.
\end{equation}
The function $\bl$ extends continuously over $\bar\sM_{\ft}$
if we set $\bl([A,\Phi,\bx])=\|\Phi\|_{L^2}^2$. For generic, positive,
small $\eps$ we have \cite[Definition 3.7]{FL2a}
\begin{equation}
  \label{eq:ASDLink}
\bL^{w,\eps}_{\ft,\ka}= \bl^{-1}(\eps)\cap\bar\sM_{\ft}/S^1,
\end{equation}
and denote $\bL^{w,\eps}_{\ft,\ka}$ by $\bL^{w}_{\ft,\ka}$ when the value
of $\eps$ is not relevant.  To compute the pairing
\eqref{eq:ASDLinkIntersection} with $\bL^{w,\eps}_{\ft,\ka}$, we must first
describe $\bar\sV(z)$ in a neighborhood of $\iota(M^w_{\ka})$ in $\sM_{\ft}$.

\begin{lem}
\label{lem:DeformingV}
Let $\ft$ be a \spinu structure over a four-manifold $X$ with
$w_2(\ft)\equiv w\pmod{2}$, where $w\in H^2(X;\ZZ)$ and $w\pmod{2}$ is
good.  Suppose that $\deg(z)\ge \dim M^w_{\ka}$ and $z$ is
intersection-suitable, as defined before Lemma
\ref{lem:IntersectionSuitable}. Denote $\bar\sV(z)\cap \barM^w_{\ka}
=\{[\hat A_i]_{i=1}^N\}$. Then for each $[\hat A]\in 
\bar\sV(z)\cap \barM^w_{\ka}$, there is an open neighborhood
$\sO_A'\subset\sO_A$ of the origin in $T_{[\hat A]}M^w_{\ka}\oplus \Ker
D_{A,\vartheta}$, where $\sO_A$ is the open neighborhood defining the
Kuranishi model \eqref{eq:ASDKuranishiEmbedding}, such that:
\begin{enumerate}
\item
There is a smooth, $S^1$-invariant map,
$$
f_A: \sO_A'\cap(\{0\}\oplus\Ker D_{A,\vartheta}) \to T_{[\hat A]}M^w_\ka,
$$
with $f_A(0)=0$ and $(Df_A)_0=0$ such that
$$
\bga_{A}^{-1}(\sV(z))\cap \sO_A' = \{(f_A(\phi),\phi): \phi\in
\sO_A'\cap (\{0\}\oplus\Ker D_{A,\vartheta})\}.
$$
\item
There is a positive constant $\eps_0$ such that for all $\eps<\eps_0$,
$$
\bar\sV(z)\cap \bL^{w,\eps}_{\ft,\ka} \subset
\bigcup_{i=1}^N \bga_{A_i}\left( \sO_{A_i}'\right)/S^1,
$$
where the union on the right above is disjoint.
\item
For each $\eps\in (0,\eps_0)$ with $\eps_0$ as in (2), there is a positive
constant $\delta$ such that all 
$(a,\phi)\in\bga_A^{-1}\left(\sV(z)\cap\bL^{w,\eps}_{\ft,\ka}\right)\cap\sO_A'$
satisfy $\|\phi\|_{L^2}^2>\delta$.
\end{enumerate}
\end{lem}

\begin{proof}
Consider $\sV(z)$ as a smooth submanifold of $\sB^{w,*}_{\ka}$.  If
$\pi_{\sB}:\sC^*_{\ft}\to\sB^{w,*}_{\ka}$ is the projection, then the
composition $\pi_{\sB}\circ \bga_A$ is a smooth map from $\sO_A$ to
$\sB^{w,*}_{\ka}$.  The manifolds $M^w_{\ka}$ and $\sV(z)$ intersect
transversely in $\sB^{w,*}_\ka$ at $[\hat A]=\pi_{\sB}\circ\bga_A(0,0)$.
The restriction of $\bga_A$ to
$$
\sO_A\cap\left(T_{[\hat A]}M^w_{\ka}\oplus\{0\}\right)
$$
is an embedding onto an open neighborhood of $[\hat A]$ in $M^w_{\ka}$, so
the composition $\pi_{\sB}\circ\bga_A$ is transverse to $\sV(z)$ at the
origin.  Thus, restricted to a sufficiently small open neighborhood
$\sO_A''\subset \sO_A$ of the origin $(0,0)$ in 
$T_{[\hat A]}M^w_{\ka}\oplus \Ker D_{A,\vartheta}$.
the map $\pi_{\sB}\circ\bga_A$ is
transverse to $\sV(z)$ in $\sB^{w,*}_{\ka}$.  Hence
$\bga_A^{-1}(\sV(z))\cap\sO_A''$ is a smooth manifold.
By shrinking the neighborhood $\sO_A''$ we can assume that
$\bga_A^{-1}\left(\sV(z)\right)\cap \sO_A''$ and $T_{[\hat
A]}M^w_{\ka}\oplus\{0\}$ intersect only at the origin, 
since $\sV(z)\cap M_\kappa^w \cap\bga_A(\sO''_A)= [\hat A]
=\pi_{\sB}\circ\bgamma_A(0,0)$. We now prove that 
\begin{equation}
\label{eq:TangentV}
T_{(0,0)}
\left(\bga_A^{-1}(\sV(z))\cap \sO_A''\right)
=
\{0\}\oplus \Ker D_{A,\vartheta}.
\end{equation}
First, note that because the derivative of $\bga_A$ at the origin is the
inclusion of $T_{[\hat A]}M^w_{\ka}\oplus \Ker D_{A,\vartheta}$ into $\Ker
d^{0,*}_{A,0}$ by construction of the Kuranishi model, we have the
inclusion:
$$
\{0\}\oplus \Ker D_{A,\vartheta}
\subset
\Ker (D(\pi_{\sB}\circ\bga_A))_{(0,0)}
\subset
T_{(0,0)}\left(\bga_A^{-1}(\sV(z))\cap \sO_A''\right).
$$
Because $(D(\pi_{\sB}\circ\bga_A))_{(0,0)}$ maps $T_{[\hat A]}
M^{w,*}_{\ka}\oplus \{0\}$ onto the normal bundle of $\sV(z)$ in
$\sB^{w,*}_{\ka}$, the above inclusion is an equality.  Equation
\eqref{eq:TangentV} implies that if
$$
\pi_{K,A}: T_{[\hat A]}M^w_{\ka}\oplus \Ker D_{A,\vartheta} 
\to \Ker D_{A,\vartheta}
$$
is the projection onto $\Ker D_{A,\vartheta}$ then the derivative of the
restriction of $\pi_{K,A}$,
\begin{equation}
\label{eq:RestrictedpiKA}
\pi_{K,A}:\bga_A^{-1}(\sV(z))\cap \sO_A''
\to \Ker D_{A,\vartheta}
\end{equation}
at the origin is an isomorphism.  Therefore, for small enough $\sO_A''$
the map \eqref{eq:RestrictedpiKA}
is a diffeomorphism onto a neighborhood $\sO_{K,A}$ of the
origin in $\Ker D_{A,\vartheta}$.  If
$\sO_A'=\sO_A''\cap\pi_{K,A}^{-1}(\sO_{K,A})$, then $\sO_A'\cap
\{0\}\oplus\Ker D_{A,\vartheta}\subset\sO_{K,A}$ and $\pi_{K,A}$ restricted
to $\bga_A^{-1}(\sV(z))\cap \sO_A'$ is still a diffeomorphism onto
$\sO_{K,A}$ with inverse as described.

By shrinking the sets $\sO_A'$, we can assume the images $\bga_A(\sO_A')$ are
disjoint, proving the final statement in Assertion (2).  Suppose
$\eps_\alpha$ is a sequence of positive numbers converging to zero.  If
Assertion (2) were not true, there would be a sequence
$$
\{[A_\alpha,\Phi_\alpha,\bx_\alpha]\}_{\alpha=1}^\infty
\subset \bar\sV(z)\cap \bar\sM_{\ft},
$$
satisfying $\bl([A_\alpha,\Phi_\alpha,\bx_\alpha])=\eps_\alpha$ and
$[A_\alpha,\Phi_\alpha,\bx_\alpha]$ not in any $\bga_{A_i}(\sO_{A_i}')$.
Since $\bar\sV(z)\cap\bar\sM_{\ft}$ is compact, there would be a
convergent subsequence, also denoted
$\{[A_\alpha,\Phi_\alpha,\bx_\alpha]\}_{\alpha=1}^\infty$, converging to
$[A_\8,\Phi_\8,\bx_\8]$.  Because $\bl$ is continuous on $\bar\sM_{\ft}$,
we would have $\bl([A_\8,\Phi_\8,\bx_\8])=0$, so
$[A_\8,\Phi_\8,\bx_\8]\in\iota(\barM^w_{\ka})$ by Lemma
\ref{lem:ASDClosureEquality}, with $\Phi_\8=0$. Then, 
$$
[A_\8,0,\bx_\8]\in\bar\sV(z)\cap\iota(\barM^w_\ka)
=
\sV(z)\cap\iota(M^w_\ka),
$$ 
and thus $\bx=\emptyset$ and $[\hat A_\8]\in \{[\hat A_1,\dots,\hat
A_N]\}$.  The images $\bga_{A_i}(\sO_{A_i}')$ contain open neighborhoods of
the points $\{[\hat A_1],\dots,[\hat A_N]\}$ in $\bar\sM_{\ft}$, so for
large enough $\alpha$ the sequence must lie in the union of these images,
contradicting the assumption that $[A_\alpha,\Phi_\alpha,\bx_\alpha]$
is not in any $\bga_{A_i}(\sO_{A_i}')$.  This proves Assertion (2).

We use contradiction to prove Assertion (3): if it were not true, then
there would be a sequence $\{(a_\alpha,\phi_\alpha)\}_{\alpha=1}^\8$ in
$\bga_A^{-1}(\sV(z))\cap\sO_A'$ with
$\bl(\bga_A(a_\alpha,\phi_\alpha))=\eps$ and $\phi_\alpha\in \sO_{K,A}$
with $\lim_{\alpha\to\8}\|\phi_\alpha\|_{L^2}^2=0$.  By Assertion (1) this
sequence could be written
$(a_\alpha,\phi_\alpha)=(f_A(\phi_\alpha),\phi_\alpha)$.  Because the
sequence $\{\phi_\alpha\}\subset \Ker D_{A,\vartheta}$ converged to zero in
$L^2$ and $D_{A,\vartheta}$ is elliptic, it would converge to zero in
$L^2_\ell$ (for $\ell\geq 2$). Since $f_A$ is continuous on $\sO_{K,A}$, we
would have
$\lim_{\alpha\to\8}a_\alpha=\lim_{\alpha\to\8}f_A(\phi_\alpha)=0$ and,
as $\bl\circ\bgamma_A$ is continuous on $\sO_{K,A}$,
$$
\lim_{\alpha\to\8}\bl(\bga_A(a_\alpha,\phi_\alpha))
=
\bl(\bga_A(0,0))
=
\bl(0)
=
0,
$$
contradicting our assumption that for all $\alpha$ we have
$\bl(\bga_A(a_\alpha,\phi_\alpha))=\eps$. This proves Assertion (3) and
completes the proof of the lemma.
\end{proof}

Because the spaces $\bga_A(\sO_A')$ suffice to cover the intersection
$\bar\sV(z)\cap\bL^{w,\eps}_{\ft,\ka}$ for $\eps$ sufficiently small by
Lemma \ref{lem:DeformingV}, we shall henceforth restrict the domain of
$\bga_A$ to $\sO_A'$.

We define a link of the submanifold $T_{[\hat
  A]}M^{w}_{\ka}\oplus\{0\}\subset T_{[\hat A]}M^w_{\ka}\oplus \Ker
D_{A,\vartheta}$ by setting
\begin{equation}
  \label{eq:bKAdelta}
\bK_{A,\delta}
=
\{(a,\phi)\in T_{[\hat A]}M^w_{\ka}\oplus \Ker D_{A,\vartheta}:
\|\phi\|_{L^2}^2 = \delta \}.
\end{equation}
The link $\bK_{A,\delta}/S^1$ is more convenient to work with than the
level sets of $\bl\circ\bga_A$ defining
$\bga_A^{-1}\left(\bL^{w,\eps}_{\ft,\ka}\right)$.  We will in Lemma
\ref{lem:Cobordism} that the two
links are related by an oriented cobordism. Therefore, prior to showing
this equivalence, we first discuss the orientation of the spaces.

An orientation $O$ for $\sM^{*,0}_{\ft}$ determines an orientation
\begin{equation}
\label{eq:ASDLinkOrientation1}
\om(\bL,\rd O)
\end{equation}
for $\bL^{w,\eps}_{\ft,\ka}$ by considering $\bL^{w,\eps}_{\ft,\ka}$ as the
boundary of the subspace $\ell^{-1}([\eps,\infty)) \subset \sM_{\ft}/S^1$ and
using the convention \eqref{eq:BoundaryOrientation}.

An orientation $o$ of $T_{[\hat A]}M^w_{\ka}$ determines an orientation of
$\bK_{A,\delta}/S^1$ by identifying $\bK_{A,\delta}/S^1$ with $T_{[\hat
A]}M^w_{\ka}\times\CC\PP^{k-1}$, where $k=\dim_\CC\Ker D_{A,\vartheta}$ and
taking the complex orientation of $\CC\PP^{k-1}$, denoted by
\begin{equation}
\label{eq:ComplexOrientationForK}
\om(\bK,o).
\end{equation}
We now describe a convention for orienting smooth submanifolds.

\begin{conv}
\label{conv:NormalBundleOrientation}
Suppose $Z$ and $M$ are manifolds which intersect transversely. Then an
orientation $O$ for $M$ and an orientation $\om(Z)$ for the normal bundle
of $Z$ determine an orientation for $M\cap Z$, which we write as
$O/\om(Z)$.
\end{conv}

{}From equation \eqref{eq:ASDKuranishiMap} and the fact that $\varphi_A$
vanishes transversely on $(\sO_A' - (T_{[\hat A]}M^w_{\ka}\oplus
\{0\})/S^1$, the fibers of 
the normal bundle of $\sZ_A/S^1$ in $(\sO_A' - (T_{[\hat A]}M^w_{\ka}\oplus
\{0\})/S^1$ are naturally identified with 
$\Coker D_{A,\vartheta}$. Let $\omega(\sZ)$ be the orientation of this
normal bundle of $\sZ_A/S^1$ obtained by giving $\Coker D_{A,\vartheta}$
the complex orientation. 

By Convention \ref{conv:NormalBundleOrientation}, the orientation
$\om(\bK,o)$ of $\bK_{A,\delta}/S^1$ and orientation $\om(\sZ)$ of the
normal bundle of $\sZ_A/S^1$ determine a ``complex orientation'' of
$\sZ_A\cap\bK_{A,\delta}/S^1$,
\begin{equation}
\label{eq:ComplexOrientOfZK}
\om(\bK,o)/\om(\sZ),
\end{equation}
where $o$ is the orientation for $T_{[\hat A]}M^w_{\ka}$, by analogy with
Definition \ref{defn:ComplexOrientation} which this construction matches.

However, to compare the orientations of $\bL^{w,\eps}_{\ft,\kappa}$ and
$\sZ_A\cap\bK_{A,\delta}/S^1$ with those of other links in a cobordism
formula such as equation \eqref{eq:CompactReductionFormula1}, it is natural
to orient $\sZ_A\cap\bK_{A,\delta}/S^1$ as a boundary of the cobordism.  If
$O$ is an orientation of $\sM_{\ft}^{*,0}/S^1$, we obtain an orientation 
\begin{equation}
\label{eq:BoundaryOrientationKZ}
\om(\sZ\cap\bK,\rd O)
\end{equation}
for $\sZ_A\cap\bK_{A,\delta}/S^1$ by identifying this manifold, via the map
$\bga_A$, with the boundary of
$$
\sM_{\ft}\less \bga_A\left( T_{[\hat A]}M^w_\ka\times
B_A(0,\delta)\cap\sO_A'\right),
$$
where $B_A(0,\delta)\subset \Ker D_{A,\vartheta}$ is the ball of radius
$\delta$, and using convention
\eqref{eq:BoundaryOrientation}. The proof of the following lemma is the
same as that of Lemma \ref{lem:ReducibleOrientation}.

\begin{lem}
\label{lem:ComplexASDOrientationComparison}
Let $w$ be an integral lift of $w_2(\ft)$.  Fix an orientation
$o=o(\Omega,w)$ of $M^w_{\ka}$.  Let $O=O^{\asd}(\Om,w)$ be the orientation
for $\sM^{*,0}/S^1$ in Definition \ref{defn:ASDOrient}. Then, the
orientations \eqref{eq:ComplexOrientOfZK} and
\eqref{eq:BoundaryOrientationKZ} for $\sZ_A\cap\bK_{A,\delta}/S^1$ agree,
that is
$$
\om(\bK,o)/\om(\sZ) =
\om(\sZ\cap\bK,\rd O).
$$
\end{lem}

By the definition of a geometric representative (Definition
\ref{defn:GeomRepresentative}), the normal bundle of $\sV(z)$ in
$\sB^{w,*}_{\ka}$ has an orientation, which we denote by $\om(\sV)$.
Because $\sV(z)$ intersects
$$
\pi_{\sB}\circ\bga_A(\bK_{A,\delta}/S^1),\quad
\pi_{\sB}\circ\bga_A(\sZ_A\cap\bK_{A,\delta}/S^1), \quad\text{and}\quad
\pi_{\sB}(\bL^{w,\eps}_{\ft,\ka})
$$
transversely, the orientations $\om(\sV)$, $\om(\bK,o)$, $\om(\sZ)$,
$\om(\sZ\cap\bK,\rd O)$ and $\om(\bL,\rd O)$ determine orientations
\begin{equation}
\label{eq:VZKOrientation}
\begin{aligned}
\om(\bK,o) / \om(\sV)
    &\quad\text{for}\quad
    \bga_A^{-1}(\sV(z))\cap\bK_{A,\delta}/S^1,
\\
\left.\begin{aligned}
\left(\om(\bK,o) / \om(\sZ)\right) / \om(\sV)&\\
\om(\sZ\cap\bK,\rd O)/ \om(\sV)
\end{aligned}\right\}
    &\quad\text{for}\quad
    \bga_A^{-1}(\sV(z))\cap \sZ_A\cap\bK_{A,\delta}/S^1,
\\
\om(\bL,\rd O) / \om(\sV)
    &\quad\text{for}\quad
    \sV(z)\cap\bL^{w,\eps}_{\ft,\eps}.
\end{aligned}
\end{equation}
Observe that Lemma \ref{lem:ComplexASDOrientationComparison} implies that
\begin{equation}
\label{eq:CplxBoundaryOrientationEqualityWithV}
\left(\om(\bK,o) / \om(\sZ)\right) / \om(\sV)=\om(\sZ\cap\bK,\rd O)/ \om(\sV),
\end{equation}
if $o=o(\Om,w)$ and $O=\Omega^{\asd}(\Om,w)$.

\begin{lem}
\label{lem:OrientationOfV}
Let $w$ be an integral lift of $w_2(\ft)$.  Let $\eps(A)=\pm 1$ be the
signed intersection number of $\sV(z)$ and $M^w_{\ka}$ at $[\hat A]$, where
$M^w_\ka$ is given the orientation $o(\Omega,w)$.
Then the following map is a diffeomorphism,
\begin{equation}
\label{eq:DefineCPDiffeom}
g_A=f_A\times\id_{\Ker D_{A,\vartheta}}:
\CC\PP^{k-1} \to \bga_{A}^{-1}\left(\sV(z)\right)\cap \bK_{A,\delta}/S^1,
\end{equation}
where $\CC\PP^{k-1}=\PP(\Ker D_{A,\vartheta})$ and
$f_A$ is defined in Assertion (1) of Lemma \ref{lem:DeformingV}. If
$\CC\PP^{k-1}$ has the complex orientation
and $\bga_{A}^{-1}\left(\sV(z)\right)\cap \bK_{A,\delta}/S^1$ has
the orientation $\om(\bK,o)/ \om(\sV)$ of \eqref{eq:VZKOrientation}
for $o=o(\Om,w)$,
then $g_A$ preserves orientation if and only if $\eps(A)=1$.
\end{lem}

\begin{proof}
By Assertion (1) of Lemma \ref{lem:DeformingV}, the intersection
$\bga_{A}^{-1}\left(\sV(z)\right)\cap \bK_{A,\delta}/S^1$ is given by
the $S^1$ quotient of the graph of $f_A$ restricted to $S^{2k-1}\subset
\Ker D_{A,\vartheta}$.  Thus, $g_A$ gives the desired diffeomorphism.
Because $\sV(z)$ and $M^w_\ka$ intersect transversely in $\sB^{w,*}_{\ka}$
at $[\hat A]$, the normal bundle of $\sV(z)$ in $\sB^{w,*}_{\ka}$ is
identified with $T_{[\hat A]}M^w_{\ka}$ but their orientations agree if and
only if $\eps(A)=1$.  The result then follows from the definition
\eqref{eq:ComplexOrientationForK} of the
orientation $\om(\bK,o)$ of $\bK_{A,\delta}/S^1$ determined by the
orientation $o(\Omega,w)$ of $T_{[\hat A]}M^w_\ka$.
\end{proof}

\begin{lem}
\label{lem:Cobordism}
Continue the assumptions and notation of Lemma
\ref{lem:DeformingV}.  For $\eps$ sufficiently small and $\delta$
as in Assertion (3) of Lemma \ref{lem:DeformingV} and generic,
there is a smooth, compact, and oriented cobordism between
\begin{equation}
\label{eq:Link1}
\bga_A^{-1}\left( \sV(z)\cap \bL^{w,\eps}_{\ft,\ka}\right),
\end{equation}
with the orientation $\om(\bL,\rd O)/\om(\sV)$ of
\eqref{eq:VZKOrientation}
for $O=O^{\asd}(\Omega,w)$, and the manifold
\begin{equation}
\label{eq:Link2}
\sZ_{A}\cap \bga_A^{-1}\left(\sV(z)\right)\cap\bK_{A,\delta}/S^1,
\end{equation}
with the orientation $\left(\om(\bK,o) / \om(\sZ)\right) / \om(\sV)$
of \eqref{eq:VZKOrientation}, where $o=o(\Om,w)$.
\end{lem}

\begin{proof}
As before, we let $B_{A}(0,\delta)\subset\Ker D_{A,\vartheta}$ be the open
ball of radius $\delta$.  Assertion (3) of Lemma \ref{lem:DeformingV}
yields the inclusion
$$
\bga_A
\left(
\left(T_{[\hat A]}M^w_{\ka}\times B_A(0,\delta)\right)\cap \sO_A'/S^1
\right)\cap\bar\sV(z)\cap\bar\sM_{\ft}/S^1\ \subset\ \bl^{-1}([0,\eps)).
$$
Then for generic values of $\eps$ and $\delta$,
\begin{equation}
\label{eq:Cobordism1}
\bar\sM_{\ft}/S^1\cap\bar\sV(z) \cap \bl^{-1}([0,\eps])\cap\bga_A(\sO_A')
-
\bga_A\left(T_{[\hat A]}M^w_{\ka}\times B_A(0,\delta)/S^1\right)
\end{equation}
is a smooth manifold with boundaries given by the manifold \eqref{eq:Link1}
and by
\begin{equation}
\label{eq:CobordismBoundary}
\bar\sM^0_{\ft}\cap\bar\sV(z)
\cap
\bga_A
\left(
   \left(
      T_{[\hat A]}M^w_{\ka}\times \rd \bar B_A(0,\delta)
   \right)
\cap
\sO_A'/S^1
\right),
\end{equation}
which is diffeomorphic, via $\bga_A$, to the manifold
\eqref{eq:Link2}.

By Assertion (2) of Lemma \ref{lem:DeformingV}, the compact set
$\bar\sV(z)\cap\bar\sM_{\ft}\cap\bl^{-1}([0,\eps])$ is contained in a
finite, disjoint union $\cup_{i=1}^N \bga_{A_i}\left(\sO_{A_i}'/S^1\right)$,
where $\bar\sV(z)\cap \barM_\kappa^w = \{[\hat A_i]\}_{i=1}^N$. Thus, each
component
$$
\bar\sV(z)\cap\bar\sM_{\ft}\cap\bl^{-1}([0,\eps])
\cap\bga_A\left(\sO_A'/S^1\right),
$$
of this disjoint union is compact. Therefore, the space \eqref{eq:Cobordism1}
is a compact and smooth cobordism between the manifold \eqref{eq:Link1} and
the manifold \eqref{eq:CobordismBoundary} which, as
previously noted, is diffeomorphic to
the manifold \eqref{eq:Link2}.

Let $\sM^{*,0}_{\ft}/S^1\cap\sV(z)$ have the orientation determined by the
orientation $O^{\asd}(\Omega,w)$ of $\sM_{\ft}/S^1$ and the orientation
$\om(\sV)$ of the normal bundle of $\sV(z)$.  The manifold
\eqref{eq:Cobordism1} has codimension zero in
$\sM^{*,0}_{\ft}/S^1\cap\sV(z)$ and thus inherits an orientation from
$\sM^{*,0}_{\ft}/S^1\cap\sV(z)$.  Hence, the manifold \eqref{eq:Cobordism1}
defines an oriented cobordism.  The orientation of the manifold
\eqref{eq:Link1} given by viewing it as a component of the boundary of
the oriented manifold \eqref{eq:Cobordism1} is then equal to
$-\om(\bL,\rd O)/\om(\sV)$ (as defined in \eqref{eq:VZKOrientation}) for
$O=O^{\asd}(\Om,w)$.  The
negative sign arises because the orientation $\om(\bL,\rd O)$ of
\eqref{eq:ASDLinkOrientation1} is defined by viewing
$\bL^{w,\eps}_{\ft,\ka}$ as the boundary of $\ell^{-1}([\eps,\infty))$.

The orientation of the manifold \eqref{eq:Link2} given by considering it as
a component of the boundary of the oriented manifold \eqref{eq:Cobordism1}
is then equal to $\om(\sZ\cap\bK,\rd O)/ \om(\sV)$ (as defined in
\eqref{eq:VZKOrientation}) for
$O=O^{\asd}(\Om,w)$ which, by equation
\eqref{eq:CplxBoundaryOrientationEqualityWithV} is equal to
$\left(\om(\bK,o) / \om(\sZ)\right) / \om(\sV)$ for $o=o(\Om,w)$.

Recall that two oriented manifolds $(M_i,\om_i)$
for $i=0,1$ are cobordant if there is an oriented manifold $W$ whose
oriented boundary is $(M_0,-\om_0)\cup (M_1,\om_1)$ \cite[p. 170]{Hirsch}.
Therefore, the manifold \eqref{eq:Cobordism1} gives the desired cobordism.
\end{proof}

We now give  a cohomological description of the zero locus
$\sZ_A$ of the obstruction map:

\begin{lem}
\label{lem:HomologyOfMonopolesNearASD}
Continue the hypotheses and notation of Lemmas \ref{lem:DeformingV} and
\ref{lem:OrientationOfV}.
Assume that $n_a(\ft)=\ind_\CC D_{A,\vartheta}$ is positive and let
$c=\dim_\CC\Coker D_{A,\vartheta}$.
If $\CC\PP^{k-1}\cong \PP(\Ker D_{A,\vartheta})$
has the complex orientation and $h\in
H^2(\CC\PP^{k-1};\ZZ)$ is the positive generator, then for generic
$\delta>0$, there is a smooth submanifold $T$ of $\CC\PP^{k-1}$ which
is Poincar\'e dual to $h^c$ such that the restriction of the map $g_A$ of
definition \eqref{eq:DefineCPDiffeom} to $T$ gives a diffeomorphism,
$$
g_A: T \simeq
\sZ_A\cap
\bga_A^{-1}\left(\sV(z)\right)
\cap \bK_{A,\delta}/S^1.
$$
If $T$ is oriented as the Poincar\'e dual of $h^c$ and
$\sZ_A\cap\bga_A^{-1}\left(\sV(z)\right)\cap \bK_{A,\delta}/S^1$
is oriented by $\om(\sZ\cap\bK,o)/\om(\sV)$ from
\eqref{eq:VZKOrientation} where $o=o(\Om,w)$, then
the restriction of $g_A$ to $T$ is orientation preserving
if and only if $\eps(A)=1$,
where $\eps(A)$ is defined in Lemma \ref{lem:OrientationOfV}.
\end{lem}

\begin{proof}
As noted before Lemma \ref{lem:DeformingV}, the Kuranishi map $\bvarphi_A$
in \eqref{eq:ASDKuranishiMap} vanishes transversely on
$\sO_A'-(T_{[\hat A]}M^w_{\ka}\oplus\{0\})$.
For generic $\delta$, the map $\bvarphi_A$
vanishes transversely on $\bga_A^{-1}\left(\sV(z)\right)\cap
\bK_{A,\delta}/S^1$ because $\sV(z)$ is transverse to
$\sM^{*,0}_{\ft}/S^1$ by construction.  This implies that the zero locus of
$\bvarphi_A$ is Poincar\'e dual to the Euler class of the vector bundle
\eqref{eq:ASDObstructionVectorBundle} of which $\bvarphi_A$ is a section.
We define a smooth submanifold,
$$
T=g_A^{-1}(\bvarphi_A^{-1}(0)) \subset \CC\PP^{k-1},
$$
and observe that the diffeomorphism $g_A$ in equation
\eqref{eq:DefineCPDiffeom} restricts to the desired diffeomorphism of $T$.

{}From the definition of $\bvarphi_A$ in equation \eqref{eq:ASDKuranishiMap}
and of $g_A$ in equation \eqref{eq:DefineCPDiffeom}, we see that the
composition $\bvarphi_A\circ g_A$ can be viewed as an $S^1$-equivariant map
$$
\bvarphi_A\circ g_A:S^{2k-1}\to \CC^c,
$$
where $\CC^c\cong\Coker D_{A,\vartheta}$,
and thus a section of the vector bundle
\begin{equation}
\label{eq:ASDObstructionVectorBundle}
S^{2k-1}\times_{S^1}\CC^c \to \CC\PP^{k-1},
\end{equation}
where the $S^1$ action is diagonal since $(\bvarphi_A\circ
g_A)(e^{i\theta}z)=e^{i\theta}(\bvarphi_A\circ g_A)(z)$, for $z\in
S^{2k-1}$ and $e^{i\theta} \in S^1$.  Because the action is diagonal, the
Euler class of this bundle is $h^c$ (see \cite[Lemma 3.27]{FL2a} for a
further explanation of the sign).

By Lemma \ref{lem:OrientationOfV}, the diffeomorphism
$g_A$ defines an orientation-preserving diffeomorphism from
$\CC\PP^{k-1}$ to $\bga_A^{-1}\left(\sV(z)\right)\cap\bK_{A,\delta}/S^1$
(with the orientation $\om(\bK,o)/\om(\sV)$ of \eqref{eq:VZKOrientation})
if and only if $\eps(A)=1$.
Recall that the orientation $\left(\om(\bK,o) / \om(\sZ)\right) / \om(\sV)$
from \eqref{eq:VZKOrientation} of
$\sZ_A\cap\bga_A^{-1}\left(\sV(z)\right)\cap \bK_{A,\delta}/S^1$
is given by the orientation $\om(\bK,o)/\om(\sV)$  of
$\bga_A^{-1}\left(\sV(z)\right)\cap\bK_{A,\delta}/S^1$ and
the orientation $\om(\sZ)$ of the normal bundle of $\sZ_A/S^1$.
Thus, if $T$ is oriented as the Poincar\'e dual of $h^c$ and
thus has the orientation determined by the complex orientation
of $\CC\PP^{k-1}$ and the complex orientation of the normal
bundle of $T$, then the restriction of $g_A$ to $T$ is
orientation preserving if and only if $\eps(A)=1$.
\end{proof}

The final tool needed to compute intersection numbers with
$\bL^w_{\ft,\ka}$ in equation \eqref{eq:LinkOfASD} is the following
description of $\bga_A^{-1}(\sW)$:

\begin{lem}
\label{lem:CohomOnASDNormalSlice}
Continue the notation and assumptions of Lemmas
\ref{lem:DeformingV}, \ref{lem:OrientationOfV},
and \ref{lem:HomologyOfMonopolesNearASD}.
Then $(\bga_A\circ g_A)^*\mu_{c}=2h$, where $\mu_{c}$ is the cohomology
class \eqref{eq:DefineC1}.
\end{lem}

\begin{proof}
  Recall that $\mu_{c}$ is the first Chern class of the line bundle
  $\LL_{\ft}$ in definition \eqref{eq:DefineC1LineBundle}.
  The embedding $\bga_A$ and the map $g_A$ are $S^1$-equivariant so,
  noting that $\LL_{\ft}$ is defined by the $S^1$ action in equation
  \eqref{eq:C1ClassS1Action}, we have
\begin{equation}
\label{eq:PulledBackDetBundle}
(\bga_A\circ g_A)^*\LL_{\ft}
\cong
S^{2k-1}\times_{(S^1,\times -2)}\CC
\to
 S^{2k-1}/S^1,
\end{equation}
The bundle \eqref{eq:PulledBackDetBundle} has first Chern class $2h$, the
sign being positive because the $S^1$ action is diagonal
(see \cite[Lemma 3.27]{FL2a}).
\end{proof}

Using Lemma \ref{lem:CohomOnASDNormalSlice} we can prove the assertion of
Lemma \ref{lem:CyclesExtension} that $\iota(\barM^w_{\ka}) \subset \bar\sW$:

\begin{proof}[Proof of Assertion (2d) in Lemma \ref{lem:CyclesExtension}]
  Lemma \ref{lem:CohomOnASDNormalSlice} shows that $\sW$ will have
  non-trivial intersection with the normal cone of any point in
  $\iota(M^w_{\ka})\subset \sM_{\ft}$, where by ``normal cone'' we mean
  $\bga_A(\sZ_A\cap \{0\}\oplus\Ker D_{A,\vartheta})/S^1$.  Therefore, the
  closure of $\sW$ will contain all points in $\iota(M^w_{\ka})$ and thus
  $\iota(\barM^w_{\ka}) \subset \bar\sW$.
\end{proof}

We can now compute the intersection with the link.

\begin{prop}
\label{prop:LinkOfASD}
Let $\ft$ be a \spinu structure on a four-manifold $X$, with $w$ an
integral lift of $w_2(\ft)$ and $w\pmod{2}$ is good. We further assume
that $d_a(\ft) = \dim M^w_{\ka}\ge 0$ and that $n_a(\ft)=\ind_\CC D_A >0$.
Let $\delta_{c}$ be a non-negative integer such that 
$$
\deg(z)+2\delta_{c}=d_a+2n_a-2=\dim(\sM_{\ft}^{*,0}/S^1)-1.
$$  
Suppose $z\in\AAA(X)$ has degree $\deg(z)\geq d_a$ and is
intersection-suitable.  If $\bL^{w,\eps}_{\ft,\kappa}\cap
\iota(M_\kappa^w)$ is oriented as the boundary of
$\sM^{*,\ge\eps}_{\ft}/S^1$, where $\sM_{\ft}^{*,0}/S^1$ is given the
orientation $O^{\asd}(\Omega,w)$, then there is a positive constant
$\varepsilon_0$ such that for generic $\varepsilon\in(0,\varepsilon_0)$,
\begin{equation}
  \label{eq:LinkOfASD}
\#\left(\bar\sV(z)\cap \bar\sW^{\delta_{c}}\cap \bL^{w,\eps}_{\ft,\kappa}
\right)
=
\begin{cases}
2^{n_a-1}\#(\bar\sV(z)\cap \barM^w_{\ka}),
& \text{if $\deg(z) = d_a$},\\
0,
&
\text{if $\deg(z) > d_a$}.
\end{cases}
\end{equation}
Moreover, these intersection numbers are independent of
the choice of generic $\varepsilon<\varepsilon_0$.
\end{prop}

\begin{proof}
  By Assertion (2) of Lemma \ref{lem:DeformingV}, the pairing
  \eqref{eq:LinkOfASD} is a sum of local terms,
  \begin{equation}
\label{eq:ASDIntersectionSum}    
\begin{aligned}
{}&\#\left(\bar\sV(z)\cap \bar\sW^{\delta_{c}}\cap \bL^{w,\eps}_{\ft,\kappa}
\right)
\\
&=
\sum_{[A]\in\bar\sV(z)\cap\barM^w_{\ka}}
\#\left(\bar\sV(z)\cap \bar\sW^{\delta_{c}}\cap \bL^{w,\eps}_{\ft,\kappa}
\cap\bga_A\left(\sO_A'\right)
\right),
\end{aligned}
  \end{equation} where $\sO_A'$ is the neighborhood defined in Lemma
\ref{lem:DeformingV}.  If $\deg(z)>\dim M^w_\ka$, the intersection
$\bar\sV(z)\cap\barM^w_{\ka}$ is empty and so the sum is trivial. Hence, we
can assume $\deg(z)=\dim M^w_\ka$, so $\delta_{c}=n_a(\ft)-1$.
Let $c=\dim_\CC\Coker D_{A,\vartheta}$, so $k=\dim\Ker D_{A,\vartheta}=n_a+c$.
If $\bar\sV(z)$ has multiplicity $q$ (in the sense of Definition
\ref{defn:GeomRepresentative}), then
we can evaluate the terms in the sum
in equation \eqref{eq:ASDIntersectionSum} as
\begin{align*}
{}&\#
\left(\bar\sV(z)\cap \bar\sW^{n_a-1}\cap \bL^{w,\eps}_{\ft,\kappa}
\cap\bga_A\left(\sO_A'\right)\right)
\\
&=\#\left(
\bga_{A}^{-1}\left(\sV(z)\cap\bar\sW^{n_a-1}\cap \bl^{-1}(\eps)\right)
\cap \sZ_A/S^1\right)
\quad\text{(Equations \eqref{eq:KurHomeoZeroLocusToPU2Mod}
  \& \eqref{eq:ASDLink})} 
\\
&=\#
\left(
\bga_A^{-1}\left(\bar\sW^{n_a-1}\right)
\cap\bga_A^{-1}\left(\sV(z)\right)
   \cap\bK_{A,\delta}\cap\sZ_A/S^1
   \right)
   \quad\text{(Lemma \ref{lem:Cobordism})}
\\
&=q\#
\left( (\bga_A\circ g_A)^{-1}\left(\bar\sW^{n_a-1}\right)
      \cap g_A^{-1}(\sZ_A/S^1)\cap
      \CC\PP^{n_a+c-1}\right)
      \quad\text{(Lemma \ref{lem:OrientationOfV})}
\\
&=q\eps(A)
\left\langle  (2h)^{n_a-1}\smile h^c ,[\CC\PP^{n_a+c-1}]\right\rangle
\quad\text{(Lemmas \ref{lem:HomologyOfMonopolesNearASD}
            and \ref{lem:CohomOnASDNormalSlice})}
            \\
&=q\eps(A)2^{n_a-1}.
 \end{align*}
Hence, equation \eqref{eq:ASDIntersectionSum} simplifies to give
$$
\#  \left(\bar\sV(z)\cap \bar\sW^{\delta_{c}}\cap \bL^{w,\eps}_{\ft,\kappa}
\right)
=
q2^{n_a-1}\sum_{[A]\in\bar\sV(z)\cap\barM^w_{\ka}}\eps(A)
=
2^{n_a-1} \#\left(\bar\sV(z)\cap\barM^w_{\ka}\right),
$$
completing the proof of the proposition.
\end{proof}

As an application of Lemma \ref{lem:HomologyOfMonopolesNearASD}, we explain
why the moduli space $\sM_{\ft}$ contains solutions to the $\PU(2)$
monopole equations \cite[Equation (2.32)]{FL2a} which are distinct from the
anti-self-dual or reducible solutions.  Lemma
\ref{lem:HomologyOfMonopolesNearASD} yields the following analogue of
Taubes' existence theorem for solutions to the anti-self-dual equation for
$\SO(3)$ connections:

\begin{prop}
\label{prop:PU(2)MonopoleExist}
Let $\ft$ be a \spinu structure on a four-manifold $X$, where we allow
$b_2^+(X)\geq 0$, and suppose $w_2(\ft)\equiv w\pmod{2}$, for $w\in
H^2(X;\ZZ)$.  Assume that $w\pmod{2}$ is good.  If $n_a(\ft)>0$, then for a
generic, $C^\8$ pair $(g,\rho)$, consisting of a Riemannian metric and
Clifford map, and generic, $C^\8$ parameters $(\tau,\vartheta)$, the moduli
space $\sM_{\ft}^{*,0}(g,\rho,\tau,\vartheta)$ of irreducible,
non-zero-section $\PU(2)$ monopoles is non-empty if the moduli space
$M_{\ka}^{w,*}$ of irreducible, anti-self-dual $\SO(3)$ connections on
$\fg_{\ft}$ is non-empty.
\end{prop}

\begin{proof}
This follows from the fact that the Euler class of the obstruction bundle
in the Kuranishi model of $[A,0]\in \sM_{\ft}$ is non-trivial by Lemma
\ref{lem:HomologyOfMonopolesNearASD}.
\end{proof}

\subsection{The cobordism formula}
\label{subsec:CompactReducFormula} If $z\in\AAA(X)$ is
intersection-suitable and $\deg(z)+2\delta_{c}=\dim
(\sM_{\ft}^{*,0}/S^1)-1=d_a+2n_a-2$, then Corollary
\ref{cor:IntersectionOfSmoothLowerStrata} tells us that the
intersection
\begin{equation}
\label{eq:IntersectionWithBoundary}
\bar\sV(z)\cap\bar\sW^{\delta_{c}}\cap\bar\sM^{*,\ge \eps}_{\ft}/S^1
\end{equation}
is a union of smooth one-dimensional manifolds in $M^{*,\ge \eps}_{\ft}/S^1$.
The boundaries of these one-manifolds will lie either on
$\bL^{w}_{\ft,\kappa}$ or in a neighborhood of some reducible
monopole, possibly in a lower level. Proposition \ref{prop:LinkOfASD}
describes the intersection of this family of one-manifolds with
the component $\bL^{w}_{\ft,\kappa}$ of the boundary of
$\bar\sM^{*,\geq\eps}_{\ft}/S^1$.  In particular, we see there are
finitely many points in this boundary.

If $w_2(\ft)$ is good then (noting that we always assume $b_2^+(X)>0$), for
any splitting $\ft=\fs\oplus \fs\otimes L$, the space $M_{\fs}$ contains no
zero-section monopoles, \cite[Corollary 3.3]{FL2a}.  In \cite[Definition
3.22]{FL2a} we constructed a homology class $[\bL_{\ft,\fs}]$ of the link
of the family of reducibles $M_{\fs}$ contained in the top level
$\sM_{\ft}^0/S^1$.
By Lemma \ref{lem:OrientComparison}, the orientations $O^{\asd}(\Omega,w)$ and
$O^{\red}(\Omega,\ft,\fs)$ differ by $\frac{1}{4}(w-c_1(L))^2$.  {}From the
definition of the geometric representatives and Lemma
\ref{lem:ReducibleOrientation}, we have:

\begin{lem}
\label{lem:CohomOnReducibleLink}
Let $\ft$ be a \spinu structure over a four-manifold $X$, with $w$ an
integral lift of $w_2(\ft)$.  Assume that $w\pmod{2}$ is good. If
$z\in\AAA(X)$ is intersection-suitable and
$\deg(z)+2\delta_{c}=d_a+2n_a-2$, then
$$
\#
\left(\bar\sV(z)\cap\bar\sW^{\delta_{c}}\cap\bL_{\ft,\fs}\right)
= (-1)^{\textstyle{\frac{1}{4}}(w-c_1(L))^2} \langle
\mu_{p}(z)\smile\mu_{c}^{\delta_{c}},[\bL_{\ft,\fs}]\rangle,
$$
where $\bL_{\ft,\fs}$ is given the boundary orientation
determined by $O^{\asd}(\Omega,w)$ on the left hand side
and the complex orientation of Definition
\ref{defn:ComplexOrientation} on the right hand side of the above
identity.
\end{lem}

We now characterize the \spinc structures, $\fs$, for which $\ft=\fs\oplus
\fs'$.

\begin{lem}
\label{lem:SplittingSpinu}
A \spinu structure $\ft$ on $X$ admits a splitting, $\ft=\fs\oplus \fs'$,
if and only if
\begin{equation}
\label{eq:SplittingSpinu}
\left(c_1(\ft)-c_1(\fs)\right)^2
=
p_1(\ft).
\end{equation}
\end{lem}

\begin{proof}
  Assume $\ft = \fs\oplus \fs'$.  We may write $\fs'=\fs\otimes L$ for some
  line bundle $L$, with $\ft=(\rho,V)$ and $\fs=(\rho,W)$.  Then
  $V=W\otimes E$, where $E=\underline{\CC}\oplus L$, and $\fg_{\ft}\cong
  i\ubarRR\oplus L$.  Thus, $c_1(L)=c_1(\ft)-c_1(\fs)$ obeys
  $c_1(L)^2=p_1(\fg_{\ft})=p_1(\ft)$, as desired.
  
  Conversely, suppose $c_1(\ft)-c_1(\fs)$ obeys condition
  \eqref{eq:SplittingSpinu}. Let $L$ be a complex line bundle with
  $c_1(L)=c_1(\ft)-c_1(\fs)$.  {}From Lemma 2.3 in \cite{FL2a} we know
  that $V\cong W\otimes E$ for a complex rank-two bundle $E$ determined up
  to isomorphism by $\fs$ and $\ft$. Then
$$
c_1(E)
=
\textstyle{\frac{1}{2}}c_1(V^+)-c_1(W^+)
=
c_1(\ft) - c_1(\fs)
=
c_1(L),
$$
while 
$$
c_2(E) 
= 
-\textstyle{\textstyle{\frac{1}{4}}}(p_1(\su(E))-c_1(E)^2)
=
-\textstyle{\textstyle{\frac{1}{4}}}(p_1(\ft)-c_1(L)^2)
=
0,
$$
where the final equality follows from the fact that $p_1(\ft) =
c_1(L)^2$ by hypothesis.  Hence, $E\cong \underline{\CC}\oplus L$ and so
$V\cong W\oplus W\otimes L$, as desired.
\end{proof}

If $\fs$ is a \spinc structure with $c_1(\fs)$ obeying condition
\eqref{eq:SplittingSpinu}, there is a topological embedding $M_\fs\embed
\sM_{\ft}$ of $M_\fs$ into the top level of $\bar\sM_{\ft}$  
\cite[Lemma 3.13]{FL2a}. More generally, if $c_1(\fs)$ obeys
\begin{equation}
\label{eq:LowerLevelReducible}
(c_1(\ft)-c_1(\fs))^2 = p_1(\ft) + 4\ell,
\end{equation}
for some non-negative integer $\ell$, then there is a topological embedding
of $M_\fs$ into the lower-level $\PU(2)$-monopole moduli space
$\sM_{\ft_\ell}\times\Sym^\ell(X)$, where $\ft_\ell$ is a \spinu structure
with $p_1(\ft_\ell)=p_1(\ft)+4\ell$ \cite[Equation (2.44)]{FL2a}.

If the reducibles in $\bar\sM_{\ft}$ appear only in the top Uhlenbeck level
$\sM_{\ft}$, then $\bar\sM^{*,\geq\varepsilon}_{\ft}/S^1$ is a cobordism
between the link $\bL^{w,\eps}_{\ft,\kappa}$ of the stratum defined by the
anti-self-dual moduli space, $\iota(M^w_{\ka})$, and the links
$\bL_{\ft,\fs}$ of the strata of reducibles, $\iota(M_{\fs})$.  Counting
the points in the boundary of the oriented, one-dimensional manifold
\eqref{eq:IntersectionWithBoundary} then gives the identity:
\begin{equation}
\label{eq:RawCobordismFormula}
\#\left(\bar \sV(z)\cap\bar\sW^{\delta_{c}}\cap \bL_{\ft,\kappa}^w\right)
=
-\sum_{\{\fs:\ \fs\oplus\fs'=\ft\}}
\#\left(\bar \sV(z)\cap\bar\sW^{\delta_{c}}\cap \bL_{\ft,\fs}\right).
\end{equation}
Let $w$ be an integral lift of $w_2(\ft)$ defining the orientation
$O^{\asd}(\Omega,w)$ of $\sM_{\ft}^{*,0}/S^1$.  Lemmas
\ref{lem:OrientComparison} and
\ref{lem:ReducibleOrientation} imply that the orientation $O^{\asd}(\Omega,w)$
and the complex orientation for the link $\bL_{\ft,\fs}$ differ by
\begin{equation}
\label{eq:OrientationFactor}
\begin{aligned}
(-1)^{o_{\ft}(w,\fs)},\quad\text{where we define }
o_{\ft}(w,\fs)
&=
\textstyle{\textstyle{\frac{1}{4}}}(w -c_1(L))^2
\\
&=
\textstyle{\textstyle{\frac{1}{4}}}(w- c_1(\ft) +c_1(\fs))^2.
\end{aligned}
\end{equation}
Equation \eqref{eq:RawCobordismFormula}, Proposition \ref{prop:LinkOfASD}
and Lemma \ref{lem:CohomOnReducibleLink} then yield the following result.

\begin{thm}
\label{thm:CompactReductionFormula}
Let $\ft$ be a \spinu structure on an oriented, smooth four-manifold $X$
with $b_2^+(X)>0$ and $w_2(\ft)\equiv w\pmod{2}$, for $w\in H^2(X;\ZZ)$.
Assume that $w\pmod{2}$ is good. Suppose $z\in\AAA(X)$ has degree
\begin{equation}
\label{eq:degzLowerUpperBounds}
d_a(\ft) \leq \deg(z) \leq d_a(\ft) + 2n_a(\ft) - 2,
\end{equation}
and is intersection-suitable.  Assume that the set of isomorphism classes
of \spinc structures, $\fs\in\Spinc(X)$, defining reducible $\PU(2)$
monopoles in $\bar\sM_{\ft}$ all obey condition \eqref{eq:SplittingSpinu},
and so non-empty Seiberg-Witten moduli strata $\iota(M_\fs)$ appear
only in the top level, $\sM_{\ft}$.  

\medskip
\noindent{\em{\bf (a)}} If
for all $\fs\in \Spinc(X)$ with $M_{\fs}$ non-empty we have
\begin{equation}
\label{eq:NoReducibles}
(c_1(\ft) - c_1(\fs))^2 < p_1(\ft),
\end{equation}
so $\bar\sM_{\ft}$ contains no reducible monopoles, then
\begin{equation}
\label{eq:DVanishing}
\#\left(\bar\sV(z)\cap \bar M_{\kappa}^w(X)\right)
=
0.
\end{equation}
\noindent{\em{\bf (b)}} 
If $\deg(z) = d_a$ and $o_{\ft}(w,\fs)$ is as defined
in equation \eqref{eq:OrientationFactor}, then
\begin{equation}
\#\left(\bar \sV(z)\cap \bar M_{\kappa}^w(X)\right)
=
- 2^{1-n_a}\sum_{\{\fs:\ \fs\oplus\fs'=\ft\}}
(-1)^{o_{\ft}(w,\fs)}
\langle\mu_{p}(z)\smile\mu_{c}^{n_a-1},[\bL_{\ft,\fs}]\rangle,
\label{eq:CompactReductionFormula1}
\end{equation}
where the class $[\bL_{\ft,\fs}]$ is defined by the complex orientation of
Definition \ref{defn:ComplexOrientation} on $\bL_{\ft,\fs}$.

\noindent{\em{\bf (c)}} If $d_a<\deg(z)\leq d_a+2n_a-2$ and
$\delta_{c}\in\NN$ is defined by $\deg(z)+2\delta_{c} = d_a + 2n_a - 2$, then
\begin{equation}
\label{eq:CompactReductionFormula1IsZero}
\sum_{\{\fs:\ \fs\oplus\fs'=\ft\}}
(-1)^{o_{\ft}(w,\fs)} 
\langle\mu_{p}(z)\smile\mu_{c}^{\delta_{c}},[\bL_{\ft,\fs}]\rangle
= 0.
\end{equation}
\end{thm}

The sums in equations \eqref{eq:CompactReductionFormula1} and
\eqref{eq:CompactReductionFormula1IsZero} are necessarily finite, because
there are only finitely many \spinc structures $\fs$ with $M_\fs$ non-empty
and thus $\bL_{\ft,\fs}$ non-empty.  Theorem
\ref{thm:CompactReductionFormula} can be strengthened to a more useful form
if we assume

\begin{conj}
\label{conj:Multiplicity}
\cite[Conjecture 3.1]{FKLM}
Continue the notation of Theorem \ref{thm:CompactReductionFormula}. Suppose
that $\ft_\ell=\fs\oplus\fs'$, where $p_1(\ft_\ell) = p_1(\ft)+4\ell$ and
$\iota(M_{\fs})$ is contained in the level
$\sM_{\ft_\ell}\times\Sym^\ell(X)$, for some natural number $\ell\geq 0$.
Then the pairing $\#(\bar\sV(z)\cap\bar\sW^{\delta_{c}} \cap
\bL_{\ft,\fs})$ is a multiple of $\SW_X(\fs)$ and thus vanishes if
the Seiberg-Witten function $\SW_{X,\fs}$ is trivial.
\end{conj}

See \S \ref{subsec:SWinvariants} for a definition of the Seiberg-Witten
invariants.  The motivation for this conjecture is discussed in
\cite{FLGeorgia} and almost certainly does hold, one of our current goals
being to provide a complete proof in the near future. The difficulty lies
in the appropriate
construction of the link $\bL_{\ft,\fs}$ of a Seiberg-Witten moduli
space when $\ell=\ell(\ft,\fs)>0$.  We show that Conjecture
\ref{conj:Multiplicity} holds when $\ell=0$ in Theorem
\ref{thm:DegreeZeroFormula} and when $\ell=1$ in \cite{FLLevelOne}. By
adapting Leness's proof of the wall-crossing formula in \cite{LenessWC}, we
can also see that the conjecture holds when $\ell=2$.

\begin{cor}
\label{cor:CompactReductionFormula}
Given Conjecture \ref{conj:Multiplicity}, we can relax the hypothesis of
Theorem \ref{thm:CompactReductionFormula} 
---that {\em non-empty\/} Seiberg-Witten moduli spaces $M_{\fs}$ appear
only in the top level $\sM_{\ft}$---to the weaker requirement that
Seiberg-Witten moduli spaces $M_{\fs}$ with {\em non-trivial Seiberg-Witten
functions $SW_{X,\fs}$\/} appear only in the top level $\sM_{\ft}$.
Then the conclusions of Theorem \ref{thm:CompactReductionFormula} hold
without change.
\end{cor}

\begin{rmk}
\label{rmk:WhereDoTheSWSpacesLie}
The Seiberg-Witten stratum $\iota(M_{\fs})$
corresponding to a splitting $\ft_{\ell}=\fs\oplus\fs'$ lies in level
$$
\ell(\ft,\fs) = \textstyle{\frac{1}{8}}(d_a-2r(\Lambda,\fs))
$$ 
of the space of ideal $\PU(2)$ monopoles containing $\bar\sM_{\ft}$, from
the definition
\eqref{eq:SWInTopLevelFunction} of $r(\Lambda,\fs)$, where $d_a$ is the
dimension of the anti-self-dual moduli space $M_\kappa^w\embed\sM_{\ft}$.
Then, by definition
\eqref{eq:SWInTopLevelFunction} of $r(\Lambda)$, the Seiberg-Witten strata
with non-trivial invariants are contained only in levels $\ell$ of this
space of ideal $\PU(2)$ monopoles, where
$$
0\leq \ell\leq \textstyle{\frac{1}{8}}(d_a-2r(\Lambda)). 
$$
See the proof of Theorem \ref{thm:Main} in \S
\ref{subsec:ProofsOfMainTheorems}.
\end{rmk}


\section{Intersection with the link of a stratum of top-level reducibles}
\label{sec:DegreeZero} In this section we calculate the pairings
appearing on the right-hand-side of equation
\eqref{eq:CompactReductionFormula1} in Theorem
\ref{thm:CompactReductionFormula} under the additional assumption
that for all $\alpha,\alpha'\in H^1(X;\ZZ)$ one has $\alpha\smile\alpha'=0$.
We begin by giving a definition of the Seiberg-Witten invariants which is
appropriate for the perturbations we use in our version of the
Seiberg-Witten equations \cite[\S 2.3]{FL2a}.  Recall that the link
$\bL_{\ft,\fs}$ constructed in
\cite[\S 3.5]{FL2a} is diffeomorphic via a map $\bga$
to the zero-locus of a section of an obstruction bundle
$\bgamma^*\Xi/S^1$ over 
the complex projectivization of a bundle $N_{\ft}(\Xi,\fs)\to M_{\fs}$.
Thus, in \S \ref{subsec:ClassesOnReducibleLink}, we calculate pullbacks
of the classes $\mu_{p}$ and $\mu_{c}$ by $\bga$
to $\PP N_{\ft}(\Xi,\fs)$.  In \S
\ref{subsec:CharacteristicsOfRedNormal} we compute the total
Segre class---the formal inverse of the total Chern class as defined in
Lemma \ref{lem:Segre}---of the virtual normal bundle $N_{\ft}(\Xi,\fs)$ of
$M_{\fs}$ and the Euler class of the obstruction bundle
$\bgamma^*\Xi/S^1$. Finally,  
in \S \ref{subsec:Computation} we perform the actual computation and
complete the proofs of Theorems \ref{thm:DSWSeriesRelation}, \ref{thm:Main},
\ref{thm:FLthm}, and Corollary \ref{cor:FLthm}. 

\subsection{A definition of the Seiberg-Witten invariants}
\label{subsec:SWinvariants}
In this subsection we give a definition of the Seiberg-Witten invariants
for a closed, smooth four-manifold $X$; we allow $b_1(X)\geq 0$
and $b_2^+(X)\geq 1$.

Recall that $X$ is equipped with an orientation for which $b_2^+(X)>0$ and
that we have fixed an orientation for $H^1(X;\RR)\oplus H^+(X;\RR)$: the
Seiberg-Witten moduli spaces are then oriented according to the conventions
of \cite[\S 6.6]{MorganSWNotes} (or see our \S \ref{subsec:OrientSW}).
In \cite[\S 2.4.2]{FL2a} we defined a universal complex line bundle,
\begin{equation}
\label{eq:UniversalSWDefinition}
\LL_{\fs}
=
\tsC^0_{\fs}\times_{\sG_{\fs}}\underline{\CC}
\to
\sC^0_{\fs}\times X,
\end{equation}
where $\underline{\CC}=X\times\CC$ and the action of $\sG_{\fs}$ is given
for $s\in\sG_{\fs}$, $x\in X$ and $z\in\CC$ by
\begin{equation}
\label{eq:SWUniversalGAction}
\left( s, (B,\Psi),(x,z)\right)
\mapsto
\left( s(B,\Psi),(x,s(x)^{-1}z)\right).
\end{equation}
We then defined cohomology classes on $\sC^0_{\fs}$ by
\begin{equation}
\mu_{\fs}: H_\bullet(X;\RR)\to H^{2-\bullet}(\sC^0_{\fs};\RR),
\quad
\mu_{\fs}(\alpha)
=
c_1(\LL_{\fs})/\alpha,
\label{eq:SWMuMap}
\end{equation}
where $\alpha$ is either the positive generator $x\in H_0(X;\ZZ)$ or a
class $\gamma\in H_1(X;\RR)$.

If $z\in\BB(X)$ is a monomial $\alpha_1\cdots\alpha_p$ with
$\alpha_i\in H_0(X;\ZZ)$ or $H_1(X;\RR)$, then it has total
degree $\deg(z) = \sum_{i=1}^p\deg(\alpha_i)$. If
$z=x^m\gamma_1\cdots\gamma_n \in \BB(X)$, we set
\begin{equation}
  \label{eq:SWMuMonomial}
\mu_{\fs}(z)
=
\underbrace{\mu_{\fs}(x)\smile\cdots\smile\mu_{\fs}(x)}_{\text{$m$
times}}\smile\mu_{\fs}(\gamma_1)\smile\cdots\smile\mu_{\fs}(\gamma_n),
\end{equation}
and define $\mu_{\fs}(z)$ for arbitrary elements $z\in\BB(X)$ by
$\RR$-linearity.

Let $\tilde X=X\#\overline{\CC\PP}^2$ be the blow-up of $X$ with
exceptional class $e\in H_2(\tilde X;\ZZ)$ and denote its Poincar\'e dual
by $\PD[e]\in H^2(\tilde X;\ZZ)$. Let $\fs^\pm=(\tilde\rho,\tilde W)$ denote
the \spinc structure on $\tilde X$ with $c_1(\fs^\pm)=c_1(\fs)\pm \PD[e]$
obtained by splicing the \spinc structure $\fs=(\rho,W)$ on $X$ with the
\spinc structure on $\overline{\CC\PP}^2$ with first Chern class $\pm \PD[e]$.
(See \S \ref{subsec:LinkBlowUp} for a more detailed explanation of the
relation between \spinc and \spinu structures on $X$ and $\tilde X$, as
well as the blow-up formula for Seiberg-Witten invariants, which we shall
invoke below.) One easily checks that $\dim M_{\fs^\pm}(\tilde X) = \dim
M_{\fs}(X)$, where the Seiberg-Witten moduli spaces $M_{\fs}(X)$ and
$M_{\fs^\pm}(\tilde X)$ are defined in \cite[Equation (2.55) \& Lemma
3.12]{FL2a} with perturbation $\eta=F^+(A_{\Lambda})$. Here, $2A_\Lambda$
is the fixed connection on $\det(V^+)\cong\det(\tilde V^+)$ and
$\tilde\ft=(\tilde\rho,\tilde V)$ is the \spinu structure on $\tilde X$
defined in Lemma \ref{lem:BlowUpSpinu}, with $c_1(\tilde\ft)=c_1(\ft)$,
$p_1(\tilde\ft)=p_1(\ft)-1$, and $w_2(\tilde\ft)=w_2(\ft)+\PD[e]$. Now
$$
c_1(\fs^\pm)-c_1(\tilde\ft)=c_1(\fs)\pm \PD[e]-\Lambda \in H^2(\tilde X;\ZZ)
$$
is not a torsion class and so---for $b_2^+(X)>0$, generic Riemannian
metrics $g$ on $X$ and related metrics on the connect sum $\tilde X$---the
moduli spaces $M_{\fs^\pm}(\tilde X)$ contain no zero-section pairs
\cite[Proposition 6.3.1]{MorganSWNotes}. Thus, for our choice of generic
perturbations, the moduli spaces $M_{\fs^\pm}(\tilde X)$ are compact,
oriented, smooth manifolds.

Noting that $\BB(X)\cong\BB(\tilde X)$, we define the {\em Seiberg-Witten
  invariants for $(X,\fs)$\/} as an $\RR$-linear function
\eqref{eq:IntroSWFunction} by setting
\begin{equation}
  \label{eq:SWFunctionDefn}
  SW_{X,\fs}(z) = \langle\mu_{\fs^+}(z),[M_{\fs^+}(\tilde X)]\rangle,
\end{equation}
with $SW_{X,\fs}(z) = 0$ when $\deg(z) \neq \dim M_{\fs}(X)$.
The blow-up formula for Seiberg-Witten invariants (Theorem
\ref{thm:SWBlowUpFormula}) implies that
\begin{equation}
  \label{eq:SWDefnBlowUp}
\langle\mu_{\fs}(z),[M_{\fs}(X)]\rangle =
\langle\mu_{\fs^\pm}(z),[M_{\fs^\pm}(\tilde X)]\rangle,
\end{equation}
when the pairing on the left is well-defined, that is, when $M_{\fs}(X)$
contains no zero-section monopoles. For example, with our version of the
Seiberg-Witten equations \cite[Equation (2.55)]{FL2a}, this situation
arises when $c_1(\fs)-\Lambda\in H^2(X;\ZZ)$ is not a torsion class and
thus $M_\fs(X)$ contains no zero-section pairs if the metric $g$ is
generic and $b_2^+(X)>0$ \cite[Proposition 6.3.1]{MorganSWNotes}.
Therefore, our definition of the Seiberg-Witten
invariants coincides with the usual one \cite{MorganSWNotes} in this case,
but has the advantage that it is valid even when $c_1(\fs)-\Lambda$ is
torsion and one cannot perturb the Seiberg-Witten equations by a generic
two-form $\eta$ (see Remark 2.14 in \cite{FL2a}).
When $b^+_2(X)>1$, the pairing on the right-hand-side of
equation \eqref{eq:SWFunctionDefn} is independent of the metric,
\cite[Lemma 6.7.1]{MorganSWNotes}.

When $b^+_2(X)=1$, however, the pairing on the right-hand-side of definition
\eqref{eq:SWFunctionDefn} depends on the period point
$\omega(\tilde g)$ of the metric $\tilde g$ on $\tilde X$, as in the case
of the Donaldson invariants (see \S \ref{subsubsec:DefnDInvariants}). To
explain this dependence when the Seiberg-Witten moduli spaces are defined
as in \cite[\S 2.3]{FL2a} with the perturbation parameters
$\eta=F^+(A_{\Lambda})$ described above, we note that
the moduli space $M_{\fs^+}(\tilde g)$ contains zero-section pairs
if and only if the period point $\omega(\tilde g)$ lies on the
$(c_1(\fs^+)-\Lambda)$-wall in the positive cone of $H^2(\tilde X;\RR)$,
that is
\begin{equation}
\label{eq:DefiningChamber}
\omega(\tilde g)\smile (c_1(\fs^+)-\Lambda)=0.
\end{equation}
When $\omega(\tilde g)$ does not lie on the wall, the pairing in
definition \eqref{eq:SWFunctionDefn} may depend on the sign of
$\omega(\tilde g)\smile ( c_1(\fs^+)-\Lambda)$.
The chambers for the Seiberg-Witten invariants of $M_{\fs^+}$ are thus
connected components of the complement of the
$(c_1(\fs^+)-\Lambda)$-wall in the positive cone of $H^2(\tilde X;\RR)$,
which we call $(c_1(\fs^+)-\Lambda)$-chambers.

By an argument which is the same as the one we gave for the Donaldson
invariants in \S \ref{subsubsec:DefnDInvariants}, if $c_1(\fs)-\Lambda$ is
not torsion then each $(c_1(\fs)-\La)$-chamber in the positive cone of
$H^2(X;\RR)$ is contained in a unique $(c_1(\fs^+)-\La)$-chamber in
the positive cone of
$H^2(\tilde X;\RR)\cong \RR[e]\oplus H^2(X;\RR)$, the {\em related
chamber\/}.  The Seiberg-Witten invariant associated to a
$(c_1(\fs)-\La)$-chamber is then defined by evaluating the pairing in
equation \eqref{eq:SWFunctionDefn} with a metric whose period point lies in
the related $(c_1(\fs^+)-\La)$-chamber.

Suppose $w_2(X)-\Lambda\pmod{2}$ is good, in the sense of Definition
\ref{defn:Good}. For any \spinc structure $\fs$ over $X$, we have
$c_1(\fs)\equiv w_2(X)\pmod 2$ and so $c_1(\fs)-\Lambda\pmod{2}$ is
good. Then $c_1(\fs)-\Lambda$ is not torsion and the Seiberg-Witten
invariants for $M_{\fs}$ depend only on the metric $g$ through the
$(c_1(\fs)-\Lambda)$-chamber for $\omega(g)$.

If $w_2(X)-\Lambda\pmod{2}$ is not good, then $c_1(\fs)-\Lambda$ may be
torsion and in this situation
$$
\omega(\tilde g)\smile (c_1(\fs^+)-\Lambda)
=
\omega(\tilde g)\smile \PD[e],
$$
so the sign of the cup-product would depend on the sign of $\omega(\tilde
g)\smile \PD[e]$, which converges to zero as the neck is stretched. Hence,
the definition
\eqref{eq:SWFunctionDefn} of the Seiberg-Witten invariant in this case
requires a more delicate analysis of the sign of $\omega(\tilde g)\smile
\PD[e]$ as the length of the neck converges to infinity \cite{Yang}, which
we shall not consider here.

Thus, when $b_2^+(X)=1$, we shall assume that $w_2(X)-\Lambda\pmod{2}$ is
good. Since
$$
w \equiv w_2(X)-\Lambda \pmod{2},
$$
this coincides with the constraint we used to define the Donaldson
invariants in \S \ref{subsubsec:DefnDInvariants}.  

We now compare the Donaldson and Seiberg-Witten chamber structures:

\begin{lem}
\label{lem:IdentifyChambers}
Let $\ft$ be a \spinu structure on a four-manifold $X$ with $b_2^+(X)=1$,
where $w$ is an integral lift of $w_2(\ft)$ and $w\pmod{2}$ is good, and
$c_1(\ft)=\Lambda$.  Then there is a one-to-one correspondence between
the set of $(w,p_1(\ft))$-walls and the set of $(c_1(\fs)-\La)$-walls, where
$M_{\fs}$ is contained in the space of ideal $\PU(2)$ monopoles,
$\cup_{\ell=0}^\infty\left( \sM_{\ft_\ell}\times\Sym^\ell(X)\right)$.
\end{lem}

\begin{proof}
A $(w,p_1(\ft))$-wall is defined by class $\alpha\in H^2(X;\ZZ)$ with
$\alpha\equiv w\pmod 2$ and $\alpha^2=p_1(\ft)+4\ell$ for $\ell\ge 0$ (see
equation \eqref{eq:WallCondition}).  Because $\alpha$ is an integral lift
of $w_2(\ft)$, the class $\La-\alpha$ is characteristic.  Hence, there is a
\spinc structure $\fs$ with $c_1(\fs)=\La-\alpha$.  By Lemma
\ref{lem:SplittingSpinu} and the identity $\alpha^2=p_1(\ft)+4\ell$, a
\spinu structure $\ft_\ell$ with $c_1(\ft_\ell)=\La$ and
$p_1(\ft_\ell)=p_1(\ft)+4\ell$ admits a splitting $\ft_\ell=\fs\oplus \fs'$
for any such \spinc structures.  Conversely, given a \spinc structure $\fs$
with $M_{\fs}$ contained in the space of ideal monopoles, Lemma
\ref{lem:SplittingSpinu} implies that $(c_1(\fs)-\La)^2=p_1(\ft)+4\ell$ and
$c_1(\fs)-\La\equiv w_2(\ft)\pmod 2$, so the $(c_1(\fs)-\La)$-wall is a
$(w,p_1(\ft))$-wall.
\end{proof}

\begin{rmk}
\label{rmk:SameChambers}
Lemma \ref{lem:IdentifyChambers} implies that in formulas such as
\eqref{eq:RawGeneralCobordismFormula} and \eqref{eq:GeneralFormulaForPairing}
derived using the cobordism $\sM_{\ft}^{*,0}/S^1$ to compare Donaldson and
Seiberg-Witten invariants, if a period point $\om(g)$ crosses
a wall for the Donaldson invariant, it will also cross a wall
for one of the Seiberg-Witten invariants in the formula; thus
both sides of the identity will change.
\end{rmk}

\subsection{Pullbacks of cohomology classes to the link of a
stratum of reducibles}
\label{subsec:ClassesOnReducibleLink}
We now compute the pullbacks  of the cohomology classes
$\mu_{p}(z)$ and $\mu_{c}$ by
$\bga:N^\eps_{\ft}(\Xi,\fs)/S^1\to \sC_{\ft}/S^1$, the
restriction of the $S^1$-equivariant embedding
defined in \cite[Equation (3.44)]{FL2a} to the $\eps$-sphere of the bundle
$N_{\ft}(\Xi,\fs)$.
The result will be expressed in terms of the Seiberg-Witten
$\mu_{\fs}$-classes and one additional cohomology class:

\begin{defn}
\label{defn:nu}
Let $\nu\in H^2(\PP N_{\ft}(\Xi,\fs);\ZZ)$ be the
negative of the first Chern class of the $S^1$ bundle
$N^\eps_{\ft}(\Xi,\fs)\to \PP N_{\ft}(\Xi,\fs)$.
Restricted to each fiber of $\PP N_{\ft}(\Xi,\fs)$, the class
$\nu$ is the positive generator of the cohomology. With the conventions of
\cite[\S 3.1]{Fulton}, the class $\nu$ is the first Chern class of
the line bundle $\sO_{\PP N_{\ft}(\Xi,\fs)}(1)$, the
dual of the tautological bundle.
\end{defn}

Let $\tN_{\ft}(\Xi,\fs)\to\tM_{\fs}$ be the
pullback of $N_{\ft}(\Xi,\fs)$ by the projection $\tM_{\fs}\to M_{\fs}$.
To compute the pullbacks by $\bga$ of the cohomology classes $\mu_{p}(\beta)$
to $\PP N_{\ft}(\Xi,\fs)$, we first compute the pullback of the universal
bundle $\FF_{\ft}$ defined in equation \eqref{eq:UniversalSO(3)Bundle}.

\begin{lem}
\label{lem:UniversalReduction}
Let $\ft$ be a \spinu structure which admits a splitting $\ft=\fs\oplus
\fs\otimes L$.  Assume $M_{\fs}$ contains no zero-section pairs.  Let
$\bgamma:N_{\ft}^\eps(\Xi,\fs)\embed \sC_{\ft}$ be the embedding constructed in
\cite[\S 3.5.3]{FL2a}, and let $\tN^\eps_{\ft}(\Xi,\fs)$ denote the
$\eps$-sphere bundle of $\tN_{\ft}(\Xi,\fs)$. Then, we have an isomorphism
of $\SO(3)$ bundles over $\PP N_{\ft}(\Xi,\fs)\times X$,
\begin{equation}
\label{eq:UniversalBundleReduction}
(\bga\times\id_X)^*\FF_{\ft}
\cong
\tN^\eps_{\ft}(\Xi,\fs)\times_{\sG_{\fs}\times S^1} (i\ubarRR\oplus L),
\end{equation}
where $s\in\sG_{\fs}$ and $e^{i\theta}\in S^1$ act on
$\tN^\eps_{\ft}(\Xi,\fs)\times (i\ubarRR\oplus L)$ by
\begin{equation}
\label{eq:ReducibleUniversalAction}
((B,\Psi,\eta),f\oplus z)
\mapsto
\left(s(B,\Psi, e^{i\theta}\eta), f \oplus s^{-2}e^{i\theta}z\right),
\end{equation}
where $(B,\Psi,\eta)\in\tN_{\ft}(\Xi,\fs)$, $(B,\Psi)\in \tM_{\fs}$,
and $f\oplus z\in i\ubarRR\oplus L$.
\end{lem}

\begin{proof}
Since $V=W\oplus W\otimes L$, where $\fs=(\rho,W)$ and $\ft=(\rho,V)$, we
have an isomorphism of $\SO(3)$ bundles $\fg_{\ft}\cong
i\underline{\RR}\oplus L$ and the definition
\eqref{eq:UniversalSO(3)Bundle} of $\FF_{\ft}$ yields
isomorphism of $\SO(3)$ bundles over $\sC_{\ft}^*\times X$,
$$
\FF_{\ft}
\cong
\tsC_{\ft}^*\times_{\sG_{\ft}\times S^1} (i\underline{\RR}\oplus L).
$$
{}From \cite[\S 3.5.4]{FL2a} we recall that the embedding
$\bgamma:N_{\ft}(\Xi,\fs)\to \sC_{\ft}$ lifts to a map
$$
\tilde\bgamma:\tN_{\ft}(\Xi,\fs)\to \tsC_{\ft},
$$
which is $\sG_{\fs}$ equivariant when $s\in\sG_{\fs}$ acts
on $\tsC_{\ft}$ via the embedding
\begin{equation}
\label{eq:DefnGInclusion}
\varrho:\sG_{\fs}\embed\sG_{\ft},
\quad
s\mapsto \varrho(s) = s\,\id_{W}\oplus s^{-1}\id_{W\otimes L},
\end{equation}
while $s$ acts on the base $\tilde M_{\fs}$ by the usual action of
$\sG_{\fs}$ and on the fibers of $\tN_{\ft}(\Xi,\fs)\to \tilde M_{\fs}$
by the action on $L^2_k(\Lambda^1\otimes L)\oplus L^2_k(W^+\otimes L)$
induced by the isomorphisms $\fg_{\ft}\cong i\underline{\RR}\oplus L$ and
$V=W\oplus W\otimes L$ and the action of $\varrho(s)$ on $V$.

We also recall from \cite[\S 3.5.4]{FL2a} that the map $\tilde\bgamma$ is $S^1$
equivariant with respect to the action on the complex fibers of
$\tN_{\ft}(\Xi,\fs)$ by scalar multiplication and the trivial action on the
base $\tM_{\fs}$, while $S^1$ acts on $\tsC_{\ft}$ through
\begin{equation}
\label{eq:DefnCircleMultWL}
\varrho_L: S^1\times V\to V,
\quad\text{where}\quad
\varrho_L(e^{i\theta}) = \id_{W}\oplus e^{i\theta}\,\id_{W\otimes L}.
\end{equation}
Therefore, we have an isomorphism of $\SO(3)$ bundles:
$$
(\tilde\bgamma\times\id_X)^*
\left(\tsC_{\ft}^*\times_{\sG_{\ft}} (i\underline{\RR}\oplus L)\right)
\cong
\tN_{\ft}(\Xi,\fs)\times_{\sG_{\fs}} (i\underline{\RR}\oplus L).
$$
We obtain $(\tilde\bgamma\times\id_X)^*\FF_{\ft}$ on the left above after
we take the $S^1$ quotient, with $S^1$ acting on $\tsC_{\ft}$ through
complex multiplication on $V$ and trivially on $\fg_{\ft}$. Given
$$
[A,\Phi,f,z]
=
[\tilde\bgamma(B,\Psi,\eta),f,z]
\in
\tsC_{\ft}^*\times_{\sG_{\ft}} (i\underline{\RR}\oplus L),
$$
and noting that
\begin{equation}
\label{eq:S1ActionRelation}
e^{i\theta}\,\id_V = \varrho(e^{i\theta})\varrho_L(e^{2i\theta}),
\quad
e^{i\theta} \in S^1,
\end{equation}
then we can identify the pull-back of the $S^1$ action: 
\begin{align*}
[e^{i\theta}(A,\Phi),f,z]
&=
[e^{i\theta}\tilde\bgamma(B,\Psi,\eta),f,z]
\\
&=
[\varrho(e^{i\theta})\varrho_L(e^{2i\theta})\tilde\bgamma(B,\Psi,\eta),f,z]
\quad\text{(Equation \eqref{eq:S1ActionRelation})}
\\
&=
[\varrho(e^{i\theta})\tilde\bgamma(B,\Psi,e^{2i\theta}\eta),f,z]
\quad\text{(see \cite[\S 3.5.4]{FL2a})}
\\
&=
[\tilde\bgamma(B,\Psi,e^{2i\theta}\eta),\varrho(e^{-i\theta})(f,z)]
\\
&=
[\tilde\bgamma(B,\Psi,e^{i2\theta}\eta),f,e^{2i\theta}z]
\quad\text{(see \cite[\S 3.4.2]{FL2a})}.
\end{align*}
The final equality follows from the observation that the action of
$s\in\sG_{\fs}$ induced by the embedding $\varrho:\sG_{\fs}\to\sG_{\ft}$,
the homomorphism $\Ad:\Aut(V)\to\su(V)$, the projection
$\su(V)\to\fg_{\ft}$, and the isomorphism $\fg_{\ft}\cong
i\underline{\RR}\oplus L$, is given by $(f,z)\mapsto (f,s^{-2}z)$; see
\cite[\S 3.4.2]{FL2a} for details.
\end{proof}

{}Lemma \ref{lem:UniversalReduction} shows that we can compute
$\bga^*p_1(\FF_{\ft}) \in H^4(\PP N_{\ft}(\Xi,\fs)\times X;\RR)$ once we
know the Chern class of the line-bundle component of the $\SO(3)$ bundle
\eqref{eq:UniversalBundleReduction}:

\begin{lem}
\label{lem:UniversalS1Bundle}
Continue the hypotheses of Lemma \ref{lem:UniversalReduction} and let $\nu$
be the cohomology class in Definition \ref{defn:nu}.
Then the complex line bundle,
\begin{equation}
\label{eq:UniversalReduced}
\tN^\eps_{\ft}(\Xi,\fs)\times_{\sG_{\fs}\times S^1} L
\to
\PP N_{\ft}(\Xi,\fs) \times X,
\end{equation}
has first Chern class
$$
\pi_{\PP N}^*\nu + 2(\pi_{\fs}\times\pi_X)^*c_1(\LL_{\fs}) + \pi_X^*c_1(L)
\in
H^2(\PP N_{\ft}(\Xi,\fs)\times X;\ZZ),
$$
where $\pi_{\PP n}$, $\pi_{\fs}$ and $\pi_X$ are the projections
from $\PP N_{\ft}(\Xi,\fs)\times X$ to
$\PP N_{\ft}(\Xi,\fs)$, $M_{\fs}$, and $X$, respectively,
and $\LL_{\fs}$ is the universal Seiberg-Witten line bundle
\eqref{eq:UniversalSWDefinition}.
\end{lem}

\begin{proof}
The projection $\tN^{\eps}_{\fs}(\Xi,\fs)\to \PP N_{\ft}(\Xi,\fs)$ is a
principal $\sG_{\fs}\times S^1$ bundle, where $S^1$ acts by scalar
multiplication on the fibers of $\tN_{\ft}(\Xi,\fs)$. One has an
isomorphism of $\sG_{\fs}\times S^1$ bundles
\begin{equation}
\label{eq:SplittNBundle}
\tN^{\eps}_{\ft}(\Xi,\fs)
\cong
N^{\eps}_{\ft}(\Xi,\fs)\times_{M_{\fs}} \tM_{\fs},
\end{equation}
defined, for $(B,\Psi)\in \tM_{\fs}$ and $(B,\Psi,\eta)\in
\tN_{\ft}^\eps(\Xi,\fs)$, by the map
$$
(B,\Psi,\eta) \mapsto \left([B,\Psi,\eta], (B,\Psi)\right),
$$
where the square brackets indicate equivalence modulo
$\sG_{\fs}$.  Applying the isomorphism
\eqref{eq:SplittNBundle} to the $\sG_{\fs}\times S^1$ bundle
\eqref{eq:UniversalReduced} yields an isomorphism of complex line bundles,
\begin{equation}
\label{eq:UniversalReduced2}
\tN^\eps_{\ft}(\Xi,\fs)\times_{\sG_{\fs}\times S^1} L
\cong
(N^{\eps}_{\ft}(\Xi,\fs)\times_{M_{\fs}}(\tM_{\fs}\times_{\sG_{\fs}} L))/S^1
\end{equation}
where, as in definition \eqref{eq:ReducibleUniversalAction}, an element
$s\in\sG_{\fs}$ acts on $\tM_{\fs}\times L$ by $(B,\Psi,z)\mapsto
(s(B,\Psi),s^{-2}z)$ and $S^1$ acts diagonally on
$N^\eps_{\ft}(\Xi,\fs)\times L$.

The isomorphism in \cite[Equation (6.68)]{FL2a} gives
\begin{equation}
\label{eq:TwistedSWUniversal}
\tM_{\fs}\times_{\sG_{\fs}}L
\cong
\LL_{\fs}^{\otimes 2}\otimes \pi_X^*L.
\end{equation}
Substituting the isomorphism \eqref{eq:TwistedSWUniversal} into equation
\eqref{eq:UniversalReduced2} yields an isomorphism of complex line bundles,
\begin{equation}
\tN^\eps_{\ft}(\Xi,\fs)\times_{\sG_{\fs}\times S^1} L
\cong
\label{eq:UniversalReduced3}
(N^{\eps}_{\ft}(\Xi,\fs)\times_{M_{\fs}}
\LL_{\fs}^{\otimes 2}\otimes\pi_X^*L)/S^1.
\end{equation}
The proof is completed by applying Lemma 3.27 in \cite{FL2a} to compute the
first Chern class of a fiber product with an $S^1$ action, the observation
that the $S^1$ action in
\eqref{eq:UniversalReduced3} is diagonal, and the fact that
$\nu = -c_1(N^{\eps}_{\ft}(\Xi,\fs))$.
\end{proof}

The reduction of $(\bga\times\id_X)^* \FF_{\ft}$ in Lemma
\ref{lem:UniversalReduction}, the computation in Lemma
\ref{lem:UniversalS1Bundle}, and the fact that $c_1(L) = c_1(\ft)-c_1(\fs)$
give the following expression for $(\bga\times\id_X)^*p_1(\FF_{\ft})$.

\begin{cor}
\label{cor:UniversalPontrjaginClassPullback}
Continue the hypotheses and notation of Lemmas \ref{lem:UniversalReduction}
and \ref{lem:UniversalS1Bundle}.
Then,
\begin{equation}
\label{eq:PulledBackUniversalPontrjagin}
\begin{aligned}
(\bga\times\id_X)^* p_1(\FF_{\ft})
&=
\left(\pi_{\PP N}^*\nu + 2(\pi_{\fs}\times\pi_X)^*c_1(\LL_\fs)
+ \pi_X^*(c_1(\ft)-c_1(\fs))\right)^2
\\
&\in
H^2(\PP N_{\ft}(\Xi,\fs)\times X;\ZZ).
\end{aligned}
\end{equation}
\end{cor}

We compute the pullbacks of the cohomology classes $\mu_{p}(\beta)$ in
$\sC_{\ft}^{*,0}/S^1$ to $\PP N_{\ft}(\Xi,\fs)$:

\begin{cor}
\label{cor:CohomologyOnReducibleLink}
Continue the hypotheses of Lemma \ref{lem:UniversalReduction}.
Let $\{\ga_i\}$ be a basis for $H_1(X;\ZZ)/\Tor$ and
let $\{\ga_i^*\}$ be the dual basis for $H^1(X;\ZZ)$.
Then if $x\in H_0(X;\ZZ)$ is the positive generator,
$\gamma\in H_1(X;\RR)$, $h\in H_2(X;\RR)$, and $[Y]\in H_3(X;\RR)$, the
pullbacks of the cohomology classes $\mu_{p}(\beta)$ in
$H^\bullet(\sC_{\ft}^{*,0}/S^1;\RR)$ by the embedding $\bga:\PP
N_{\ft}(\Xi,\fs) \embed \sC_{\ft}^{*,0}/S^1$ to cohomology classes
$\bgamma^*\mu_{p}(\beta)$ in $H^\bullet(\PP N_{\ft}(\Xi,\fs);\RR)$
are given by
\begin{equation}
\label{eq:PulledBackCohomology1}
\begin{aligned}
\bga^*\mu_{p}([Y])
&= \sum_{i=1}^{b_1(X)}
\left\langle(c_1(\fs)-c_1(\ft))\smile \gamma_i^*,[Y]\right\rangle
    \mu_{\fs}(\ga_i),
\\
\bga^*\mu_{p}(h)
& = {\textstyle{\frac{1}{2}}}\langle c_1(\fs)-c_1(\ft),h\rangle
       \left(2\mu_{\fs}(x)+\nu\right)
    -2\sum_{i<j}
    \langle\gamma_i^*\smile\gamma_j^*,h\rangle
      \mu_{\fs}(\ga_i\ga_j),
\\
\bga^*\mu_{p}(\gamma)
&= -\sum_{i=1}^{b_1(X)}\langle \gamma_i^*,\gamma\rangle
          \left(2\mu_{\fs}(x)+\nu\right))\smile\mu_{\fs}(\ga_i),
\\
\bga^*\mu_{p}(x)
& = -\textstyle{\frac{1}{4}}\left(2\mu_{\fs}(x)+\nu\right)^2,
\end{aligned}
\end{equation}
where we have written $\mu_{\fs}(\beta)$ for the pullback of this
class to $\PP N_{\ft}(\Xi,\fs)$.
\end{cor}

\begin{proof}
Recall from \cite[Lemma 2.24]{FL2a} that
$$
c_1(\LL_{\fs})
=
\mu_{\fs}(x)\times 1 + \sum_{i=1}^{b_1(X)} \mu_{\fs}(\ga_i)\times \ga^*_i
\in
H^2(\sC_{\fs}^0\times X;\RR).
$$
The identities \eqref{eq:PulledBackCohomology1} then follow from
equation \eqref{eq:PulledBackUniversalPontrjagin}, the definition
\eqref{eq:DefineMuMap} of the
cohomology classes $\mu_{p}(\beta)=-\quarter p_1(\FF_{\ft})/\beta$, and
standard computations (compare the proof of \cite[Proposition 5.1.21]{DK}).
\end{proof}

Finally, we compute the pullback of the class $\mu_{c}$ in
$\sC_{\ft}^{*,0}/S^1$ to $\PP N_{\ft}(\Xi,\fs)$:

\begin{lem}
\label{lem:PullbackOfc1}
Continue the hypotheses of Lemma \ref{lem:UniversalReduction}
and let $\nu$ be the cohomology class in Definition \ref{defn:nu}. Then
\begin{equation}
\label{eq:c1Pullback}
\bga^* \mu_{c} = \nu
\in
H^\bullet(\PP N_{\ft}(\Xi,\fs);\RR).
\end{equation}
\end{lem}

\begin{proof}
We compute $c_1(\bga^*\LL_\ft)$, where $\LL_\ft$ is the line bundle
\eqref{eq:DefineC1LineBundle} with $c_1(\LL_{\ft})=\mu_{c}$, so
$$
\LL_{\ft} = (\sC_{\ft}^{*,0}\times\CC)/S^1,
$$
with circle action given by
\begin{equation}
\label{eq:DetBundleRewrite}
\left( e^{i\theta}, ([A,\Phi], z)\right) \to
\left( e^{i\theta}[A,\Phi], e^{2i\theta}z\right)
=
\left( \varrho_L(e^{2i\theta})[A,\Phi], e^{2i\theta}z\right),
\end{equation}
where $\varrho_L$ is given by definition \eqref{eq:DefnCircleMultWL} and
the preceding equality follows from the relation
\eqref{eq:S1ActionRelation} between the circle actions.  The embedding
$\bga:N_{\ft}^\eps(\Xi,\fs)\to \sC_{\ft}^{*,0}$ is $S^1$ equivariant with
respect to scalar multiplication on the fibers of $N_{\ft}(\Xi,\fs)$ and
the action induced by $\varrho_L:S^1\times V\to V$ on $\sC^{*,0}_{\ft}$.
Hence, by definition of $\LL_{\ft}$ and the $S^1$ equivariance of
$\bgamma$, we obtain an isomorphism of complex line bundles
$$
\bgamma^*\LL_{\ft}
\cong
(N^{\eps}_\ft(\Xi,\fs)\times\CC)/S^1,
$$
where the circle acts by
$$
(e^{i\theta},([B,\Psi,\eta],z))
\mapsto
([B,\Psi,e^{i\theta}\eta),e^{i\theta}z).
$$
Therefore, $\bga^*\mu_{c}=\bga^*c_1(\LL_{\ft})=c_1(\bga^*\LL_{\ft})=\nu$,
as desired.
\end{proof}

Corollary \ref{cor:CohomologyOnReducibleLink} completes the proof
of Lemma \ref{lem:CyclesExtension}:

\begin{proof}[Proof of Assertions 2 (a), (b) and (c)
in Lemma \ref{lem:CyclesExtension}]
Because the class $\nu$ is
non-trivial on the fiber of  the projection
$\PP N_{\ft}(\Xi,\fs)\to M_{\fs}$,  the closures of the
geometric representatives $\sV(\beta)$, $\sV(\gamma)$, and $\sV(x)$
dual to the cohomology classes $\mu_{p}(\beta)$,
$\mu_{p}(\gamma)$, and $\mu_{p}(x)$ must contain each
point  in $\iota(M_{\fs})$.
Hence, the closures $\bar\sV(\beta)$,
$\bar\sV(\gamma)$, and $\bar\sV(x)$ contain $\iota(M_{\fs})$.
\end{proof}

\subsection{Euler and Segre classes}
\label{subsec:CharacteristicsOfRedNormal}
In \cite[Equation (3.48)]{FL2a}, the homology class of the link is given by
$$
[\bL_{\ft,\fs}] =e(\bga^*(\Xi/S^1))\cap [\PP N_{\ft}(\Xi,\fs)],
$$
where $\Xi/S^1$ is an obstruction bundle over an open neighborhood of
$\iota(M_{\fs})\subset\sC_{\ft}$,  
as defined in \cite[Theorem 3.19]{FL2a}.
In this section, we compute the Euler class of this obstruction bundle
and then compute the Segre classes of $N_{\ft}(\Xi,\fs)$ in order to
relate intersection pairings on $\PP N_{\ft}(\Xi,\fs)$ with pairings on
$M_{\fs}$. 

The obstruction bundle $\Xi$ is given by
$$
\Xi\cong \sU\times\CC^{r_{\Xi}} \to \sU,
$$
where (see \cite[Theorem 3.19]{FL2a})
$\sU$ is a neighborhood of $M_{\fs}$ in $\sC_{\ft}$ and the
$S^1$ action is given by, for $[A,\Phi]\in\sU$, $z\in \CC^{r_{\Xi}}$,
and $\varrho_L$ is the map \eqref{eq:DefnCircleMultWL},
\begin{equation}
\label{eq:ObstructionS1Action}
\left( [A,\Phi],z\right)\to
\left( \varrho_L(e^{i\theta})[A,\Phi], e^{i\theta}z\right).
\end{equation}
Because the embedding $\bga:N^\eps_{\ft}(\Xi,\fs)\to \sU\subset\sC_{\ft}$
is $S^1$ equivariant with respect to scalar multiplication on
the fibers of $N_{\ft}(\Xi,\fs)$ and the action induced by the
map $\varrho_L$ on $\sC_{\ft}$, we have an isomorphism:
$$
\bga^*\left(\Xi/S^1\right)
\cong
N^\eps_{\ft}(\Xi,\fs)\times_{(S^1,\times -1)}\CC^{r_{\Xi}},
$$
where the factor $-1$ indicates that the $S^1$ action in equation
\eqref{eq:ObstructionS1Action} is diagonal.  Thus, we can calculate the
Euler class of $\bga^*\Xi/S^1$:

\begin{lem}
\label{lem:EulerClassOfObstruction}
The vector bundle $\bga^*\Xi/S^1$ has Euler class
$e(\bga^*\Xi/S^1)=\nu^{r_\Xi}$, where $\nu$ is the cohomology class in
Definition \ref{defn:nu}, $r_\Xi = \rank_\CC\Xi$, and
$$
[\bL_{\ft,\fs}]=\nu^{r_\Xi}\cap [\PP N_{\ft}(\Xi,\fs)],
$$
where $\bL_{\ft,\fs}$ is given the complex orientation of
Definition \ref{defn:ComplexOrientation} and
$[\PP N_{\ft}(\Xi,\fs)]$ is given the orientation defined by the
orientation of $TM_{\fs}$ and the complex orientation of
the fibers.
\end{lem}

\begin{proof}
The obstruction section  vanishes transversely, so its zero
locus, $\bL_{\ft,\fs}$ is Poincar\'e dual to
$e(\bga^*\Xi/S^1)=\nu^{r_\Xi}$ by
\cite[Proposition 12.8]{BT}. Note that the top Chern class is the Euler class
associated to the complex orientation of the fibers of a complex vector
bundle \cite[Lemma 14.1 \& Definition, p. 158]{MilnorStasheff}.
\end{proof}

Although the definition of Segre classes is well-known
\cite[p. 69]{Fulton}, we include a definition here via the following lemma
in order to make our conventions clear.

\begin{lem}
\label{lem:Segre}
Let $N$ be a complex rank-$r$ vector bundle with
Chern classes $c_i=c_i(N)$ over an oriented, real $m$-dimensional
manifold $M$.  Let $N^\eps\subset N$ be the associated
$\eps$-sphere bundle.  Define Segre classes $s_i=s_i(N)\in
H^{2i}(X;\ZZ)$ by the relation
\begin{equation}
\label{eq:FormalInverseDefinition}
(1+c_1+c_2+\cdots+c_r)(s_0+s_1+\cdots)=1.
\end{equation}
Let $\pi:\PP(N)\rightarrow M$ be the
projectivization of $N$ and
$h$ the negative of the first Chern class of the bundle
$N^\eps\to \PP(N)$. Then, for any $\alpha\in H^{m-2i}(M;\ZZ)$,
\begin{equation}
  \label{eq:Segre}
\langle h^{r+i-1}\smile\pi^*\alpha,[\PP(N)]\rangle = \langle
s_i\smile \alpha,[M]\rangle,
\end{equation}
where $\PP(N)$ is given the orientation arising from that of $M$
and the complex orientation of the fibers of $\pi$.
\end{lem}

\begin{proof}
The cohomology ring of $\PP(N)$ is given by \cite[Equation (20.7)]{BT},
$$
H^\bullet(\PP(N);\ZZ) =
\pi^*H^\bullet(M;\ZZ)[h]/(h^r+\pi^*c_1h^{r-1}+ \cdots +\pi^*c_r).
$$
We then have $h^r=-\sum_{i=1}^r \pi^*c_ih^{r-i}$. Suppose
$\alpha\in H^m(M;\ZZ)$, so $i=0$ in the assertion of the lemma,
and $\alpha$ is dual to $\langle \alpha,[M]\rangle p \in
H_0(M;\ZZ)$, where $p$ is a point. Then
\begin{align*}
\langle h^{r-1}\smile\pi^*\alpha,[\PP(N)]\rangle
&=
\langle h^{r-1},[\pi^{-1}(p)]\rangle \langle \alpha,[M]\rangle
\\
&= \langle h^{r-1},\CC\PP^{r-1}]\rangle \langle \alpha,[M]\rangle
= \langle \alpha,[M]\rangle.
\end{align*}
Because $s_0(N)=1$ by equation \eqref{eq:FormalInverseDefinition}, equation
\eqref{eq:Segre} holds for $i=0$.  We now use induction on $i$ and consider
$\alpha\in H^{m-2i}(M,\ZZ)$:
\begin{align*}
\langle h^{r+i-1}\smile\pi^*\alpha,[\PP(N)]\rangle &=
-\Bigl\langle \Bigl(\sum_{j=1}^i
\pi^*c_jh^{r+i-1-j}\Bigr)\smile\pi^*\alpha, [\PP(N)]\Bigr\rangle
\\
&= -\sum_{j=1}^i\langle c_js_{i-j}\smile\alpha,[M]\rangle =
\langle s_{i}\smile\alpha,[M]\rangle.
\end{align*}
The last equality follows from the identity $(\sum_i c_i)(\sum_j
s_j) = 1$. This gives the desired relation.
\end{proof}

Next we compute the Segre classes of $N_{\ft}(\Xi,\fs)$. Recall from
\cite[Equation (3.72)]{FL2a} that
\begin{equation}
\label{eq:RedefineNormalIndices}
\begin{aligned}
n_s'(\ft,\fs)
&=
- (c_1(\ft)-c_1(\fs))^2 -\textstyle{\frac{1}{2}}(\chi+\sigma),
\\
n_s''(\ft,\fs)
&=
\textstyle{\frac{1}{8}}\left( (2c_1(\ft)-c_1(\fs))^2 -\sigma\right),
\end{aligned}
\end{equation}
where $n_s=n_s'+n_s''$ is the index of the elliptic ``normal deformation
operator'' for the Seiberg-Witten stratum $\iota(M_{\fs})\subset
\sM_{\ft}$.

\begin{lem}
\label{lem:SegreOfN}
Suppose that for all $\alpha,\alpha'\in H^1(X;\ZZ)$ one has
$\alpha\smile\alpha'=0$. Let $\ft$ be a \spinu structure with
$\ft=\fs\oplus \fs\otimes L$ and assume $M_{\fs}$ contains no zero-section
pairs. Then the bundle $N_{\ft}(\Xi,\fs)\to M_{\fs}$ has Segre classes
\begin{equation}
  \label{eq:SegreOfN}
s_i(N_{\ft}(\Xi,\fs))
=
\mu_{\fs}(x)^i\sum_{j=1}^i 2^{j} {\binom {-n_s'} {j}}{\binom
{-n_s''} {i-j}}, \qquad i =0,1,2,\dots.
\end{equation}
\end{lem}

\begin{proof}
  With the hypothesis on $H^1(X;\ZZ)$, Corollary 3.30 in
  \cite{FL2a} asserts that $N_{\ft}(\Xi,\fs)$ has total Chern class
$$
c(N_{\ft}(\Xi,\fs)) = (1+2\mu_{\fs}(x))^{n_s'}(1+\mu_{\fs}(x))^{n_s''}.
$$
As described in Lemma \ref{lem:Segre}, the total Segre class
$s=s_0+s_1+s_2+\cdots$ is the
formal inverse of the total Chern class $c=1+c_1+c_2+\cdots$, so
$$
s(N_{\ft}(\Xi,\fs)) = (1+2\mu_{\fs}(x))^{-n_s'}(1+\mu_{\fs}(x))^{-n_s''}.
$$
The lemma follows by computing the formal power series expansions for
the above expression, using equation (2) in Lemma
\ref{lem:CombinatorialIdentities} to simplify before multiplying the two
series.
\end{proof}

\begin{rmk}
\label{rmk:MoreGeneralSegreClasses}
The assumption $\alpha\smile\alpha'=0$ in Lemma \ref{lem:SegreOfN} is used
to simplify the expression for the Chern character of $N_{\ft}(\Xi,\fs)$
computed in \cite[Theorem 3.29]{FL2a}.  Without this assumption, the
universal expression of the Segre classes in terms of the Chern character
given in \cite[Equation (2.11)]{MacDonald} and Theorem 3.29 of \cite{FL2a}
still show that the Segre classes of $N_{\ft}(\Xi,\fs)$ are expressible
in terms of the $\mu_{\fs}$-classes, though not as explicitly.
\end{rmk}

\subsection{Computation of the intersection pairing}
\label{subsec:Computation}
We now compute intersection pairings with $\bL_{\ft,\fs}$ of the type
encountered in the cobordism formula \eqref{eq:CompactReductionFormula1}.
Some combinatorial factors appearing in this computation can
be expressed in terms of the {\em Jacobi
  polynomials\/} \cite[\S 8.96]{GR}, which are defined by
\begin{equation}
P^{a,b}_n(\xi)
=
{\frac{1}{ 2^n}}
\sum_k {\binom {a+n}{n-k}}{\binom {n+b}k}(\xi-1)^k(\xi+1)^{n-k},
\quad \xi\in\CC.
\label{defn:JacobiPolyn}
\end{equation}
Functional relations, relations with other special functions, and
the generating function for the Jacobi polynomials can be found in
\cite[pp. 1034--1035]{GR}.

\begin{thm}
\label{thm:DegreeZeroFormula}
Let $X$ be a four-manifold with $\alpha\smile\alpha'=0$ for every
$\alpha,\alpha'\in H^1(X;\ZZ)$, and let $\Omega$ be a homology orientation.
Let $\fs$ and $\ft$ be a \spinc and \spinu structure on $X$ for which
$\ft=\fs\oplus\fs\otimes L$, and assume $M_{\fs}$ contains no zero-section
pairs.  Give $\bL_{\ft,\fs}$ the complex orientation, determined by the
orientation for $M_{\fs}$ fixed by the homology orientation $\Omega$ as in
Definition \ref{defn:ReducibleOrient}.  Let $z\in\AAA(X)$
and let $\delta_{c}$ be a non-negative integer satisfying
$$
\deg(z) +2\delta_{c} = \dim(\sM_{\ft}^{*,0}/S^1)-1.
$$
If $z=z'Y$ for some $Y\in H_3(X;\ZZ)$ and $z'\in\AAA(X)$, then
\begin{equation}
\label{eq:Homology3Case}
\langle \mu_{p}(z)\smile \mu_{c}^{\delta_{c}},
[\bL_{\ft,\fs}]\rangle=0.
\end{equation}
If $z=x^{\delta_0}\vartheta h^{\delta_2}$,
where $h\in H_2(X;\RR)$, $\vartheta\in \Lambda^{\delta_1}(H_1(X;\RR))$, and
$x\in H_0(X;\ZZ)$ is the positive generator, then
$d_s(\fs)\equiv\delta_1\pmod{2}$ and if we set
\begin{equation}
  \label{eq:DefnDeltap1}
  2d = d_s(\fs)-\delta_1,
\end{equation}
then
\begin{equation}
\label{eq:GoodDelta1Case}
\begin{aligned}
\left\langle\mu_{p}(z)\smile \mu_{c}^{\delta_{c}},
[\bL_{\ft,\fs}]\right\rangle
&= (-1)^{\delta_0+\delta_1} 2^{-\delta_2-2\delta_0}
C_{\chi,\si}\left( \deg(z),\delta_{c},d_a(\ft),d_s(\fs),\delta_1\right)
\\
&\quad\times
\left\langle\mu_{\fs}(x^d\vartheta),
[M_{\fs}]\right\rangle\langle c_1(\fs)-c_1(\ft),h\rangle^{\delta_2},
\end{aligned}
\end{equation}
where
$$
C_{\chi,\si}\left( \deg(z),\delta_{c},d_a(\ft),d_s(\fs),\delta_1\right)
=
(-2)^{d}P^{a,b}_{d}(0),
$$
for
$$
a =\delta_{c}-d-1 \quad\text{and}\quad
b = \textstyle{\frac{1}{2}}\left(\deg(z)-d_a(\ft)-d_s(\fs)\right) 
-\textstyle{\frac{1}{4}}\left(\chi+\si\right).
$$
If $d=0$, then $C_{\chi,\si}=1$.
\end{thm}

\begin{rmk}
If $\delta_1> d_s(\fs)$, then the pairing $\langle\mu_{\fs}(x^d\vartheta),
[M_{\fs}]\rangle$ vanishes and so the pairing
\eqref{eq:GoodDelta1Case} also vanishes.
\end{rmk}

\begin{proof}[Proof of Theorem \ref{thm:DegreeZeroFormula}]
By the multilinearity of the pairing, we can assume that $\vartheta$ is a
monomial, $\vartheta=\ga_1\cdots\ga_{\delta_1}$, where
$\{\ga_1,\dots,\ga_{\delta_1}\}$ is a subset of a basis for
$H_1(X;\ZZ)/\Tor$.  Extend it to a basis and let $\{\ga_i^*\}$ be a dual
basis for $H^1(X;\RR)$, so $\langle \ga_i^*,\ga_j\rangle =\delta_{ij}$.

Suppose $z=z'Y$ for $Y\in H_3(X;\ZZ)$ and $z'\in\AAA(X)$.  The expression
for $\mu_{p}([Y])$ in equation \eqref{eq:PulledBackCohomology1} is a sum
of terms of the form
$$
\langle (c_1(\fs)-c_1(\ft))\smile\gamma_i^*,[Y]\rangle \mu_{\fs}(\ga_i)
=
\langle (c_1(\fs)-c_1(\ft))\smile\gamma_i^*\smile\PD [Y],[X]\rangle \mu_{\fs}(\ga_i),
$$
which vanish by our hypothesis on $H^1(X;\ZZ)$. This
yields identity \eqref{eq:Homology3Case}.

The integers $\delta_i$ and $\delta_{c}$ satisfy
\begin{equation}
  \label{eq:IndicesEquality1}
\begin{aligned}
2\delta_2+3\delta_1+4\delta_0+2\delta_{c}
&=
\deg(z)+2\delta_{c}
=\dim (\sM^{*,0}_{\ft}/S^1)-1
\\
&=
d_s(\fs)+2n_s'+2n_s''-2.
\end{aligned}
\end{equation}
Thus $\delta_1 \equiv d_s(\fs)\pmod 2$ and $d=\half(d_s(\fs)-\delta_1)$ is
an integer.  We use Corollary
\ref{cor:CohomologyOnReducibleLink}, the comments following Lemma 3.23 in
\cite{FL2a}, and Lemma \ref{lem:EulerClassOfObstruction} to write the
pairing in equation \eqref{eq:GoodDelta1Case} as
\begin{equation}
\label{eq:Pairing0}
\begin{aligned}
{}&\left\langle
    \mu_{p}(h^{\delta_2}\vartheta x^{\delta_0})
    \mu_{c}^{\delta_{c}},
[\bL_{\ft,\fs}]\right\rangle
\\
&= 2^{-\delta_2}\langle c_1(\fs)-c_1(\ft),h\rangle^{\delta_2}
\\
&\quad\times
\left\langle
(2\mu_{\fs}(x)+\nu)^{\delta_2}
    \left(
        \prod_{k=1}^{\delta_1}(-2\mu_{\fs}(x)-\nu)\mu_{\fs}(\ga_k)
    \right)
\left(-\textstyle{\frac{1}{4}} (2\mu_{\fs}(x)+\nu)^{2}\right)^{\delta_0}
\nu^{\delta_{c}},\right.
\\
&\quad\qquad\left.\nu^{r_\Xi}\cap [\PP N_{\ft}(\Xi,\fs)]\right\rangle.
\end{aligned}
\end{equation}
Write $C_1=(-1)^{\delta_1+\delta_0}2^{-\delta_2-2\delta_0}
\langle c_1(\fs)-c_1(\ft),h\rangle^{\delta_2}$
and define 
\begin{equation}
\label{eq:defndelta_p}
\delta_{p}
=
\textstyle{\frac{1}{2}}(\deg(z) -\delta_1)=\delta_2+\delta_1+2\delta_0. 
\end{equation}
The pairing \eqref{eq:Pairing0} is then equal to
\begin{equation}
\label{eq:Pairing1}
\begin{aligned}
{}&C_1 \left\langle (2\mu_{\fs}(x) +\nu)^{\delta_{p}}
\mu_{\fs}(\vartheta)
\nu^{\delta_{c}+r_\Xi}, [\PP N_{\ft}(\Xi,\fs)]\right\rangle
\\
&\quad = C_1\biggl\langle
\sum_{i=0}^{\delta_{p}}
2^i{\binom {\delta_{p}}{i}}
\mu_{\fs}(x^i\vartheta)
\nu^{\delta_{p}+\delta_{c}+r_\Xi-i} , [\PP N_{\ft}(\Xi,\fs)]\biggl\rangle.
\end{aligned}
\end{equation}
Then
$\mu_{\fs}(x^i\vartheta) \in
H^{d_s(\fs)-2(d-i)}(M_{\fs};\ZZ)$ and
\begin{align*}
\delta_{p}+\delta_{c}+r_\Xi-i
&=
\delta_2+\delta_1+2\delta_0+\delta_{c}+r_\Xi-i
\quad \text{(as $\delta_{p}=\delta_2+\delta_1+2\delta_0$)}
\\
&=
n_s'+n_s''+r_\Xi+\textstyle{\frac{1}{2}} (d_s(\fs)-\delta_1) -1-i
\quad \text{(Equation \eqref{eq:IndicesEquality1}) }
\\
&= n_s'+n_s''+r_\Xi+(d-i)-1
\quad \text{(as $d=\textstyle{\frac{1}{2}} (d_s(\fs)-\delta_1)$).}
\end{align*}
We use the preceding equation, the Segre class relation \eqref{eq:Segre}, and
the formulas \eqref{eq:SegreOfN} for the Segre classes
$s_i(N_{\ft}(\Xi,\fs))$ to calculate
\begin{equation}
\label{eq:SegreReduction}
\begin{aligned}
\left\langle \mu_{\fs}(x^i\vartheta)
\nu^{\delta_{p}+\delta_{c}+{r_\Xi}-i},
[\PP N_{\ft}(\Xi,\fs)]\right\rangle
&=
\biggl\langle\mu_{\fs}(x^i\vartheta)
    s_{d-i},[M_{\fs}]\biggr\rangle
\\
&=
\sum_{j=0}^{d-i} 2^j{\binom {-n_s'}{j}}{\binom {-n_s''}{d-i-j}}
\left\langle\mu_{\fs}(x^d\vartheta),
[M_{\fs}]\right\rangle.
\end{aligned}
\end{equation}
Writing $C_2 = C_1\langle\mu_{\fs}(x^d\vartheta),
[M_{\fs}]\rangle$ and substituting the formula \eqref{eq:SegreReduction}
into equation \eqref{eq:Pairing1} yields a simplified expression for that
pairing:
\begin{equation}
\begin{aligned}
\label{eq:Pairing2}
{}&C_2 \sum_{i=0}^{d}2^i{\binom
{\delta_{p}}{i}} \sum_{j=0}^{d-i}
2^{j}{\binom {-n_s'}{j}}{\binom {-n_s''}{d-i-j}}
\\
&\quad = C_2 \sum_{i=0}^{d}\sum_{j=0}^{d-i}2^{i+j} {\binom
{\delta_{p}}{i}} {\binom {-n_s'}{j}} {\binom
{-n_s''}{d-i-j}}.
\end{aligned}
\end{equation}
If we write $u=i+j$, then the pairing \eqref{eq:Pairing2} becomes
\begin{equation}
\label{eq:Pairing3}
\begin{aligned}
C_2\sum_{u=0}^{d}\sum_{i=0}^{u} 2^u
{\binom {\delta_{p}}{i}} {\binom {-n_s'}{u-i}} {\binom
{-n_s''}{d-u}}
&=
C_2 \sum_{u=0}^{d}2^u{\binom {-n_s''}{d-u}}
\sum_{i=0}^u {\binom {\delta_{p}}{i}}
{\binom {-n_s'}{u-i}}
\\
&=
C_2 \sum_{u=0}^{d} 2^u {\binom {-n_s''}{d-u}} {\binom
{\delta_{p}-n_s'}{u}},
\end{aligned}
\end{equation}
where the second equality follows from the Vandermonde convolution
identity (see equation (5) in Lemma \ref{lem:CombinatorialIdentities}).
Equation \eqref{eq:GoodDelta1Case} will follow from
equation \eqref{eq:Pairing3} and  Lemma
\ref{lem:JacobiPolynomialInterp} which expresses the last
sum in equation \eqref{eq:Pairing3} in terms of the Jacobi polynomial given
by $C_{\chi,\si}\left( \deg(z),\delta_{c},d_a(\ft),d_s(\fs),\delta_1\right)$.
This completes the proof of the theorem.
\end{proof}

\begin{rmk}
\label{rmk:LevelZeroMultConj}
The proof of Theorem \ref{thm:DegreeZeroFormula} implies that Conjecture
3.1 of \cite{FKLM} holds for level-zero reducibles, even without the
assumption that $\alpha\smile\alpha'=0$ for every $\alpha,\alpha'\in
H^1(X;\ZZ)$. We only used this condition on $H^1(X;\ZZ)$ in equation 
\eqref{eq:SegreReduction} in order to apply
the Segre class computations of Lemma \ref{lem:SegreOfN}.  If the condition
on $H^1(X;\ZZ)$ is is omitted, then---as noted in Remark
\ref{rmk:MoreGeneralSegreClasses}---the Segre classes  
can still be computed in terms of $\mu_{\fs}$-classes, though less
explicitly. In the general situation, the pairing 
\eqref{eq:GoodDelta1Case} could still be expressed in terms of a pairing of
$\mu_{\fs}$-classes with $M_{\fs}$ and,
when $c_1(\fs)$ is not a basic class, the pairing
\eqref{eq:GoodDelta1Case} would be zero.
\end{rmk}

Before proving the relation between the combinatorial expression in
equation \eqref{eq:Pairing3} and the Jacobi polynomial used at the end of
the proof of Theorem \ref{thm:DegreeZeroFormula}, it is convenient to
collect the combinatorial identities we shall need here.  For $a\in\RR$ and
$n\in\NN$, define
\begin{equation}
  \label{eq:Comb(a)Defn}
(a)_n = a(a+1)\cdots (a+n-1),\quad\text{and}\quad (a)_0=1.
\end{equation}
We then have:

\begin{lem}
\label{lem:CombinatorialIdentities} \cite[p. 9]{Luke}
Let $a, b\in\RR$ and let $k, n\in\NN$ be non-negative integers. Then:
\begin{align}
\tag{1} (a)_k &= (-1)^k(1-a-k)_k,
\\
\tag{2} {\binom {a}{k}} &= {\frac{(-1)^k(-a)_k}{k!}},
\\
\tag{3} (a)_{n-k} &= {\frac {(-1)^k(a)_n}{(1-a-n)_k}},\quad n\ge
k,
\\
\tag{4} (n-k)! &= {\frac {n!}{(-1)^k(-n)_k}},
\\
\tag{5} \sum_{i=0}^u {\binom {a}{i}}{\binom {b}{u-i}} &= {\binom
{a+b}{u}}.
\end{align}
\end{lem}

Identity (5) in Lemma
\ref{lem:CombinatorialIdentities} (the Vandermonde convolution identity)
follows by comparing coefficients in binomial expansions of the two sides
of the identity $(x+y)^a(x+y)^b=(x+y)^{a+b}$.

\begin{rmk}
\label{rmk:BinomCoeff}
Note that Identity (2) in Lemma \ref{lem:CombinatorialIdentities} allows
one to extend the definition of a binomial coefficient $\binom{n}{r}$ to
the case when $n\le 0$.  If $r=0$, the identity $\binom{n}{r}=1$ still
holds.
\end{rmk}

We now prove the relation between the combinatorial
expression in equation \eqref{eq:Pairing3} and the Jacobi polynomial.

\begin{lem}
\label{lem:JacobiPolynomialInterp}
Continue the notation of Theorem \ref{thm:DegreeZeroFormula}. Then
$$
\sum_{u=0}^{d} 2^u {\binom {-n_s''}{d-u}} {\binom
{\delta_{p}-n_s'}{u}}
=
(-2)^{d}P^{a,b}_{d}(0),
$$
where $\delta_p$ is defined in equation \eqref{eq:defndelta_p} and
$$
a =\delta_{c}-d-1 \quad\text{and}\quad
b = \textstyle{\frac{1}{2}}\left(\deg(z)-d_a(\ft)-d_s(\fs)\right)
-\textstyle{\frac{1}{4}}\left(\chi+\si\right).
$$
\end{lem}

\begin{proof}
  We first recall that the {\em hypergeometric functions\/} \cite[\S
  9.10]{GR} are defined by
\begin{equation}
\HG(a,b;c;\xi)
=
\sum_{k=0}^\infty {\frac {(a)_k(b)_k}{(c)_k k!}}\xi^k,\quad \xi\in\CC.
\label{defn:HyperGeomFunction}
\end{equation}
We shall use the following identities \cite[Equation (23), p. 40 \&
Equation (2), p. 274]{Luke}:
\begin{equation}
\label{eq:FunctionRelations}
\begin{aligned}
\HG(-m,b;c;\xi)&= {\frac {(b)_m(-1)^m\xi^m}{(c)_m}}
\HG(-m,1-m-c;1-m-b;\xi^{-1}),
\\
P^{a,b}_n(\xi)&= {\frac {(-1)^n(b+1)_n}{n!}}
\HG(-n,n+a+b+1;b+1;\textstyle{\frac{1}{2}}(1+\xi)).
\end{aligned}
\end{equation}
By equation (2) of Lemma \ref{lem:CombinatorialIdentities},
the combinatorial expression in equation \eqref{eq:Pairing3}
can be written as
\begin{equation}
\label{eq:Combinatorial1}
\begin{aligned}
C(n_s'',n_s',\delta_{p},d)
&:=
\sum_{u=0}^{d} 2^u {\binom {-n_s''}{d-u}} {\binom
{\delta_{p}-n_s'}{u}}
\\
&=
\sum_{u=0}^{d} 2^u (-1)^{d-u}{\frac
{(n_s'')_{d-u}}{(d-u)!}} (-1)^u
{\frac {(n_s'-\delta_{p})_u}{u!}}.
\end{aligned}
\end{equation}
Applying equation (3) of Lemma \ref{lem:CombinatorialIdentities} to
$(n_s'')_{d-u}$ and equation (4) to $(d-u)!$ yields
\begin{equation}
\label{eq:Combinatorial1a}
\begin{aligned}
C(n_s'',n_s',\delta_{p},d)
&=
(-1)^{d}\sum_{u=0}^{d}2^u {\frac {(-1)^u
(n_s'')_{d}(-1)^u(-d)_u(n_s'-\delta_{p})_u}
{(1-n_s''-d)_u(d)!u!}}
\\
&= {\frac {(-1)^{d}(n_s'')_{d}}{(d)!}} \sum_{u=0}^{d}
{\frac {(-d)_u(n_s'-\delta_{p})_u}{(1-n_s''-d)_uu!}}2^u
\\
&= {\frac {(-1)^{d}(n_s'')_{d}}{(d)!}}
\HG(-d,n_s'-\delta_{p};1-n_s''-d;2).
\end{aligned}
\end{equation}
Substituting the first identity in equation \eqref{eq:FunctionRelations}
into equation \eqref{eq:Combinatorial1a} gives
\begin{equation}
\label{eq:Combinatorial2}
C(n_s'',n_s',\delta_{p},d)
=
{\frac{(n_s'')_{d}2^{d}(n_s'-\delta_{p})_{d}}
{(d)!(1-n_s''-d)_{d}}} \HG
(-d,n_s'';\delta_{p}-n_s'-d+1;\textstyle{\frac{1}{2}}).
\end{equation}
By substituting equation (1) from Lemma \ref{lem:CombinatorialIdentities}
into equation \eqref{eq:Combinatorial2} we obtain
\begin{equation}
\label{eq:Combinatorial3}
C(n_s'',n_s',\delta_{p},d)
=
{\frac {2^{d}(-1)^{d}(n_s'-\delta_{p})_{d}} {(d)!}} \HG
(-d,n_s'';\delta_{p}-n_s'-d+1;\textstyle{\frac{1}{2}}).
\end{equation}
Applying the second identity in \eqref{eq:FunctionRelations} to equation
\eqref{eq:Combinatorial3} yields
\begin{equation}
\label{eq:Combinatorial4}
C(n_s'',n_s',\delta_{p},d)
=
{\frac {2^{d}(n_s'-\delta_{p})_{d}(d)!}
{(d)!(\delta_{p}-n_s'-d+1)_{d}}}
P^{n_s''+n_s'-\delta_{p}-1,\delta_{p}-n_s'-d}_{d}(0).
\end{equation}
By applying equation (1) from Lemma \ref{lem:CombinatorialIdentities}, we
can then simplify the right-hand side of equation \eqref{eq:Combinatorial4}
to give
\begin{equation}
\label{eq:Combinatorial5}
C(n_s'',n_s',\delta_{p},d)
=
(-2)^{d}P^{n_s''+n_s'-\delta_{p}-1,\delta_{p}-n_s'-d}_{d}(0).
\end{equation}
Then, the equality
$$
d_s(\fs) +2n_s' +2n_s'' -2
=
\dim (\sM^{*,0}_{\ft}/S^1) -1
=
\deg(z) +2\delta_{c}
=
2\delta_{p}+\delta_1+2\delta_{c}
$$
and $d_s(\fs)=2d +\delta_1$ imply
\begin{equation}
\label{eq:IndexEquality1}
n_s''+n_s'-\delta_{p}-1=\delta_{c}-d-1.
\end{equation}
The definition \eqref{eq:RedefineNormalIndices} of $n'_s$ implies that
$$
n'_s = \textstyle{\frac{1}{2}}d_a(\ft) + \textstyle{\frac{1}{4}}(\chi+\si),
$$
which, together with the identities
$d=\frac{1}{2}(d_s(\fs)-\delta_1)$ and
$\delta_{p}=\frac{1}{2}(\deg(z)-\delta_1)$, yields
\begin{equation}
\label{eq:IndexEquality2}
\delta_{p}-n_s'-d
=
\textstyle{\frac{1}{2}}\left(\deg(z)-d_a(\ft)-d_s(\fs)\right) 
-\textstyle{\frac{1}{4}}\left(\chi+\si\right).
\end{equation}
Substituting equations \eqref{eq:IndexEquality1}
and \eqref{eq:IndexEquality2} into equation \eqref{eq:Combinatorial5} then
completes the proof.
\end{proof}

\subsection{A blow-up formula for Seiberg-Witten link pairings}
\label{subsec:LinkBlowUp}
Our formula \eqref{eq:GoodDelta1Case} in Theorem
\ref{thm:DegreeZeroFormula} for pairings with Seiberg-Witten links
$[\bL_{\ft,\fs}]$ only applies when $M_{\fs}$ contains no zero-section
pairs. In the same vein, our formula \eqref{eq:CompactReductionFormula1} in
Theorem \ref{thm:CompactReductionFormula} for $\#(\bar \sV(z)\cap\bar
M_\kappa^w)$ in terms of pairings with $[\bL_{\ft,\fs}]$ only applies when
$w_2(\ft)\equiv w\pmod{2}$ and $w\pmod{2}$ is good; when
$\ft=\fs\oplus\fs'$ and $w_2(\ft)$ is good, then $M_{\fs}$ contains no
zero-section pairs. The definition
\eqref{eq:DefineDInvarBlowUp} of the Donaldson invariants $D_X^w(z)$
incorporates the blow-up formula in order to remove any such constraint on
$w$. Therefore, we derive a ``blow-up'' formula for the
Seiberg-Witten link pairings which, together with equations
\eqref{eq:CompactReductionFormula1} and \eqref{eq:GoodDelta1Case},
will allow us to compute $D_X^w(z)$ for arbitrary $w\in H^2(X;\ZZ)$.

As before, we let $\tilde X=X\#\overline{\CC\PP}^2$ be the blow-up of $X$
with exceptional class $e\in H_2(\tilde X;\ZZ)$ and denote its
Poincar\'e dual by $\PD[e]\in H^2(\tilde X;\ZZ)$. We first need to relate
\spinc structures on $X$ with those on $\tilde X$.  Because
$\overline{\CC\PP}^2$ is simply-connected, the following map is a
bijection:
$$
\Spinc(\overline{\CC\PP}^2)
\to
\{(2k-1)\PD[e]: k\in\ZZ\}\subset H^2(\overline{\CC\PP}^2;\ZZ),
\quad
\fs \mapsto c_1(\fs).
$$
Let $\fs_{2k-1}$ denote the spinc structure on $\overline{\CC\PP}^2$ with
$c_1(\fs_{2k-1})=(2k-1)\PD[e]$. By the discussion in \cite[\S
12.4]{SalamonSWBook}, a \spinc structure $\fs$ on $X$ and $\fs_{2k-1}$ on
$\overline{\CC\PP}^2$ can be spliced together to yield a \spinc structure
$\fs\#\fs_{2k-1}$ on $\tilde X$ with
\begin{equation}
  \label{eq:C1BlownUpSpinStructure}
c_1(\fs\#\fs_{2k-1})
=
c_1(\fs)+(2k-1)\PD[e].
\end{equation}
Moreover, every \spinc structure on $\tilde X$ can be realized in this
way. The dimensions of the Seiberg-Witten moduli spaces are related by
\begin{equation}
\label{eq:DimBlownUpModuli}
\begin{aligned}
d_s(\fs\#\fs_{2k-1})
&=
\textstyle{\frac{1}{4}}
\left( c_1(\fs\#\fs_{2k-1})^2 -2\tilde\chi-3\tilde\sigma\right)
\\
&=
d_s(\fs)-k(k-1),
\end{aligned}
\end{equation}
where $\tilde\chi=\chi+1$ denotes the Euler characteristic of $\tilde X$
and $\tilde\sigma=\sigma-1$ is the signature of $\tilde X$.  We now define
a \spinu structure $\tilde\ft$ on $\tilde X$ related to a
\spinu structure $\ft$ on $X$, and relate reducible
$\PU(2)$ monopoles in $\bar\sM_{\tilde\ft}$ to those in
$\bar\sM_{\ft}$.

\begin{lem}
\label{lem:BlowUpSpinu}
Let $\ft$ be a \spinu structure on $X$ with the property that reducible
$\PU(2)$ monopoles in $\bar\sM_{\ft}$ appear only in the top level
$\sM_{\ft}$.  For $\tilde X=X\#\overline{\CC\PP}^2$, let $e\in H_2(\tilde
X;\ZZ)$ be the exceptional class, and let $\PD[e]\in H^2(\tilde X;\ZZ)$ be its
Poincar\'e dual. Then we have:
\begin{enumerate}
\item
There is a \spinu structure $\tilde\ft$ on $\tilde X$ satisfying
$$
c_1(\tilde\ft)=c_1(\ft), \quad p_1(\tilde\ft)=p_1(\ft)-1, \quad\text{and}\quad
  w_2(\tilde\ft)\equiv w_2(\ft)+\PD[e]\pmod{2}.
$$
\item 
The reducible $\PU(2)$ monopoles in $\bar\sM_{\tilde\ft}$ appear only in
  the top level, $\sM_{\tilde\ft}$, and are defined by \spinc structures
  $\fs^\pm$ on $\tilde X$ with $c_1(\fs^\pm)=c_1(\fs)\pm \PD[e]$, where
  $\iota(M_{\fs})\subset\sM_{\ft}$.
\item
Suppose we relax the assumption that reducible $\PU(2)$ monopoles in
$\bar\sM_{\ft}$ appear only in the top level $\sM_{\ft}$, to the assumption
that reducible $\PU(2)$ monopoles in $\bar\sM_{\ft}$ with non-trivial
Seiberg-Witten functions appear only in the top level $\sM_{\ft}$. Then
reducible $\PU(2)$ monopoles in $\bar\sM_{\tilde\ft}$ with non-trivial
Seiberg-Witten functions appear only in the top-level $\sM_{\tilde\ft}$.
\end{enumerate}
\end{lem}

\begin{proof}
  Suppose $\ft=(\rho,V)$. By Lemma 2.3 in \cite{FL2a}, we may assume
  $V=W\otimes E$, where $\fs=(\rho,W)$ is a \spinc structure on $X$ and
  $E\to X$ is a complex rank-two vector bundle. Let $\tilde E\to\tilde X$
  be the complex, rank-two bundle with $c_1(\tilde E)=c_1(E)+\PD[e]$
  and $c_2(\tilde E)=c_2(E)$, and let $\tilde\fs=(\tilde\rho,\tilde W)$ be
  the \spinc structure on $\tilde X$ with $c_1(\tilde\fs)=c_1(\fs)-\PD[e]$.
  Then set $\tilde V = \tilde W\otimes\tilde E$ and $\tilde\ft =
  (\tilde\rho,\tilde V)$, and observe that $\tilde\ft$ has the desired
  characteristic classes.

  By the neck-stretching argument described in \cite{FSTurkish}, the only
  non-empty Seiberg-Witten moduli spaces on $\tilde X$ are defined by \spinc
  structures $\fs\# \fs_{2k-1}$, where $M_{\fs}$ is non-empty. Since
  $\iota(M_{\fs})\subset\sM_{\ft}$ by hypothesis, equation
  \eqref{eq:SplittingSpinu} implies that $c_1(\fs)$ obeys
$$
\left( c_1(\fs) - c_1(\ft)\right)^2
=
p_1(\ft).
$$
To see which Seiberg-Witten moduli spaces $\iota(M_{\fs\# \fs_{2k-1}})$
can be contained in $\bar\sM_{\ft}$ and in which level, we need to check
the corresponding equation for $c_1(\fs\# \fs_{2k-1})$.
Equation \eqref{eq:C1BlownUpSpinStructure} and the relations between the
characteristic classes of $\ft$ and $\tilde\ft$ yields
\begin{align*}
\left( c_1(\fs\#\fs_{2k-1}) - c_1(\tilde\ft)\right)^2
&=
\left( c_1(\fs) + (2k-1)\PD[e] - c_1(\ft)\right)^2
\\
&=
p_1(\ft) -(2k-1)^2
\\
&=
p_1(\tilde\ft)-4k(k-1).
\end{align*}
Restricted to integers, the function $-4k(k-1)$ takes its maximum value at
$k=0$ and $k=1$. Hence, only the spaces $\iota(M_{\fs\#\fs_{2k-1}})$ with
$k=0,1$, appear in $\bar\sM_{\tilde\ft}$, as all other \spinc structures
$\fs\#\fs_{2k-1}$ would require an $\SO(3)$ bundle with Pontrjagin class
smaller than $p_1(\tilde\ft)$.
\end{proof}

\begin{thm}
\label{thm:SWBlowUpFormula}
\cite[Theorem 3.2]{OzsvathSzaboAdjunct}, \cite[Theorem 1.4]{FSTurkish} Let
$X$ be a four-manifold, and let $\tilde X = X\#\overline{\CC\PP}^2$ denote
its blow-up, with exceptional class $e\in H_2(\tilde X;\ZZ)$. If
$b_2^+(X)>1$, then for each \spinc structure $\tilde\fs$ on $\tilde X$ with
$d_s(\tilde \fs)\geq 0$ and each $z\in \BB(X)\cong\BB(\tilde X)$, we have
  \begin{equation}
    \label{eq:SWBlowUpFormula}
SW_{\tilde X,\tilde\fs}(z)
=
SW_{X,\fs}(x^mz),
  \end{equation}
  where $\fs$ is the \spinc structure induced on $X$ by restriction, and
  $2m = d_s(\fs)-d_s(\tilde\fs)$. If $b_2^+(X)=1$ and
  $c_1(\fs)-\La$ is not torsion, there is a one-to-one
  correspondence between $(c_1(\fs)-\La)$-chambers in
  the positive cone of $H^2(X;\RR)$ and $(c_1(\tilde\fs)-\La)$-chambers in
  the positive cone of $H^2(\tilde X;\RR)$, and the above
  relation holds provided both invariants are calculated in related chambers.
\end{thm}

\begin{rmk}
The presence of the class $\La$ in the hypotheses of Theorem
\ref{thm:SWBlowUpFormula} 
when describing the chambers arises because of the nature of the fixed
perturbation used in our definition of the Seiberg-Witten moduli spaces; see
the discussion in \S \ref{subsec:SWinvariants}.
\end{rmk}

Lemma \ref{lem:BlowUpSpinu}, the blow-up formula for Seiberg-Witten
invariants (Theorem \ref{thm:SWBlowUpFormula}), and Theorem
\ref{thm:CompactReductionFormula} then yield the following ``blow-up''
formula for Seiberg-Witten link pairings:

\begin{prop}
\label{prop:DegreeZeroBlowUp}
Continue the hypotheses and notation of Theorem \ref{thm:DegreeZeroFormula}
leading to equation \eqref{eq:GoodDelta1Case},
except we omit the requirement that $M_{\fs}$ contains no zero-section pairs
and define $z$ as given below.
Let $\tilde X=X\#\overline{\CC\PP}^2$ be the blow-up, let $e\in H_2(\tilde
X;\ZZ)$ be the exceptional class, and let $\PD[e]\in H^2(\tilde X;\ZZ)$ be its
Poincar\'e dual. Let $z = x^{\delta_0}\vartheta h^{\delta_2-k}
\in\AAA(X)\subset\AAA(\tilde X)$.
Let $\ft$ and $\tilde\ft$ be related \spinu structures on $X$ and $\tilde
X$, respectively, as in Lemma \ref{lem:BlowUpSpinu}. Then, for $k$ even,
\begin{equation}
\label{eq:DegreeZeroBlowUp}
\begin{aligned}
{}&(-1)^{o_{\tilde\ft}(w+\PD[e],\fs^+)}
\left\langle \mu_{p}(e^{k+1}z)\smile \mu_{c}^{\delta_{c}},
[\bL_{\tilde\ft,\fs^+}]\right\rangle
\\
+ &(-1)^{o_{\tilde\ft}(w+\PD[e],\fs^-)}
\left\langle \mu_{p}(e^{k+1}z)\smile \mu_{c}^{\delta_{c}},
[\bL_{\tilde\ft,\fs^-}]\right\rangle
\\
&=(-1)^{o_{\ft}(w,\fs)+\delta_0+\delta_1} 2^{-\delta_2-2\delta_0}
C_{\chi,\si}\left(\deg(z)+2k,\delta_{c},d_a(\ft),d_s(\fs),\delta_1\right)
\\
&\qquad
\times SW_{X,\fs}(x^d\vartheta)\langle c_1(\fs)-c_1(\ft),h\rangle^{\delta_2-k},
\end{aligned}
\end{equation}
while the left-hand side is zero if $k$ is odd or $z$ is replaced by $z'Y$
and $Y\in H_3(X;\ZZ)$.
\end{prop}

\begin{proof}
  The result follows by applying equation
  \eqref{eq:GoodDelta1Case} in Theorem \ref{thm:DegreeZeroFormula} to
  the links $\bL_{\tilde\ft,\fs^\pm}$ of $M_{\fs^\pm}(\tilde X)$ in
  $\sM_{\tilde\ft}(\tilde X)$, together with the following observations.

  The vanishing result in the case $z=z'Y$ follows immediately from
  equation   \eqref{eq:Homology3Case}.

  Because $d_s(\fs^\pm)=d_s(\fs)$ by equation \eqref{eq:DimBlownUpModuli},
  and $d_a(\tilde\ft)=d_a(\ft)+2$ by equation \eqref{eq:AsdDimDiracIndex}
  (noting that $p_1(\ft)=p_1(\tilde\ft)-1$ from Lemma \ref{lem:BlowUpSpinu}),
  and $\tilde\chi+\tilde\si=\chi +\si$, we have:
\begin{equation}
\label{eq:ConstantBlowUp}
C_{\tilde\chi,\tilde\si}
    (\deg(ze^{2k+1}),\delta_{c},d_a(\tilde\ft),d_s(\fs^\pm),\delta_1)
=
C_{\chi,\si}
    \left(\deg(z)+2k,\delta_{c},d_a(\ft),d_s(\fs),\delta_1\right).
\end{equation}
The proof (see \cite[\S 4]{FSTurkish}) of the blow-up formula, Theorem
\ref{thm:SWBlowUpFormula}, gives an identity
$$
\left\langle\mu_{\fs^+}(x^d\vartheta),
[M_{\fs^+}(\tilde X)]\right\rangle
=
\left\langle\mu_{\fs^-}(x^d\vartheta),
[M_{\fs^-}(\tilde X)]\right\rangle,
$$
and thus our definition \eqref{eq:SWFunctionDefn} of the Seiberg-Witten
invariants yields
$$
\left\langle\mu_{\fs^\pm}(x^d\vartheta),
[M_{\fs^\pm}(\tilde X)]\right\rangle
=
SW_{X,\fs}(x^d\vartheta).
$$
Noting that $c_1(\tilde\ft)=c_1(\ft)$ by Lemma \ref{lem:BlowUpSpinu},
the product
$$
\langle c_1(\fs^\pm)-c_1(\tilde \ft),e\rangle^{k+1}
\langle c_1(\fs^\pm)-c_1(\ft),h\rangle^{\delta_2-k}
$$
appearing in equation \eqref{eq:GoodDelta1Case} can be simplified to
\begin{equation}
\label{eq:AdditionalNegOneFactor}
\langle c_1(\fs^\pm)-c_1(\tilde\ft),e\rangle^{k+1}
\langle c_1(\fs^\pm)-c_1(\tilde\ft),h\rangle^{\delta_2-k}
=
(\mp 1)^{k+1}\langle c_1(\fs)-c_1(\ft),h\rangle^{\delta_2-k}.
\end{equation}
{}From its definition \eqref{eq:OrientationFactor}, the
orientation term is given by
\begin{equation}
\label{eq:OrientationBlowUp}
\begin{aligned}
o_{\tilde\ft}(w+\PD[e],\fs^\pm)
&=
\textstyle{\frac{1}{4}}
\left( w+\PD[e] - c_1(\tilde\ft) +c_1(\fs^\pm)\right)^2
\\
&=
\textstyle{\frac{1}{4}}
\left( w-c_1(\ft) + c_1(\fs)\right)^2 -\quarter(1\pm 1)^2
\\
&=
o_{\ft}(w,\fs)  - \textstyle{\frac{1}{4}}(1\pm 1)^2.
\end{aligned}
\end{equation}
Hence, applying equation \eqref{eq:GoodDelta1Case} to the pair
$\bL_{\tilde\ft,\fs^\pm}$, using equations
\eqref{eq:AdditionalNegOneFactor} and \eqref{eq:OrientationBlowUp} to
compute the sign differences between the pairings with
$\bL_{\tilde\ft,\fs^+}$ and $\bL_{\tilde\ft,\fs^-}$,
and using equation \eqref{eq:ConstantBlowUp} to
relate the constants yields
\begin{align*}
{}&(-1)^{o_{\tilde\ft}(w+\PD[e],\fs^+)}
\left\langle \mu_{p}(e^{k+1}z)\smile \mu_{c}^{\delta_{c}},
[\bL_{\tilde\ft,\fs^+}]\right\rangle
\\
+ &(-1)^{o_{\tilde\ft}(w+\PD[e],\fs^-)}
\left\langle \mu_{p}(e^{k+1}z)\smile \mu_{c}^{\delta_{c}},
[\bL_{\tilde\ft,\fs^-}]\right\rangle
\\
&=\left((-1)(-1)^{k+1} + 1\right)
(-1)^{o_{\ft}(w,\fs)+\delta_0+\delta_1} 2^{-\delta_2-1-2\delta_0}
C_{\chi,\si}\left(\deg(z)+2k,\delta_{c},d_a(\ft),d_s(\fs),\delta_1\right)
\\
&\qquad
\times SW_{X,\fs}(x^d\vartheta)
\langle c_1(\fs)-c_1(\ft),h\rangle^{\delta_2-k}.
\end{align*}
Since $(-1)(-1)^{k+1} + 1 = 2$ if $k$ is even and is zero if $k$ is odd,
the preceding equation reduces to the desired formula
\eqref{eq:DegreeZeroBlowUp}.
\end{proof}

\subsection{Proofs of main theorems}
\label{subsec:ProofsOfMainTheorems}
Combining Theorems \ref{thm:CompactReductionFormula} and
\ref{thm:DegreeZeroFormula} and a brief discussion
will complete the proofs of Theorems \ref{thm:DSWSeriesRelation},
\ref{thm:Main}, \ref{thm:FLthm}, and Corollary \ref{cor:FLthm}:

\begin{proof}[Proof of Theorem \ref{thm:Main}]
By hypothesis we have $w-\Lambda\equiv w_2(X)\pmod{2}$ and the invariants
$D_X^w(z)$ are zero unless $\deg(z)$ obeys the constraint \eqref{eq:Mod8}.
Let $p\in\ZZ$ satisfy $p\equiv w^2\pmod{4}$ and recall that---see the
paragraph following equation (2.20) in \cite{FL2a}---we may 
choose a \spinu structure $\ft$ on $X$ for which
$$
c_1(\ft) = \Lambda,
\quad
p_1(\ft) = p,
\quad\text{and}\quad
w_2(\ft) \equiv w \pmod{2}.
$$
Then $d_a(\ft) = -2p - \frac{3}{2}(\chi+\sigma)$ and $n_a(\ft) =
\frac{1}{4}(p+\Lambda^2-\sigma)$, by equation \eqref{eq:AsdDimDiracIndex}.

{}From equation \eqref{eq:LowerLevelReducible}, a reducible $\PU(2)$
monopole in $\bar\sM_{\ft}$ defined by a \spinc structure $\fs$ lies in the
level $\sM_{\ft_\ell}\times\Sym^\ell(X)$, where $\ell=\ell(\ft,\fs)$ and
$$
4\ell(\ft,\fs)
=
(c_1(\ft)-c_1(\fs))^2 - p_1(\ft).
$$
But $c_1(\ft)=\Lambda$ and $p_1(\ft) =
-\frac{1}{2}d_a(\ft)-\frac{3}{4}(\chi+\sigma)$ by equation
\eqref{eq:AsdDimDiracIndex}, so the definition
\eqref{eq:SWInTopLevelFunction} of
$r(\Lambda)$ and $r(\Lambda,c_1(\fs))$ implies that
\begin{equation}
  \label{eq:LevelOfReducible}
\begin{aligned}
4\ell(\ft,\fs)
&=
\textstyle{\frac{1}{2}}d_a(\ft) + (\Lambda-c_1(\fs))^2
+ \textstyle{\frac{3}{4}}(\chi+\sigma)
\\
&=
\textstyle{\frac{1}{2}}d_a(\ft) - r(\Lambda,c_1(\fs))
\\
&\leq
\textstyle{\frac{1}{2}}d_a(\ft) - r(\Lambda).
\end{aligned}
\end{equation}
Hence, when $d_a\leq 2r(\Lambda)$, the strata $\iota(M_{\fs})$ with non-trivial
Seiberg-Witten functions $\SW_{X,\fs}$ can only appear in the top level
$\sM_{\ft}$ of $\bar\sM_{\ft}$ (if at all), where they correspond to
splittings $\ft=\fs\oplus\fs'$.

{}From equation \eqref{eq:AsdDimDiracIndex}, the stratum $\iota(M_\kappa^w)$
has real codimension $2n_a(\ft)$ in $\sM_{\ft}$, with
\begin{equation}
  \label{eq:ASDCodimension}
\begin{aligned}
4n_a(\ft)
&=
p_1(\ft)+c_1(\ft)^2-\sigma
\\
&=
-\textstyle{\frac{1}{2}}d_a(\ft)
- \textstyle{\frac{3}{4}}(\chi+\sigma)+\Lambda^2-\sigma
\\
&=
-\textstyle{\frac{1}{2}}d_a(\ft)
- \textstyle{\frac{1}{4}}(7\chi+11\sigma)+\Lambda^2+\chi+\sigma
\\
&=
-\textstyle{\frac{1}{2}}d_a(\ft)  + i(\Lambda),
\end{aligned}
\end{equation}
where the second equality follows from equation \eqref{eq:AsdDimDiracIndex}
and the final one by definition \eqref{eq:PositiveDiracIndexFunction} of
$i(\Lambda)$.  Thus, our hypotheses on $\deg(z)$ and $\Lambda$ imply that
$n_a(\ft)>0$ in Cases (a) and (b) below (where $d_a=\deg(z)$), and also in
Case (c) (where $d_a<\deg(z)$), recalling that $\deg(z)$ is denoted by
$2\delta$, for $\delta\in \frac{1}{2}\ZZ$, in the hypotheses of Theorem
\ref{thm:Main}.

Therefore, provided $w\pmod{2}$ is good, we can apply Theorem
\ref{thm:CompactReductionFormula} to the cobordism $\sM_{\ft}^{*,0}$. To
eliminate this last constraint on $w$ when $b_2^+(X)>1$, we shall instead
apply Theorem
\ref{thm:CompactReductionFormula} to the cobordism $\sM_{\tilde\ft}^{*,0}$,
where $\tilde\ft$ is the related \spinu structure on the blow-up $\tilde X
= X\#\overline{\CC\PP}^2$ produced by Lemma \ref{lem:BlowUpSpinu}. When
$b_2^+(X)=1$, we assume that $w\equiv w_2(X)-\Lambda\pmod{2}$ is good so that
the Donaldson and Seiberg-Witten invariants are well-defined in this case
(see \S \ref{subsubsec:DefnDInvariants} and \S \ref{subsec:SWinvariants})
and $\sM_{\ft}$ contains no zero-section pairs.

{}From Lemma \ref{lem:BlowUpSpinu}, Seiberg-Witten strata
$\iota(M_{\fs^\pm})$ with non-trivial invariants
appear only in the top level $\sM_{\tilde\ft}$ of
$\bar\sM_{\tilde\ft}$ if and only if Seiberg-Witten
strata $\iota(M_{\fs})\subset\sM_{\ft}$ with non-trivial invariants appear
only in the top level $\sM_{\ft}$ of $\bar\sM_{\ft}$.
Since $X$ is ``effective'' by hypothesis,
we may assume Conjecture \ref{conj:Multiplicity} holds. Also,
$$
n_a(\tilde\ft) = n_a(\ft) > 0,
$$
using equation \eqref{eq:ASDCodimension} and the facts that
$p_1(\tilde\ft)=p_1(\ft)-1$ and $c_1(\tilde\ft)=c_1(\ft)$ by Lemma
\ref{lem:BlowUpSpinu} and $\sigma(\tilde X)=\sigma(X)-1$. Hence,
Corollary \ref{cor:CompactReductionFormula}
applies to $\sM_{\tilde\ft}^{*,0}$.

Theorem \ref{thm:Main} now follows by applying Proposition
\ref{prop:DegreeZeroBlowUp} in conjunction with the relation
\eqref{eq:CompactReductionFormula1} for the cobordism
$\bar\sM_{\tilde\ft}^{*,0}$. Equation \eqref{eq:CompactReductionFormula1}
(with the additional hypothesis of Corollary
\ref{cor:CompactReductionFormula}) gives
\begin{equation}
\label{eq:BlownUpCobordismFormula}
\begin{aligned}
{}&\#\left(\bar \sV(ez)\cap \bar M_{\kappa+1/4}^{w+\PD[e]}(\tilde X)\right)
\\
&=
- 2^{1-n_a}\sum_{\{\fs:\ \fs\oplus\fs'=\ft\}}
\left((-1)^{o_{\tilde\ft}(w+\PD[e],\fs^+)}
\langle\mu_{p}(ez)\smile\mu_{c}^{n_a-1},[\bL_{\tilde\ft,\fs^+}]\rangle
\right.
\\
&\qquad
\left. + (-1)^{o_{\tilde\ft}(w+\PD[e],\fs^-)}
\langle\mu_{p}(ez)\smile\mu_{c}^{n_a-1},[\bL_{\tilde\ft,\fs^-}]\rangle
\right).
\end{aligned}
\end{equation}
{}From Theorem \ref{thm:CompactReductionFormula} we see that we need to
to consider the following cases, when $n_a>0$:
\alphenumi
\begin{enumerate}
\item $\deg(z) = d_a < 2r(\Lambda)$,
\item $d_a = 2r(\Lambda)$ and $\deg(z) = d_a$,
\item $d_a = 2r(\Lambda)$ and $d_a < \deg(z) \leq d_a + 2n_a - 2$.
\end{enumerate}
{\bfseries Case (a).}
The condition $n_a>0$ is equivalent to $\delta<i(\Lambda)$, since
$n_a=\frac{1}{8}(2i(\Lambda)-d_a)$ by equation \eqref{eq:ASDCodimension}
and $\deg(z)=d_a=2\delta$ in this case.

Using the definition \eqref{eq:DefineDInvarBlowUp} of the Donaldson
invariants and using $c_1(\ft)=\Lambda$ in Theorem
\ref{thm:CompactReductionFormula} and Corollary
\ref{cor:CompactReductionFormula} yields
\begin{equation}
\label{eq:SimplifiedBlownUpCobordismFormulaIsZero}
D_X^w(z) = 0, \quad\text{for}\quad \deg(z) < 2r(\Lambda).
\end{equation}
This proves Case (a).

\noindent{\bfseries Case (b).}
The condition $n_a>0$ is again equivalent to $\delta<i(\Lambda)$, since
$\deg(z)=d_a=2\delta$ in this case. We also have $\deg(z)=2r(\Lambda)$.

Using the definition \eqref{eq:DefineDInvarBlowUp} of the Donaldson
invariants, applying our blow-up formula \eqref{eq:DegreeZeroBlowUp} (with
$k=0$) for link pairings to equation
\eqref{eq:BlownUpCobordismFormula}, and using $c_1(\ft)=\Lambda$ yields
\begin{equation}
\label{eq:SimplifiedBlownUpCobordismFormula}
\begin{aligned}
D_X^w(z)
&=
2^{1-n_a} 2^{-\delta_2-2\delta_0}(-1)^{\delta_0+\delta_1+1}
\\
&\quad\times\sum_{\{\fs:\ \fs\oplus\fs'=\ft\}} (-1)^{o_{\ft}(w,\fs)}
C_{\chi,\si}\left(\deg(z),\delta_{c},2r(\La),d_s(\fs),\delta_1\right)
\\
&\quad\times SW_{X,\fs}(x^d\vartheta)
\langle c_1(\fs)-\Lambda,h\rangle^{\delta_2}.
\end{aligned}
\end{equation}
Note that although we assume $\delta_{c}=n_a-1$ in Case (b), we allow
$\delta_{c}\leq n_a-1$ in the above sum, as 
$\delta_{c}<n_a-1$ in Case (c).  The inequality
\eqref{eq:LevelOfReducible} implies that the subset of $\fs\in\Spinc(X)$ giving
a splitting $\ft=\fs\oplus\fs'$ coincides with the subset for which
$r(\Lambda,c_1(\fs))=r(\La)$.  Hence, the sum in
\eqref{eq:SimplifiedBlownUpCobordismFormula} is over the same subset of
$\Spinc(X)$ as that in equation \eqref{eq:Main}.

We simplify the sign factor $(-1)^{o_{\ft}(w,\fs)}$ in equation
\eqref{eq:SimplifiedBlownUpCobordismFormula} by writing
\begin{equation}
\label{eq:Sign1}
\begin{aligned}
o_{\ft}(w,\fs)
&=
\textstyle{\frac{1}{4}}\left( w -\La +c_1(\fs)\right)^2
\quad\text{(from definition \eqref{eq:OrientationFactor})}
\\
&=
\textstyle{\frac{1}{2}}c_1(\fs)\cdot(w -\La)
+
\textstyle{\frac{1}{4}}\left((w -\La)^2+c_1(\fs)^2\right).
\end{aligned}
\end{equation}
Because $c_1(\fs)$ and $\La-w$ are characteristic,
we have $c_1(\fs)^2\equiv (\La-w)^2\equiv \sigma\pmod 8$.
Thus, equation \eqref{eq:Sign1} yields
\begin{equation}
\begin{aligned}
\label{eq:Sign2}
o_{\ft}(w,\fs)
&\equiv
\textstyle{\frac{1}{2}}c_1(\fs)\cdot(w-\Lambda)
+
\textstyle{\frac{1}{2}}\sigma\pmod 2
\\
&\equiv
\textstyle{\frac{1}{2}}\left(w^2 + c_1(\fs)\cdot(w-\La)\right)
+\textstyle{\frac{1}{2}}\left(\sigma - w^2\right) \pmod 2.
\end{aligned}
\end{equation}
Substituting equation \eqref{eq:Sign2} for $o_{\ft}(w,\fs)\pmod{2}$
into equation \eqref{eq:SimplifiedBlownUpCobordismFormula} implies
that the power of $(-1)$ in that formula for $D_X^w(z)$ becomes
$$
(-1)^{\delta_0+\delta_1+1}
(-1)^{\frac{1}{2}(\sigma -w^2)}
(-1)^{\frac{1}{2}(w^2 + c_1(\fs)\cdot(w-\La))},
$$
matching the power of $(-1)$ appearing in equation \eqref{eq:Main}.

Equation \eqref{eq:ASDCodimension} gives
$n_a(\ft)=\frac{1}{4}(i(\La)-\delta)$ and so the power of $2$ in equation
\eqref{eq:SimplifiedBlownUpCobordismFormula} for $D_X^w(z)$ then becomes
$$
2^{1-\frac{1}{4}(i(\La)-\delta)} 2^{-\delta_2-2\delta_0},
$$
matching the power of $2$ appearing in equation \eqref{eq:Main}.

Finally we simplify the expression for the constant $C_{\chi,\si}\left(
\deg(z),\delta_{c},2r(\La),d_s(\fs),\delta_1\right)$.  Equation
\eqref{eq:ASDCodimension}, the equality $\deg(z)+2\delta_{c}=d_a +2n_a-2$ and
the assumption that $d_a=2r(\La)$ gives
$$
\delta_{c}
=\textstyle{\frac{1}{2}}\left( d_a+2n_a-2-\deg(z)\right)
=\textstyle{\frac{1}{4}}\left(3r(\La)+i(\La)\right)
 -\textstyle{\frac{1}{2}}\deg(z) -1.
$$
(Note that this holds without the assumption $\delta_{c}=n_a-1$.)
Then, by the expression for $C_{\chi,\si}$ in Lemma
\ref{lem:JacobiPolynomialInterp},
$$
C_{\chi,\si}\left(\deg(z),\delta_{c},2r(\La),d_s(\fs),\delta_1\right)
=
H_{\chi,\si}(\La^2,\deg(z),d_s(\fs),\delta_1),
$$
where the function $H$ is defined in  equation
\eqref{eq:MainJacobiCoefficient}. This completes the proof of the formula
\eqref{eq:Main} in Case (b).

The result mentioned in Remark \ref{rmk:H1andH3} for $z=z'Y$ can be proved
by the same argument, noting that $z$ as described there is
intersection-suitable in the sense of Lemma
\ref{lem:IntersectionSuitable} and that the pairings with
$\bL_{\tilde\ft,\fs^\pm}$ all vanish by Proposition
\ref{prop:DegreeZeroBlowUp}.

\noindent{\bfseries Case (c).}
Continue to assume $d_a=2r(\Lambda)$, so the reducibles
(with non-trivial Seiberg-Witten functions) can lie in the top
level (but not in any lower level).  This case follows in exactly the same
way as Case (b), except that we now use equation
\eqref{eq:CompactReductionFormula1IsZero} in place of equation
\eqref{eq:CompactReductionFormula1} when $\deg(z)$ lies in the range
\begin{equation}
\label{eq:VanishingDNonTrivSWRelnDegRange}
d_a < \deg(z) \leq d_a + 2n_a - 2,
\end{equation}
so we obtain non-trivial relations among the Seiberg-Witten invariants from
the cobordism.  Using $n_a = \frac{1}{8}(2i(\Lambda)-d_a)$, the upper bound
in equation \eqref{eq:VanishingDNonTrivSWRelnDegRange} becomes
\begin{align*}
d_a + 2n_a - 2
&=
d_a + \textstyle{\frac{1}{4}}(2i(\Lambda)-d_a) - 2
\\
&=
\textstyle{\frac{3}{4}}d_a + \textstyle{\frac{1}{2}}i(\Lambda) - 2
\\
&=
\textstyle{\frac{3}{2}}r(\Lambda) + \textstyle{\frac{1}{2}}i(\Lambda) - 2
\\
&=
\textstyle{\frac{1}{2}}(3r(\Lambda) + i(\Lambda)) - 2.
\end{align*}
Thus, our pair of inequalities reduces to
\begin{equation}
\label{eq:BoundsOnDegz}
2r(\Lambda) < \deg(z) \leq r(\Lambda) +
\textstyle{\frac{1}{2}}(r(\Lambda) + i(\Lambda)) - 2.
\end{equation}
Therefore we obtain a non-trivial relation amongst the Seiberg-Witten
invariants and a vanishing result for Donaldson invariants when the
constraint \eqref{eq:BoundsOnDegz} on $\Lambda^2$ and $\deg(z)$ holds.
\end{proof}

\begin{proof}[Proof of Theorem \ref{thm:FLthm}]
  By hypothesis, $\La\cdot c_1(\fs)=0$ for all $\fs\in\Spinc(X)$ with
  $SW_X(\fs)\neq 0$, so from equation \eqref{eq:SWInTopLevelFunction} for
  $r(\Lambda,c_1(\fs))$ we have
\begin{align*}
r(\Lambda,c_1(\fs))
&=
-c_1(\fs)^2 - \Lambda^2 - \textstyle{\frac{3}{4}}(\chi+\sigma)
\\
&=
-(2\chi + 3\sigma) - \Lambda^2 - \textstyle{\frac{3}{4}}(\chi+\sigma)
\\
&=
-(\chi + \sigma) - \Lambda^2 + c(X)
\\
&=
r(\Lambda),
\end{align*}
using the definition of $c(X)$ (see \S \ref{subsec:Statement}) and the
definition \eqref{eq:SWInTopLevelFunction} of $r(\Lambda)$, and the formula
(\S \ref{subsec:Statement}) for $c_1(\fs)^2$ when $X$ has SW-simple type. A
reducible $\PU(2)$ monopole in $\bar\sM_{\ft}$ defined by a splitting
$\ft_{\ell}=\fs\oplus\fs'$ appears in level
$$
\ell(\ft,\fs)
=
\textstyle{\frac{1}{8}}(d_a(\ft) - 2r(\Lambda,c_1(\fs)))
=
\textstyle{\frac{1}{8}}(d_a(\ft) - 2r(\Lambda)),
$$
and thus all reducibles appear in the same level of $\bar\sM_{\ft}$. Hence,
the sum over $\fs\in\Spinc(X)$ with $r(\Lambda,c_1(\fs))=r(\Lambda)$ can
be replaced 
by a sum over $\fs\in\Spinc(X)$ when $\Lambda\in B^\perp$. We write
$\deg(z)=2\delta$, as in the hypothesis of the theorem.

{\bfseries Case (a).} In this situation, $\delta<r(\Lambda)$,
$\delta<i(\Lambda)$, and 
$$
D_X^w(z) = 0, \quad\text{when }0\leq \delta < r(\Lambda),
$$
by Theorem \ref{thm:Main}.

{\bfseries Case (b).} In this situation, $\delta=r(\Lambda)$
and $\delta < i(\Lambda)$. We can further simplify the formula
\eqref{eq:Main}.  First, using $i(\La)=2c(X)-r(\La)=2c(X)-\delta$,
$\delta_2=\delta-2m$, and $\delta_0=m$, the power of $2$ in equation
\eqref{eq:Main} becomes
$$
2^{1-\frac{1}{4}(2c(X)-2\delta) -\delta}
=
2^{1- \frac{1}{2}(c(X)+\delta)},
$$
matching the power of $2$ in equation \eqref{eq:DSWrel}.  The power of
$(-1)$ in equation \eqref{eq:Main} is
$$
(-1)^{m+1}
(-1)^{\frac{1}{2}(\sigma -w^2)}
(-1)^{\frac{1}{2}(w^2 + c_1(\fs)\cdot w)},
$$
since $\delta_0=m$ and $c_1(\fs)\cdot\Lambda=0$, and also matches that
in \eqref{eq:DSWrel}.  Finally, $d_s(\fs)=0$ because we assume $X$ is
SW-simple type, so the constant $H_{\chi,\si}(\cdot)$
is equal to one and thus equation \eqref{eq:DSWrel} follows from
equation \eqref{eq:Main}.

{\bfseries Case (c).} Using $\deg(z) = 2\delta$, $r(\Lambda)+i(\Lambda) =
2c(X)$, and equation \eqref{eq:SWSimpleUDepthDiracIndexParam} for
$r(\Lambda)$, the constraint \eqref{eq:BoundsOnDegz} simplifies to
\begin{equation}
\label{eq:ConstraintsOnLambdaDelta}
2r(\Lambda) < 2\delta \leq r(\Lambda) + c(X) - 2.
\end{equation}
The vanishing relation follows from Case (c) in Theorem
\ref{thm:Main}. This completes the proof.
\end{proof}

\begin{proof}[Proof of Corollary \ref{cor:FLthm}]
We consider the last case of Theorem \ref{thm:FLthm}, where $\delta$ and
$\Lambda^2$ obey the constraints \eqref{eq:ConstraintsOnLambdaDelta} and so
$$
r(\Lambda) < c(X) - 2.
$$
Therefore, the choice of $r(\Lambda)<c(X)-2$ giving the largest possible
integer $\delta$ (admitting a non-trivial vanishing relation) is $r(\Lambda) =
c(X)-4$, achieved when $\Lambda=\Lambda_0$ with $\Lambda_0^2 =
4-(\chi+\sigma)$. By hypothesis, such a class $\Lambda_0\in B^\perp$ exists.
Therefore, the pair
\eqref{eq:ConstraintsOnLambdaDelta} of inequalities constrains
$$
\delta = c(X) - 3.
$$
Thus, using $z=x^m h^{\delta-2m}$ with $0\leq m \leq [\delta/2]$, we obtain
for $w_0\in H^2(X;\ZZ)$ with $w_0-\Lambda_0\equiv w_2(X) \pmod{2}$
\begin{equation}
\label{eq:PrelimVanishingSWRelation}
\sum_{\fs\in\Spinc(X)}(-1)^{\frac{1}{2}(w_0^2+c_1(\fs)\cdot w_0)}
\SW_X(\fs)\langle c_1(\fs)-\Lambda_0,h\rangle^{d}
=
0, \quad 0\leq d\leq c(X)-3.
\end{equation}
Indeed, if $c(X)-3$ is even, then we may choose $m=\delta/2$ to obtain the
above relation with $d=0$ and, as we explain shortly, the relation for
$d=0$ also holds when $c(X)-3$ is odd. Hence, the degree-$d$ terms in the
Taylor expansion of $\bS\bW_X^{w_0}(h)e^{-\langle\Lambda,h\rangle}$ about
$h=0$ are zero for $0\leq d\leq c(X)-3$ and so the same holds for
$\bS\bW_X^{w_0}(h)$.

If $w$ is any integral lift of $w_2(X)$, write $w=w+\Lambda_0-\Lambda_0$
and observe that
$$
\bS\bW_X^{w}(h)
= (-1)^{\frac{1}{2}(\Lambda_0^2-2w\cdot\Lambda_0)}\bS\bW_X^{w+\Lambda_0}(h).
$$
Thus $\bS\bW_X^{w}(h)$ vanishes to the same order as $\bS\bW_X^{w_0}(h)$
with $w_0=w+\Lambda_0$ and this completes the proof, aside from the remark
below on the case of odd $c(X)-3$.

When $c(X)-3$ is odd, it only remains to show that the relation
\eqref{eq:PrelimVanishingSWRelation} still holds when $d=0$. We choose
$\Lambda_1\in B^\perp$ with $\Lambda_1^2=6-(\chi+\sigma)$ and
$r(\Lambda_1)=c(X)-6$, so that
\eqref{eq:ConstraintsOnLambdaDelta} allows $\delta_1=c(X)-4$, which must be
even and thus, taking $m=\delta_1/2$ yields
$$
\sum_{\fs\in\Spinc(X)}(-1)^{\frac{1}{2}(w_1^2+c_1(\fs)\cdot w_1)}
\SW_X(\fs)
=
0,
$$
for any $w_1\in H_2(X;\ZZ)$ with $w_1-\Lambda_1\equiv w_2(X)$. Writing
$w_1 = w_0 + \Lambda_1-\Lambda_0$ and combining
$$
\bS\bW_X^{w_1}(h)
= (-1)^{\frac{1}{2}((\Lambda_1-\Lambda_0)^2+2w_0\cdot(\Lambda_1-\Lambda_0))}
\bS\bW_X^{w_0}(h).
$$
with the previous vanishing result yields the relation
\eqref{eq:PrelimVanishingSWRelation} when $d=0$.
\end{proof}

\begin{proof}[Proof of Theorem \ref{thm:DSWSeriesRelation}]
We assume without loss that $c(X)>0$.  {}From Theorem \ref{thm:FLthm} we
know that $D_X^w(x^mh^{\delta-2m})=0$ if $\delta < r(\Lambda)$ and that the
first potentially non-zero invariant is given by equation \eqref{eq:DSWrel}
when $\delta=r(\Lambda)$. The cobordism method constrains $\Lambda^2$ by
requiring that $\delta < i(\Lambda)$. Hence, from the graphs of
$r(\Lambda)$ and $i(\Lambda)$ as functions of $\Lambda^2$ (see Figure 1 in
\cite{FKLM}) one sees that these two lines meet for $\Lambda_0\in
H^2(X;\ZZ)$ with $\Lambda_0^2=-(\chi+\sigma)$, at which point
$r(\Lambda_0)=c(X)=i(\Lambda_0)$. Therefore, we choose $\Lambda^2$ to give
the largest possible $\delta = r(\Lambda) < c(X)$. We also take $\Lambda\in
B^\perp$, to simplify the formula \eqref{eq:DSWrel} and as the SW-basic classes
$c_1(\fs)$ are characteristic, this gives (for $B$ non-empty) $\Lambda\cdot
c_1(\fs)=0$ and $\Lambda\cdot c_1(\fs)\equiv
\Lambda^2\pmod{2}$, so that $\Lambda^2$ is even. Thus we want to choose
$\Lambda\in B^\perp$ with smallest even value of $\Lambda^2 >
-(\chi+\sigma)$, namely
\begin{equation}
  \label{eq:NiceSquareLambda}
\Lambda^2 = 2 - (\chi + \sigma),
\end{equation}
because $\chi+\sigma$ is even (in fact, divisible by four since $b_1(X)=0$
and $b_2^+(X)$ is odd). By hypothesis, $\Lambda$ exists and
the formula \eqref{eq:SWSimpleUDepthDiracIndexParam} for $r(\Lambda)$ and
the definition of $c(X)$ yield
\begin{align*}
\delta
&=r(\Lambda)
=
-\Lambda^2 + c(X) - (\chi+\sigma)
\\
&=
c(X) - 2.
\end{align*}
Therefore Theorem \ref{thm:FLthm} and the fact that
$\delta=c(X)-2=r(\Lambda)$ yield
\begin{equation}
\label{eq:VanishingDInvariants}
D_X^w(x^mh^{d-2m}) = 0,
\quad 0\leq d < \delta\quad\text{and}\quad 0\leq m\leq [d/2].
\end{equation}
When $m=0$, the power of $(-1)$ in equation \eqref{eq:DSWrel} simplifies to
$$
(-1)^{\frac{1}{2}(w^2 + c_1(\fs)\cdot w)},
$$
as $\frac{1}{2}(\sigma-w^2)\equiv 1\pmod{2}$. Indeed, to see this note that
$w-\Lambda$ is characteristic, so $(w-\Lambda)\cdot\Lambda
\equiv\Lambda^2\pmod{2}$, and as $\chi+\sigma\equiv 0\pmod{4}$, we have
\begin{equation}
\label{eq:wSquaredConstraint}
\begin{aligned}
w^2
&=
(w-\Lambda)^2 + 2(w-\Lambda)\cdot\Lambda + \Lambda^2
\\
&\equiv
\sigma + \Lambda^2 \pmod{4}
\\
&\equiv
\sigma + 2 \pmod{4} \quad\text{(by Equation \eqref{eq:NiceSquareLambda})}.
\end{aligned}
\end{equation}
The power of $2$ in equation \eqref{eq:DSWrel}, when $\delta=c(X)-2$,
becomes
$$
2^{2- c(X)},
$$
as we expect from Witten's formula \eqref{eq:WittenConjSeries}.

{}From equation \eqref{eq:VanishingDInvariants}, the invariants
$D_X^w(h^d)$ and $D_X^w(xh^{d-2})$ are zero when $d<\delta=c(X)-2$ (while
the method of this article does not allow us to compute the invariants
when $d\geq \delta+4$), so (compare equation
\eqref{eq:DSeries})
\begin{equation}
  \label{eq:DSeriesLowOrder}
\begin{aligned}
\bD^{w}_{X}(h)
&\equiv
0
\pmod{h^{\delta}},
\\
\bD^{w}_{X}(h)
&\equiv
\textstyle{\frac{1}{\delta!}}D_X^w(h^\delta)
+
\textstyle{\frac{1}{2(\delta-2)!}}D_X^w(xh^{\delta-2})
\pmod{h^{\delta+2}}.
\end{aligned}
\end{equation}
For a monomial $z\in\AAA(X)$, the invariant $D_X^w(z)$ is zero unless
$\deg(z)\equiv -2w^2-\frac{3}{2}(\chi+\sigma)\pmod{8}$. Therefore, as
$\delta\equiv w^2-\frac{3}{4}(\chi+\sigma)\pmod{4}$, the invariants
$D_X^w(h^{\delta+2})$ and $D_X^w(xh^{\delta})$ are necessarily zero and the
next potentially non-zero invariant of higher degree in $h$ would be
$D_X^w(xh^{\delta+2})$.

For the terms $D_X^w(h^\delta)$, equation \eqref{eq:DSWrel} yields
\begin{equation}
  \label{eq:DLambdaPlus}
{\textstyle{\frac{1}{\delta!}}}D_X^w(h^\delta)
=
2^{2- c(X)}
\sum_{\fs\in\Spinc(X)}
(-1)^{\frac{1}{2}(w^2 + c_1(\fs)\cdot w)}
SW_X(\fs)\textstyle{\frac{1}{\delta!}}
\langle c_1(\fs)-\Lambda,h\rangle^{\delta}.
\end{equation}
Since $\bS\bW_X^{w-\Lambda}(h)\equiv 0\pmod{h^{\delta}}$ by Corollary
\ref{cor:FLthm}, we have
\begin{equation}
\label{eq:SWSeriesVanishing}
\bS\bW_X^w(h)
=
(-1)^{\frac{1}{2}(\Lambda^2 + 2(w-\Lambda)\cdot\Lambda)}\bS\bW_X^{w-\Lambda}(h)
\equiv
0\pmod{h^\delta}.
\end{equation}
Therefore, using equation \eqref{eq:DLambdaPlus} and noting
that the terms in $e^{\frac{1}{2}Q(h,h)}$ and
$e^{\langle-\Lambda,h\rangle}$ of lowest degree in $h$ are both $1$ and the
lowest-degree non-zero term in $\bS\bW_X^{w}(h)$ has degree $\delta$ in $h$
by equation \eqref{eq:SWSeriesVanishing}, we see that
\begin{equation}
\label{eq:IWantMyDeltaD}
\textstyle{\frac{1}{\delta!}}D_X^w(h^\delta)
=
\left[2^{2- c(X)}\bS\bW_X^w(h)e^{\langle-\Lambda,h\rangle}\right]_{\delta}
=
\left[2^{2- c(X)}e^{\frac{1}{2}Q(h,h)}\bS\bW_X^w(h)\right]_{\delta},
\end{equation}
where $[\,\cdot\,]_\delta$ denotes the term of degree $\delta$ in $h$ in
the power series.

For the term $D_X^w(xh^{\delta-2})$, equation \eqref{eq:DSWrel} yields
\begin{equation}
  \label{eq:DxLambdaPlus}
\begin{aligned}
D_X^w(xh^{\delta-2})
&=
-2^{2- c(X)}
\sum_{\fs\in\Spinc(X)}
(-1)^{\frac{1}{2}(w^2 + c_1(\fs)\cdot w)}
SW_X(\fs)\langle c_1(\fs)-\Lambda,h\rangle^{\delta-2}
\\
&=
-\left[2^{2- c(X)}\bS\bW_X^w(h)e^{\langle-\Lambda,h\rangle}\right]_{\delta-2}
\cdot(\delta-2)!
\\
&=0,
\end{aligned}
\end{equation}
where the final equality follows from the fact that the term in
$e^{\langle-\Lambda,h\rangle}$ of lowest degree in $h$ is $1$ and the terms
in $\bS\bW_X^{w}(h)$ of degree $\delta-2$ or lower in $h$ are zero.
Combining equations \eqref{eq:DSeriesLowOrder}, \eqref{eq:IWantMyDeltaD}, and
\eqref{eq:DxLambdaPlus} thus completes the proof.
\end{proof}

\ifx\undefined\bysame
\newcommand{\bysame}{\leavevmode\hbox to3em{\hrulefill}\,}
\fi

\end{document}